\documentclass[11pt]{amsart}
\usepackage{amssymb}
\usepackage{amsmath}
\usepackage{amsrefs}
\usepackage[colorlinks=true]{hyperref}
\usepackage{cleveref}
\usepackage{amscd}
\usepackage{amsthm}
\usepackage{mathtools}
\usepackage{tikz-cd}
\usepackage{bm}
\usepackage{mathrsfs}
\usepackage{enumerate}
\usepackage[margin=2.93cm]{geometry}

\usepackage{xcolor}
\hypersetup{
    colorlinks,
    linkcolor={red!50!black},
    citecolor={green!50!black},
    urlcolor={blue!80!black}
}
\allowdisplaybreaks

\newtheorem{theorem}{Theorem}[section]
\newtheorem{lemma}[theorem]{Lemma}
\newtheorem{proposition}[theorem]{Proposition}
\newtheorem{corollary}[theorem]{Corollary}

\newtheorem{conjecture}[theorem]{Conjecture}
\newtheorem{question}[theorem]{Question}
\newtheorem{fact}[theorem]{Fact}

\newtheorem*{assumption-no-number}{Assumption}
\newtheorem*{corollary*}{Corollary}

\theoremstyle{definition}
\newtheorem{definition}[theorem]{Definition}
\newtheorem{example}[theorem]{Example}

\theoremstyle{remark}
\newtheorem{remark}[theorem]{Remark}

\numberwithin{equation}{section}

\newcommand{\R}{\mathbb{R}}

\newcommand{\N}{\mathbb{N}}
\newcommand{\Z}{\mathbb{Z}}

\renewcommand{\L}{\mathcal{L}}

\newcommand{\lt}{\left}
\newcommand{\rt}{\right}
\newcommand{\tms}{\times}

\newcommand{\rmk}{\begin{remark}}
\newcommand{\ermk}{\end{remark}}
\newcommand{\cor}{\begin{corollary}}
\newcommand{\ecor}{\end{corollary}}
\newcommand{\eq}{\begin{equation}}
\newcommand{\eeq}{\end{equation}}
\newcommand{\eqs}{\begin{equation*}}
\newcommand{\eeqs}{\end{equation*}}
\newcommand{\prop}{\begin{proposition}}
\newcommand{\eprop}{\end{proposition}}
\newcommand{\thm}{\begin{theorem}}
\newcommand{\ethm}{\end{theorem}}
\newcommand{\conj}{\begin{conjecture}}
\newcommand{\econj}{\end{conjecture}}
\newcommand{\lem}{\begin{lemma}}
\newcommand{\elem}{\end{lemma}}
\newcommand{\defi}{\begin{definition}}
\newcommand{\edefi}{\end{definition}}
\newcommand{\ex}{\begin{example}}
\newcommand{\eex}{\end{example}}
\newcommand{\alis}{\begin{align*}}
\newcommand{\ealis}{\end{align*}}
\newcommand{\pf}{\begin{proof}}
\newcommand{\epf}{\end{proof}}
\newcommand{\ali}{\begin{align}}
\newcommand{\eali}{\end{align}}
\newcommand{\qus}{\begin{question}}
\newcommand{\equs}{\end{question}}

\newcommand{\mc}{\mathcal}
\renewcommand{\bf}{\textbf}

\newcommand{\C}{\mathbb{C}}

\newcommand{\sub}{\subseteq}

\newcommand{\ov}{\overline}

\newcommand{\bb}{\mathbb}

\newcommand{\op}{\operatorname}

\renewcommand{\a}{\alpha}
\renewcommand{\b}{\beta}
\renewcommand{\d}{\partial}
\newcommand{\e}{\epsilon}

\newcommand{\g}{\gamma}

\newcommand{\s}{\sigma}
\renewcommand{\t}{\theta}

\renewcommand{\l}{\lambda}
\renewcommand{\o}{\omega}

\newcommand{\fk}{\frak}

\newcommand{\G}{\Gamma}
\renewcommand{\L}{\Lambda}
\renewcommand{\O}{\Omega}

\renewcommand{\S}{\Sigma}

\newcommand{\fT}{\mathfrak{T}}

\newcommand{\Q}{\mathbb{Q}}

\newcommand{\bu}{\bullet}

\renewcommand{\ov}{\overline}

\newcommand{\Hom}{\operatorname{Hom}}

\newcommand{\kos}[2]{{#1/\!\!/#2}}
\newcommand{\bsr}{\boldsymbol{r}}

\newcommand{\mush}{\mu sh}

\DeclareMathOperator{\pd}{prod}

\newcommand{\rdim}[1]{\operatorname{Rdim}#1}
\newcommand{\sd}{{sd}}

\title[The Rouquier dimension and wrapped Fukaya categories]{On the Rouquier dimension of wrapped Fukaya categories and a conjecture of Orlov} 

\author[Shaoyun Bai]{Shaoyun Bai}
\address {Department of Mathematics, Princeton University, New Jersey 08544, USA}
\email {shaoyunb@math.princeton.edu}

\author[Laurent C\^ot\'e]{Laurent C\^ot\'e}
\address{Department of Mathematics, Harvard University, 1 Oxford St, Cambridge, MA 02138}
\email{lcote@math.harvard.edu}


\begin{document}

\begin{abstract}

We study the Rouquier dimension of wrapped Fukaya categories of Liouville manifolds and pairs, and apply this invariant to various problems in algebraic and symplectic geometry.

On the algebro-geometric side, we introduce a new method based on symplectic flexibility and mirror symmetry to bound the Rouquier dimension of derived categories of coherent sheaves on certain complex algebraic varieties and stacks. These bounds are sharp in dimension at most $3$. As an application, we resolve a well-known conjecture of Orlov for new classes of examples (e.g.\ toric $3$-folds, certain log Calabi--Yau surfaces). We also discuss applications to non-commutative motives on partially wrapped Fukaya categories.

On the symplectic side, we study various quantitative questions such as: (1) given a Weinstein manifold, what is the minimal number of intersection points between the skeleton and its image under a generic compactly-supported Hamiltonian diffeomorphism? (2) what is the minimal number of critical points of a Lefschetz fibration on a Liouville manifold with Weinstein fibers? We give lower bounds for these quantities which are to our knowledge the first to go beyond the basic flexible/rigid dichotomy.  
\end{abstract}
\maketitle

\setcounter{tocdepth}{1}


\section{Introduction}
\subsection{Overview}\label{subsection:intro-motivation}
In his seminal work \cite{rouquier}, Rouquier introduced a notion of dimension for triangulated categories. Here is a sketch of the definition: if $\mc{T}$ is a triangulated category and $G$ is a (split-)generator, then $G$ is said to generate $\mc{T}$ in time $t$ if every object of $\mc{T}$ can be built from $G$ by taking finite direct sums, shifts, summands, and \emph{at most $t$} cones. The Rouquier dimension of $\mc{T}$ is by definition the minimal generation time over all generators of $\mc{T}$. (If $\mc{T}=0$, then we define the dimension to be $-1$). 

In symplectic topology, an important class of triangulated categories are derived Fukaya categories in their various guises. We will be mainly interested in wrapped Fukaya categories of Liouville manifolds. More precisely, given a Liouville manifold $(X, \l)$, its wrapped Fukaya category $\mc{W}(X)$ is an $A_\infty$ category and as such, it canonically determines a triangulated category $H^0(\op{Perf} \mc{W}(X))$. The Rouquier dimension of this triangulated category is a numerical invariant of $(X, \l)$ which will be the focus of this paper.

Our applications will be twofold: first, to a problem in algebraic geometry; second, to quantitative questions in symplectic topology. We discuss each of these two classes of applications in turn in the remainder of this introduction.

\subsection{An application of symplectic topology to algebraic geometry}\label{subsection:motivation-ag}
\subsubsection{Motivation}
In classical algebraic geometry, one is typically interested in spaces (varieties, schemes, stacks, etc.) which are locally modeled on the spectrum of a commutative ring. According to the philosophy of non-commutative geometry, it is often profitable to regard $A_\infty$ categories as ``non-commutative" spaces. One can then attempt to generalize classical algebraic geometry to the setting of $A_\infty$ category theory.

For example, there are standard notions smoothness, properness, dimension, etc.\ in classical algebraic geometry. It turns out that there also exist non-commutative analogs of some of these notions. Namely: an $A_\infty$ category is defined to be smooth if its diagonal bimodule is perfect; it is defined to be proper if its diagonal bimodule is proper.\footnote{Fix an $A_\infty$ category defined over $k$. A bimodule is \emph{perfect} if is split-generated by tensor products of representable bimodules $\hom(-, K) \otimes_k \hom(L, -)$. A bimodule is \emph{proper} if it takes value in $\op{Perf} k \subset \op{Mod} k$ (the subcategory of perfect $k$-linear chain complexes).} The dimension of an $A_\infty$ category is simply its Rouquier dimension.\footnote{Note that there are other reasonable notions of dimension discussed e.g.\ in \cite{elagin-lunts}, although the Rouquier dimension seems the most established in the literature.} 

If $Y$ is an ordinary algebraic variety, then one can associate to it the $A_\infty$ category $D^b\op{Coh}(Y)$. One typically regards $D^b\op{Coh}(Y)$ as the non-commutative partner of the classical object $Y$.  It is then natural to ask: do the non-commutative notions of smoothness, properness, dimension coincide with their classical counterparts?

As one might hope, it is known that $D^b\op{Coh}(Y)$ is smooth and proper if and only if $Y$ is smooth and proper in the classical sense. Concerning dimension, there is the following well-known conjecture or Orlov: 

\conj[Orlov; see Conjecture $10$ in \cite{orlov2009remarks}]\label{conj_orlov_intro}
Let $Y$ be a smooth quasi-projective scheme of dimension $n$. Then 
\eq\label{eqn:orlov} \rdim D^{b} \op{Coh}(Y) = n. \eeq
\econj

The assumption that $Y$ is smooth is necessary; see \Cref{example:hara}. Orlov's conjecture is in general wide open. It has been only established for the following rather restrictive classes of examples, thanks to the combined work of many authors (to the best of our knowledge, this list is exhaustive). 

\begin{enumerate}
\item Smooth affine varieties, projective spaces and smooth projective quadrics \cite[Prop.\ 7.18, Ex.\ 7.7, 7.8]{rouquier};
\item All smooth projective curves \cite[Thm.\ 6]{orlov2009remarks};
\item The following schemes defined over an algebraically-closed field of characteristic $0$: del Pezzo surfaces with $\op{rk} \op{Pic}(Y) \leq 7$, Fano threefolds of type $V_5$ and $V_{22}$, toric surfaces with nef anti-canonical divisor, toric Deligne--Mumford stacks defined over $\mathbb{C}$ of dimension no more than $2$ or Picard number no more than $2$, Hirzebruch surfaces \cite{ballard2012hochschild};
\item The product of varieties from the previous item \cite{yang2016note};
\item The product of two Fermat elliptic curves with the Fermat $K3$ surface \cite[Thm.\ 1.6]{ballard2014category};
\item Smooth toric Fano varieties of dimension $3$ or $4$ \cite[Prop.\ 2.15, 2.16]{ballard2019toric};
\item Certain blow-ups of the projective space \cite[Thm.\ 4.1]{pirozhkov2019rouquier};
\item Some weak del Pezzo surfaces \cite[Cor.\ 7.5]{elagin-xu-zhang};
\end{enumerate}

We remark that Rouquier proved that, for smooth $Y$,
\eq\label{equation:rouquier-bound} n \leq \rdim D^{b} \op{Coh}(Y) \leq 2n, \eeq
so the content of Orlov's conjecture lies in improving the upper bound. See \Cref{subsection:first-examples} for further discussion and references.

\subsubsection{Results} We introduce new methods from the flexible side of symplectic topology to study Orlov's conjecture. More precisely, we prove:
\thm\label{thm:orlov_gen_state_intro}
Suppose that $Y$ is a variety over $\C$ of complex dimension $n$ which admits a homological mirror given by a polarizable Weinstein pair. Then $\rdim D^b\op{Coh}(Y)=n$ if $n \leq 3$. For $n \geq 4$, we have $\rdim D^b\op{Coh}(Y)\leq 2n -3$.
\ethm
(A Weinstein pair $(X, A)$ is said to be polarizable if $X$ admits a global Lagrangian plane field $\xi$.)
\Cref{thm:orlov_gen_state_intro} implies the following new cases of Orlov's conjecture.
\cor
\Cref{conj_orlov_intro} is true for all $3$-dimensional toric varieties over $\mathbb{C}$ and the log Calabi--Yau surfaces in \Cref{ex_logCY}. 
\ecor
In fact, \Cref{thm:orlov_gen_state_intro} shows that \eqref{eqn:orlov} holds for a large class of (possibly singular) schemes and stacks over $\mathbb{C}$, see \Cref{subsection:orlov-conj} for the list. It also gives a new proof of some of the known cases of Orlov's conjecture listed above. Finally, \Cref{thm:orlov_gen_state_intro} improves Rouquier's upper bound \eqref{equation:rouquier-bound} for $\rdim D^b\op{Coh}(-)$ for many examples in dimensions $n \geq 4$ (including all smooth toric varieties and DM stacks). 

\subsubsection{Methods: symplectic flexibility} \label{subsubsection:arborealization-intro}
Although Orlov's conjecture is a problem in algebraic geometry, our proof of \Cref{thm:orlov_gen_state_intro} is entirely symplectic. As mentioned in \eqref{equation:rouquier-bound}, it is already known that $\rdim D^b\op{Coh}(Y) \geq n$. Hence, we only need to prove the claimed upper bounds. To do this, we apply homological mirror symmetry. This transforms \Cref{thm:orlov_gen_state_intro} into a statement about the Rouquier dimension of wrapped Fukaya categories of Weinstein pairs; namely, it is enough to prove that the Rouquier dimension of a (polarizable) Weinstein pair of real dimension $2n$ is bounded above by $n$ when $n \leq 3$ and by $2n-3$ when $n \geq 4$.  This is now purely a problem in symplectic topology which can be approached using symplectic tools.

According to the recent work of Ganatra--Pardon--Shende \cite{gps2}, wrapped Fukaya categories of Liouville manifolds satisfy ``Weinstein sectorial descent". Essentially, this means that if we (suitably) cover a Liouville manifold $X$ by Weinstein sectors $\{ X_\a \}_{\a \in I}$, then $\mc{W}(X)$ can be expressed as the homotopy colimit of the wrapped Fukaya categories of the form $\mc{W}(\cap_{\a \in J} X_\a)$ (where $J$ ranges over all non-empty subsets of $I$ and the homotopy colimit is indexed by the obvious diagram of inclusions). 

It is not hard to show that the Rouquier dimension of a homotopy colimit can be bounded from above by the Rouquier dimension of its individual pieces.  This bound depends on (a) the ``depth" of the diagram (i.e.\ the longest chain of arrows) and (b) the Rouquier dimension of the individual pieces; see \Cref{lem_rdim_upper} for a precise statement. This suggests a possible strategy for upper bounding the Rouquier dimension of a Weinstein manifold: construct a sectorial cover of ``small" depth, such that each piece in the cover has ``small" Rouquier dimension. 

Let us first implement this strategy in the special case where $X=T^*M$, for $M$ a closed $n$-manifold. We construct a sectorial cover as follows. First, triangulate $M$. Then set $X_\a:= T^*\op{star}(v_\a)$ where $v_\a$ ranges over the vertices of the triangulation. Observe that the intersection of any collection of stars in a triangulation is either empty, or is the star of some higher dimensional cell. From this, it follows that the cover has depth $n+1$. It also follows (after appropriate smoothing) that $\bigcap_{\a_1,\dots, \a_n} T^*\op{star}(v_{\a_i}) \simeq T^*D^n$, and hence the Rouquier dimension of the overlaps is zero. Plugging these numbers into \Cref{lem_rdim_upper}, one finds that $\rdim \mc{W}(T^*M) \leq n$. This upper bound is sharp in general.\footnote{For example, it is a basic fact in mirror symmetry that $\mc{W}(T^*\mathbb{T}^n) \simeq D^b\op{Coh}((\mathbb{C}^*)^n)$, so it follows from Rouquier's aforementioned result for affine varieties that $\rdim \mc{W}(T^*\mathbb{T}^n)=n$.} 

We would like to extend this argument to an arbitrary Weinstein manifold $X$. Unfortunately, the naive generalization \emph{does not work}. If $X$ is an arbitrary Weinstein manifold, its skeleton may be a highly singular object. Even assuming the skeleton can be triangulated (which is not automatic), it is completely unclear how to control the Rouquier dimension of the resulting sectorial cover. It is also unclear how to convert a cover of the skeleton into a sectorial cover of its thickening (in the cotangent bundle case, there is a canonical Morse--Bott form which allows one to easily perform this conversion). 

To overcome these difficulties, we will take advantage of recent advances from the flexible side of symplectic geometry.  Namely, the arborealization theorem of \'{A}lvarez-Gavela--Eliashberg--Nadler \cites{AGEN1, AGEN2, AGEN3} provides an h-principle for simplifying the singularities of the skeleton. The upshot is that, under rather mild topological assumptions, one can homotope a Weinstein manifold so that the singularities of its skeleton are of ``arboreal" type. Such homotopies do not affect the wrapped Fukaya category.

Arboreal singularities are indexed by finite rooted trees. Morally, if $\mc{S}_\s$ is a Liouville sector obtained by thickening an arboreal singularity $\s$, one expects that $\mc{W}(\mc{S}_\s)$ is Morita equivalent to $\op{Rep}(T_\s)$, the derived representation category of the corresponding tree $T_\s$. Moreover, it is known that $\rdim \op{Rep}(T_\s)=0$ if $T_\s$ is of Dynkin type (i.e.\ the underlying graph is a Dynkin diagram) and $\rdim \op{Rep}(T_\s)=1$ otherwise. This suggests that sectorial thickenings of arboreal singularities are a good substitute for $T^*D^n$, 

Arboreal skeleta admit a canonical stratification by ``singularity-order". We prove in \Cref{subsec_arbo} that this stratification is Whitney. By a well-known theorem of Goresky, any Whitney stratification can be refined to a (possibly non-Whitney!) triangulation. This suggests the following updated strategy:  given a (polarizable) Weinstein manifold $X$ of dimension $2n$, apply a homotopy to make the skeleton arboreal. It therefore admits a triangulation. Now construct a sectorial cover by taking arboreal thickenings of the stars of the vertices, like for cotangent bundles.

Unfortunately, there are significant obstructions to making this strategy rigorous. To begin with, one would need to develop a good theory of arboreal sectors (which should be closed under intersections), as well as a procedure for constructing a ``sectorial thickening" of an arboreal singularity.  Moreover, one would need a method for construct arboreal sectors from stars of (possibly non-Whitney) triangulations. Due to this difficulties, it is unclear to us whether the above strategy can feasibly be implemented in the framework of Fukaya categories.

Instead, we work in the framework of microlocal sheaf theory. The details of this approach are carried out in \Cref{sec_arbo}. After reviewing some basic properties of microlocal sheaves, we study the stratified topology of arboreal skeleta, and use this to construct a suitable triangulation. We then carry through a microlocal sheaf analog of the above covering argument to obtain an upper bound on the Rouquier dimension of the category of microlocal sheaves. Finally, we apply the main result of \cite{gps3} to pass back to wrapped Fukaya categories. 

The precise upper bound which we obtain via this argument (and under the assumption that $X$ is stably polarizable) is $\rdim \mc{W}(X^{2n}) \leq n$ for $n \leq 3$, and $\rdim \mc{W}(X^{2n}) \leq 2n-3$ for $n>3$. We know no counterexample to the inequality $\rdim \mc{W}(X^{2n}) \leq n$ for arbitrary Weinstein manifolds.

\subsubsection{A digression on motives}
As a byproduct of our method, we derive some consequences for the structure of \emph{non-commutative motives} on partially wrapped Fukaya categories. The notion of a non-commutative motive is was introduced by Tabuada \cite{tabuada-motives} and is briefly summarized in \Cref{subsection:motives}. To give an illustration, suppose that $(X,A)$ is a polarizable Weinstein pair such that $\mc{W}(X,A)$ is proper (it is always smooth according to \cite[Cor.\ 1.19]{gps2}). Then we find that:
\begin{itemize}
\item $HH_\bu(\mc{W}(X,A))$ is concentrated in degree zero; 
\item $ \op{Perf} \mc{W}(X,A)$ satisfies the \emph{non-commutative Weil conjecture} introduced by Tabuada (following proposals of Kontsevich) in \cite{tabuada-weil}.
\end{itemize}
Both of these statements follow from a more general result about non-commutative motives, which is itself an easy consequence of the arborealization method described above. We refer the reader to \Cref{subsection:motives} for details.

\subsection{Quantitative problems in symplectic topology}\label{subsection:motivation-symplectic}
\subsubsection{Motivation}
As motivation for the geometrically-inclined symplectic topologist, we state some quantitative questions of rather classical flavor. 
\begin{question}\label{question:arnold-type}
Let $(X, \l)$ be a Weinstein\footnote{Recall that a Liouville manifold is said to be \emph{Weinstein} if the Liouville vector field is gradient-like with respect to a proper Morse function.} manifold with skeleton $\fk{c}_X$. If $\phi: X \to X$ is a compactly-supported Hamiltonian diffeomorphism, what is the minimal possible number of intersection points of $\phi (\fk{c}_X)$ with $\fk{c}_X$? 
\end{question}

If $(X, \l)=(T^*M, pdq)$ for a closed manifold $M$, then $\fk{c}_X= 0_M$ and \Cref{question:arnold-type} is classical. As is now well-known, it follows from the existence of Floer cohomology (or finite dimensional analogs) that $|\fk{c}_X \cap \phi (\fk{c}_X)|$ is bounded from below by the sum of the Betti numbers of $M$. 

In general, however, the skeleton of a Weinstein manifold is a highly singular object for which there is no available notion of Floer cohomology (and it seems rather unlikely that such a notion exists). Note that \emph{even for the cotangent bundle of a closed manifold}, equipped with a Weinstein structure homotopic to the standard Morse--Bott one, the skeleton is typically singular.  

To the best of our knowledge, the only systematic result in the literature addressing \Cref{question:arnold-type} is the lower bound $|\fk{c}_X \cap \phi(\fk{c}_X)| \geq 1$, which holds whenever the Rabinowitz--Floer homology $RFH(X)$ is non-zero \cite{c-f}. The nonvanishing of $RFH(-)$ for Weinstein manifolds is now known to be equivalent to the nonvanishing of $SH(-)$ and $\mc{W}(-)$; see \cite[Sec.\ 1.6]{ritter2013topological}. 

\begin{question}\label{question:lefschetz}
Let $X$ be a Liouville manifold, considered up to homotopy. What is the minimal possible number of critical points of a Lefschetz fibration $f: X \to \C$ with Weinstein fibers?
\end{question}
Let us denote this number by $\op{Lef}_w(X) \in \N \cup \{\infty\}$. According to work of Giroux--Pardon \cite{giroux-pardon}, any Weinstein manifold admits a Lefschetz fibration with Weinstein fibers, so $\op{Lef}_w(X)<\infty$ iff $X$ is Weinstein (up to homotopy). It is also clear that $\op{Lef}_w(X)=0$ iff $X$ is subcritical (by a result of Cieliebak \cite{cieliebak-split}, this means equivalently that $X$ admits a Weinstein presentation with only handles of index strictly less than the middle dimension, or that $X$ is Weinstein deformation equivalent to the product of a lower-dimensional Weinstein manifold with $\C$). For topological reasons, one has the lower bound $\op{Lef}_w(X) \geq \op{rk} H_n(X; \Z)$ when $X$ has real dimension $2n$.  These are the only lower bounds in the literature that we are aware of.

An analog of $\op{Lef}_w(X)$ in real Morse theory was considered by Abouzaid--Seidel \cite{AScomplexity} under the name ``complexity". By definition, the complexity of a Weinstein manifold $X$, which we denote by $\op{WMor}(X)$, is the minimal number of critical points of a Weinstein Morse function on $X$, where $X$ is considered up to Weinstein homotopy.\footnote{A \emph{Weinstein Morse} function on $(X, \l)$ is a proper Morse function with respect to which the Liouville vector field is gradient-like. Here $\l$ is a particular choice of Liouville structure on $X$ of the given homotopy class.}  Surprisingly, Lazarev \cite{lazarev-gt} showed that for Weinstein manifolds having non-zero middle homology, one has $\op{WMor}(X)= \op{Mor}(X)$, where $\op{Mor}(X)$ denotes the minimal number of critical points of a \emph{smooth} Morse function on $X$. Does a similar phenomenon hold for $\op{Lef}_w(X)$?

\subsubsection{Results}
We give partial answers to the above questions. As we will explain, our results appear to be the first to go beyond the basic flexible/rigid dichotomy. 

We first consider \Cref{question:arnold-type}. To this end, suppose that $(X, \l)$ be a Weinstein manifold with properly embedded cocores (this is a generic condition). Let $R \sub SH^0(X)$ be a subalgebra which is of finite type over a field $k$, and suppose that there exists a (split-)generating Lagrangian $K \in \mc{W}(X)$ such that $HW^\bu(K,K)$ is a Noetherian $R$-module, where $R$ acts via the closed-open string map. 

\thm[= \Cref{prop_int_rdim} + \Cref{theorem:main-biko-adapt}]\label{theorem:rdim-skeleta-krull}
Under the above assumptions, if $\phi: X \to X$ is a generic compactly-supported Hamiltonian symplectomorphism, then:
\eq\label{equation:rdim-skeleta-krull} | \fk{c}_X \cap \phi(\fk{c}_X) | \geq \dim_R HW^\bu(K,K) +1. \eeq
\ethm
Here $\dim_R (-)$ denotes the \emph{Krull} dimension of an $R$-module. Both $\mc{W}(X)$ and $SH^{\bu}(X)$ are defined over $k$ and are, say, $\Z/2$-graded. Here are some examples illustrating applications of \Cref{theorem:rdim-skeleta-krull}.

\ex
For $(T^*S^n, \l)$, where $n \geq 2$ and $\l$ is homotopic to the standard Morse--Bott form $pdq$, we will see that $\rdim \mc{W}(T^*S^n)=1$ (say with $\Z/2$-gradings and $\Q$-coefficients) and hence $|\fk{c}_X \cap \phi(\fk{c}_X)| \geq 2$. This lower bound is of course sharp.
\eex

\ex
Let $G$ be a compact, simply-connected Lie group. Suppose that $(T^*G, \l)$ is a Weinstein structure which is homotopic to $(T^*G, pdq)$. Then we will see that $\rdim \mc{W}(T^*G) \geq \op{rank} G$ and hence $|\fk{c}_X \cap \phi(\fk{c}_X)| \geq \op{rank} G + 1$. 
\eex

\ex\label{example:intro-hom-section}
Let us consider the case where $2c_1(X) = 0$ so that both $SH^{\bu}(X)$ and $\mc{W}(X)$ are $\Z$-graded. Following \cite{pomerleano2021intrinsic}, an object $L_0 \in \mc{W}(X)$ is called a \emph{homological section} if the restriction of the closed-open map $\mc{CO}$ to the degree $0$ part of $SH^{\bu}(X)$
\eq \mc{CO}: SH^{0}(X) \to HW^{\bu}(L_0, L_0) \eeq
defines an isomorphism of algebras. If $\op{Spec}(SH^{0}(X))$ is smooth, \cite[Cor.\ 1.4]{pomerleano2021intrinsic} tells us $L_0$ is a split-generator. Then we can conclude
\eq  |\fk{c}_X \cap \phi(\fk{c}_X)| \geq \rdim \mc{W}(X) +1 \geq \dim SH^{0}(X) +1. \eeq
In practice, a large class of such Weinstein manifolds arise from the complement of the anti-canonical divisor of log Calabi--Yau varieties.
\eex

We now turn to \Cref{question:lefschetz}. We have the following result:

\thm[= \Cref{prop:lef-rdim} + \Cref{theorem:main-biko-adapt}]\label{theorem:intro-lefschetz}
If $R \sub SH^0(X)$ and $K \in \mc{W}(X)$ are as in \Cref{theorem:rdim-skeleta-krull}, then 
\eq \op{Lef}_w(X) \geq \dim_R HW^\bu(K,K) +1. \eeq
\ethm

\ex
We have $\op{Lef}_w(T^*S^n) = 2$ for $n \geq 2$. If $G$ is a simply-connected compact Lie group, then $\op{Lef}_w(T^*G) \geq \op{rank} G + 1$. If $X$ admits a homological section as in \Cref{example:intro-hom-section}, then $\op{Lef}_w(X) \geq \dim SH^{0}(X) +1$.
\eex

\ex\label{example:intro-point-like}
If $\mc{W}(X)$ is $\Z$-graded and $L \in \mc{W}(X)$ is diffeomorphic to a $n$-torus, then $L$ is a \emph{point-like object}. By work of Elagin--Lunts \cite{elagin-lunts-koszul}, the existence of such an $L$ implies that $\rdim \mc{W}(X) \geq n$ (see \Cref{example:point-like-objects}). Hence $\op{Lef}_w(X) \geq n+1$.
\eex 
As an application of \Cref{theorem:intro-lefschetz}, we prove that the natural complex analog of Lazarev's previously mentioned result $``\op{WMor}(X)= \op{Mor}(X)"$ \emph{fails}. This is a consequence of the following corollary. 

\cor[= \Cref{corollary:lef-distinguish}]\label{corollary:lef-distinguishes}
There exists a Liouville manifold $T^*S^3_{\op{exotic}}= (T^*S^3, \l)$ which is formally isotopic to $(T^* S^3, pdq)$, but where
\eq \op{Lef}_w (T^* S^3) \neq \op{Lef}_w(T^*S^3_{\op{exotic}}). \eeq
\ecor

The term ``formally isotopic" should be interpreted in the context of the h-principle: it means that $pdq$ and $\l$ are indistinguishable from the perspective of differential topology. \Cref{corollary:lef-distinguishes} provides the first example of a pair of formally isotopic Weinstein manifolds with non-zero middle homology which can be distinguished by $\op{Lef}_w(X)$.\footnote{For Weinstein manifolds with vanishing middle homology, $\op{Lef}_w(-)$ distinguishes exotic examples for trivial reasons: for example, if $X= \C^n$ and $X'$ is an exotic $\C^n$, then $\op{Lef}_w(X)=0$ but $\op{Lef}_w(X')>0$ since otherwise $X'$ would be subcritical and hence standard.} In particular, \Cref{corollary:lef-distinguishes} implies that, in contrast to $\op{WMor}(-)$, the number $\op{Lef}_w(-)$ is a truly symplectic invariant.
 
The proof of \Cref{corollary:lef-distinguishes} uses a construction of Eliashberg--Ganatra--Lazarev \cite{e-g-l}, who exhibited an exotic $T^*S^3_{\op{exotic}}$ which contains a regular Lagrangian $3$-torus. We will see that this implies that $\rdim \mc{W}(T^*S^3_{\op{exotic}}) = 3$, and hence $\op{Lef}_w(T^*S^3_{\op{exotic}}) \geq 4$. In contrast, it is known that $T^*S^3$ admits a Lefschetz fibration with $2$ critical points.

\subsubsection{Methods: bounds on the Rouquier dimension coming from geometry}
It turns out that the Rouquier dimension of wrapped Fukaya categories provides a lower bound for the quantities we wish to study. 

\prop[= \Cref{prop_int_rdim}] \label{prop:intromain-central-bound}
Let $(X, \l)$ be a Weinstein manifold with properly embedded cocores (this is a generic condition). If $\phi: X \to X$ is a generic compactly-supported Hamiltonian symplectomorphism, then: 
\eq\label{equation:rdim-skeleta} | \fk{c}_X \cap \phi(\fk{c}_X) | \geq \rdim \mc{W}(X) +1. \eeq
\eprop

(The result we prove in \Cref{prop_int_rdim} is stated slightly differently and is more general.) Let us sketch the proof of \Cref{prop:intromain-central-bound}. To begin with, it can be shown that the Rouquier dimension of any $A_\infty$ category is bounded above by the length of the shortest resolution of the diagonal bimodule by Yoneda bimodules (minus $1$). This is a purely algebraic fact which has nothing to do with Fukaya categories. 

We now set $r=|\fk{c}_X \cap \phi(\fk{c}_X)|$. Our goal is to construct a resolution of the diagonal bimodule $\Delta_{\mc{W}(X)}$ of length $r$. We do this by the following geometric argument. First, observe that the intersection points of $\fk{c}_X$ with $\phi(\fk{c}_X)$ are in bijection with the intersection points of the skeleton of $(\ov{X} \tms X, -\lambda \oplus \phi_*\l)$ with the diagonal Lagrangian. Since $\fk{c}_X \pitchfork \phi(\fk{c}_X)$, the diagonal intersects the skeleton transversally in $r$ points. By applying deep work of Chantraine--Dimitroglou-Rizell--Ghiggini--Golovko \cite{cdgg} or Ganatra--Pardon--Shende \cite{gps2}, one can then construct a resolution of the diagonal by $r$ cocores of $\ov{X} \tms X$. Morally, the idea is that as one scales the diagonal Lagrangian by the Liouville flow, it converges to cocores at each point where it intersects the skeleton. The upshot is that we have resolved the diagonal Lagrangian by $r$ cocores in $\ov{X} \tms X$, which are just products of cocores in $X$. Finally,  we can turn this resolution of the diagonal Lagrangian into a resolution of the diagonal bimodule. This last step is carried out in the Appendix, following the methods of \cite{ganatra} and \cite{gps2}.

There is also an upper bound on the Rouquier dimension coming from Lefschetz fibrations.

\prop[= \Cref{prop:lef-rdim}]\label{proposition:lefschetz-lower-bound}
Let $X$ be a Liouville manifold. Then we have
\eq \op{Lef}_w(X) \geq \rdim \mc{W}(X) +1. \eeq
\eprop

The proof of \Cref{proposition:lefschetz-lower-bound} is straightforward and (at least implicitly) well-known to experts: the Fukaya--Seidel category associated to a Lefschetz fibration admits a semi-orthogonal decomposition of length equal to the number of critical points of the fibration, which immediately implies the corresponding upper bound for the Rouquier dimension. 

\subsubsection{Methods: commutative ring actions on wrapped Fukaya categories}\label{subsubsection:biko-adapt}
In order to obtain conclusions about quantitative symplectic topology, we need a method for obtaining lower bounds on the Rouquier dimension of wrapped Fukaya categories which can then be combined with the upper bounds stated in Propositions \ref{prop:intromain-central-bound} and \ref{proposition:lefschetz-lower-bound}. These lower bounds will ultimately arise from considering a canonical action of the Hochschild cohomology ring on the Fukaya category. To explain this, it is natural to start the discussion in a more general setting.

Let $\mc{C}$ be a pre-triangulated $A_\infty$ category and let $\mc{T}= H^0(\mc{C})$ be its homotopy category. Let $Z^\bu(\mc{T})$ be the (graded-commutative) ring of graded natural transformations of the identity on $\mc{T}$. Given a graded-commutative ring $R$, a \emph{central action} of $R$ on the triangulated category $\mc{T}$ is by definition a morphism of graded-commutative rings $R \to Z^\bu(\mc{T})$. Concretely, this is the data of a morphism $R \to \hom_{H(\mc{C})}^\bu(K, K)$ for each object $K \in \mc{C}$ such that the induced left and right module structures on $\hom_{H(\mc{C})}^\bu(K, L)$ coincide up to sign. 

The upshot is that the morphism spaces of objects of $\mc{T}$ define modules over $R$. One can then study the triangulated category $\mc{T}$ using tools from (graded-)commutative algebra. This perspective leads to a very rich theory; see e.g.\ \cite{benson2008local}.

The workhorse result for lower-bounding the Rouquier dimension of wrapped Fukaya categories in this paper is the following:
\thm[= \Cref{theorem:main-biko-adapt}]\label{theorem:intro-biko-adapt}
Let $\mc{T}$ be a $\Z/m$-graded triangulated category $(1 \leq m \leq \infty)$ and let $R$ be a finite type $k$-algebra acting centrally on $\mc{T}$. Suppose that there exists a (split-)generator $G \in \mc{T}$ such that $\op{Hom}_{\mc{T}}^\bu(G,G)$ is a Noetherian $R$-module. Then $\rdim \mc{T} \geq \op{dim}_R \op{Hom}_{\mc{T}}^\bu(G,G).$
\ethm

\Cref{theorem:main-biko-adapt} can be viewed as an ``affine" analog of a beautiful theorem of Bergh--Iyengar--Krause--Oppermann \cite[Thm.\ 4.2]{b-i-k-o}. In fact, (as we will see) one could directly apply \cite[Thm.\ 4.2]{b-i-k-o} and still obtain useful lower bounds for certain examples. However, for our purposes, \cite[Thm.\ 4.2]{b-i-k-o} suffers from two limitations. First of all, it is essentially never sharp (there is a pesky $-1$ which comes, ultimately, from working in the ``projective" rather than ``affine" setting). Secondly, it requires rather strong hypotheses on $R$ and $\op{Hom}_{\mc{T}}^\bu(G,G)$ which do not hold in many cases of interest to us (for instance, \Cref{example:intro-hom-section} and \Cref{example:computation-n-torus} would fail).\footnote{There is a typo in the definition of the projective dimension in \cite[Sec.\ 2]{b-i-k-o} which propagates. As a result, one needs to add an assumption in the statement of \cite[Thm.\ 4.2]{b-i-k-o}, the most natural (to us) being that $R^0=k$.} In contrast, \Cref{theorem:main-biko-adapt} is sharp for many examples (such as cotangent bundles of $n$-spheres and $n$-tori), and more widely applicable. 

Our general approach to proving \Cref{theorem:main-biko-adapt} via \emph{Koszul objects} follows \cite{b-i-k-o}; however, the input from graded-commutative algebra is replaced by new arguments from (ungraded) commutative algebra. We crucially use the hypothesis that $R$ is of finite type. This hypothesis is rather severe from the representation theoretic perspective taken in \cite{b-i-k-o}. However, from the symplectic perspective which we adopt in this paper, this assumption is rather mild. 

Of course, \Cref{theorem:intro-biko-adapt} is useless without a good supply of central actions. Happily, for any pre-triangulated $A_\infty$ category $\mc{C}$, there is a canonical central action of the Hochschild cohomology ring $HH^\bu(\mc{C})$ on $\mc{T}= H^0(\mc{C})$ called the \emph{characteristic morphism}.  Given any object $K \in \mc{C}$, this is just the map $HH^\bu(\mc{C}) \to \hom^\bu_{H(\mc{C})}(K, K)$ that projects a Hochschild cochain onto its length zero part. In general, the characteristic morphism is neither injective nor surjective.

In the case which is relevant to us, namely when $\mc{T}=H^0(\op{Perf} \mc{W}(X))$ and $X$ is Weinstein, the characteristic morphism is known to admit a concrete description in terms of the familiar \emph{closed-open map}. To explain this a little, recall that the closed-open map defines a ring morphism $SH^\bu(X) \to HW^\bu(L, L)$, for any object $L \in \mc{W}(X)$. On the other hand, if $X$ is Weinstein, there is an isomorphism of graded-commutative rings $HH^\bu(\mc{W}(X))= SH^\bu(X)$ \cite{cdgg}. Now the point is simply that the obvious diagram formed by these maps commutes. As a result, we obtain the following corollary which is the desired lower bound. 

\cor[= \Cref{corollary:main-central-bound}]\label{corollary:main-central-bound-intro}
Suppose that $SH^0(X)$ admits a subring $R$ of finite-type over $k$, and suppose there exists a (split-)generator $K \in \mc{W}(X)$ such that $HW^\bu(K,K)$ is Noetherian over $R$, where $R$ acts via the closed-open map. Then 
\eq\label{equation:rdim-lower-bd-intro} \rdim \mc{W}(X) \geq \dim_R HW^\bu(K,K). \eeq
\ecor
\Cref{corollary:main-central-bound-intro} works with any choice of gradings, but $SH^\bu(X)$ will typically be $\Z$-graded or $\Z/2$-graded, depending on the application. The flexibility to choose different gradings is convenient in practice, since of course $SH^0(X)$ will admit different subrings depending on the choice of grading.

\subsection{Acknowledgements} 
We are indebted to Daniel \'{A}lvarez-Gavela for answering countless questions about arborealization, Vivek Shende for help with microlocal sheaves and homological mirror symmetry of toric varieties, Srikanth Iyengar for discussions and comments related to \Cref{subsection:commutative-alg-part}, Y.\ Bar\i\c{s} Kartal for telling us the proof of \Cref{prop:lef-rdim} and for useful discussions related to \Cref{subsection:motives}, and John Pardon for many helpful conversations about various aspects of the paper. In addition, we would like to thank Denis Auroux, Petter Bergh, Georgios Dimitroglou Rizell, Sheel Ganatra, Wahei Hara, Oleg Lazarev, David Nadler, Amnon Neeman and Dmitri Orlov for answering our questions at various points in this project. 

SB would like to express deep gratitude to his advisor John Pardon for constant encouragement and support. LC learned about \Cref{question:arnold-type} from Zhengyi Zhou and would like to thank him for many enjoyable and informative discussions. LC would also like to thank Y.\ Bar\i\c{s} Kartal for teaching him much of what he knows about $A_\infty$ categories. Finally, we would like to thank the anonymous referee for many useful comments which have significantly improved this paper.

Part of this work was completed while LC was a Member of the Institute of Advanced Study, and supported there by the National Science Foundation under Grant No.\ DMS-1926686. LC was also partially supported by Simons Foundation grant 385573 (Simons Collaboration on Homological Mirror Symmetry).

\section{Category theory}

\subsection{Notation and conventions}
Unless otherwise indicated, $k$ denotes a field.

Let $\mc{T}$ be a category. We will routinely abuse notation by writing $K \in \mc{T}$ to mean that $K$ is an object of $\mc{T}$. Similarly, we write $\mc{I} \sub \mc{T}$ to mean that $\mc{I}$ is a full subcategory of $\mc{T}$. Given two full subcategories $\mc{I}, \mc{J}$ of $\mc{T}$, we let $\mc{I} \cup \mc{J}$ be the full subcategory with objects the union of the objects in $\mc{I}$ and those in $\mc{J}$. 

A full subcategory $\mc{I} \sub \mc{T}$ is said to be \emph{strictly full} if is closed under isomorphisms. We will sometimes implicitly identify a set of objects with the corresponding full subcategory (resp.\ identify a set of isomorphism classes of objects with the corresponding strictly full subcategory). 

Given an object $K \in \mc{T}$, a \emph{summand} of $K$ is a triple $(Z, k, r)$ consists of an object $Z \in \mc{T}$ along with maps $r \in \hom(Z, K), k \in \hom(K, Z)$ such that $k \circ r=\op{id}$ and $r \circ k$ is an idempotent.  An idempotent $e \in \hom(K, K)$ is said to \emph{split} if $e = r \circ k$ for some summand $(Z, k, r)$. A category is said to be \emph{idempotent complete} (or Karoubi complete) if all idempotents split.

\subsection{Triangulated categories}\label{subsection:triangulated-cats}

A triangulated category $\mc{T}=(\mc{T}, \S, \Delta)$ is an additive category $\mc{T}$ equipped with an autoequivalence $\S: \mc{T} \to \mc{T}$, called the \emph{shift functor}, and a class of distinguished triangles $\Delta$. This data is required to satisfy certain axioms. A standard reference for triangulated categories is \cite{neeman}. 

We say that a triangulated category $\mc{T}$ is $\Z/m$-graded (for $m \in \N_+ \cup \{\infty\}$ and $\Z/\infty:= \Z$) if $\S^m= \op{id}$. If $\mc{T}$ is $\Z/m$-graded, we define $\op{Hom}_\mc{T}^\bu(X, Y):= \bigoplus_{i=1}^m\Hom_\mc{T}(X, \S^i Y)$ for $m<\infty$ and $\op{Hom}_\mc{T}^\bu(X, Y):= \bigoplus_{i \in \Z}\Hom_\mc{T}(X, \S^i Y)$ for $m =\infty$. 

Given full subcategories $\mc{I}_0, \mc{I}_1$ of $\mc{T}$, let $\mc{I}_0 * \mc{I}_1$ be the strictly full subcategory which is uniquely characterized by the following property: an object $K \in \mc{T}$ is contained in $\mc{I}_0 * \mc{I}_1$ if there exists a distinguished triangle $K_0 \to K \to K_1 \to$ with $K_0 \in \mc{I}_0$ and $K_1 \in \mc{I}_1$.  It is a straightforward consequence of the octahedral axiom for triangulated categories that the operation $(\mc{I}_0, \mc{I}_1) \mapsto \mc{I}_0 * \mc{I}_1$ is associative (see e.g.\ \cite[Lem.\ 1.3.10]{b-b-d}). 

Given $\mc{I} \sub \mc{T}$ a full subcategory, let $[\mc{I}]$ be the smallest strictly full subcategory of $\mc{T}$ containing $\mc{I}$ and closed under finite direct sums and shifts. Let $\langle \mc{I} \rangle$ be the smallest strictly full subcategory of $\mc{T}$ closed under finite direct sums, shifts and taking summands. Let $\op{ads}(\mc{I})$ be the smallest strictly full subcategory of $\mc{T}$ containing $\mc{I}$ and closed under \emph{arbitrary} direct sums and shifts. 

We set $[\mc{I}]_0:= 0$ and we inductively define $[\mc{I}]_n= [\mc{I}]_{n-1} * [\mc{I}]$ and $\langle \mc{I} \rangle_n = \langle \langle \mc{I} \rangle_{n-1} * \langle \mc{I} \rangle  \rangle$. We set $[\mc{I}]_\infty := \cup_{n \in \N} [\mc{I}]_n$ and $\langle \mc{I} \rangle_\infty = \cup_{n \in \N} \langle \mc{I} \rangle_n$. 

Finally, we let $\mc{I}_{sd}$ be the smallest strictly full subcategory containing all summands of $\mc{I}$ (this is often called the ``Karoubi completion"). A full subcategory of a triangulated category is called \emph{thick} if it is closed under taking summands. 

\rmk
Some sources (e.g.\ \cite{elagin-lunts}) adopt a different convention in which the indices in $\langle - \rangle_n$ and $[-]_n$ are shifted down by $1$. 
\ermk


\lem\label{lemma:operations-props}
The following properties hold for all $n \in \N \cup \{\infty\}$ and any $\mc{I} \sub \mc{T}$.
\begin{enumerate}

\item\label{item:idempotent-ops}\label{item:summand-sums} $[\mc{I} ] = [[ \mc{I} ]]$,  $\mc{I}_\sd = ( \mc{I}_\sd)_\sd$ and $[ \mc{I}_\sd ] \sub [ \mc{I} ]_\sd$
\item\label{item:rouquier-lemma} $\langle \mc{I} \rangle_n = ( [\mc{I}]_n)_\sd$ 
\item\label{item:bracket-bracket} $ [ [ \mc{I}]_n ] = [ \mc{I}]_n$ and $\langle \langle \mc{I} \rangle_n \rangle = \langle \mc{I} \rangle_n$
\item\label{item:sums-cup} Given objects $G, H \in \mc{T}$, we have $\langle G \rangle_n \sub \langle G \oplus H \rangle_n$ and $\langle H \rangle_n \sub \langle G \oplus H \rangle_n$.
\end{enumerate}
\elem

\pf
\eqref{item:idempotent-ops} is obvious. \eqref{item:rouquier-lemma} can be established by following the argument of \cite[Lem.\ 4.1]{seidel-book} (and is stated without proof in \cite[Rmk.\ 3.1]{rouquier}). \eqref{item:bracket-bracket} is straightforward using \eqref{item:rouquier-lemma}. \eqref{item:sums-cup} is obvious.
\epf
%

\lem\label{lemma:semi-orthogonality}
Let $\mc{T}$ be a triangulated category. Let $\mc{I}, \mc{J}$ be thick full triangulated subcategories such that $\mc{T}= \langle \mc{I} \cup \mc{J} \rangle_\infty$ and $\hom(K, L)=0$ for all $K \in \mc{J}$ and $L \in \mc{I}$. Then $\mc{T}= \langle \mc{J}  * \mc{I} \rangle$. 
\elem
\pf The proof is standard; see e.g.\ \cite[Lem.\ 3.1]{bondal}. Since $\langle \mc{J}  * \mc{I} \rangle$ contains $\mc{I}$ and $\mc{J}$, it is enough to verify that $\langle \mc{J}  * \mc{I} \rangle$ is closed under cones. This is a straightforward check using \cite[Prop.\ 1.1.11]{b-b-d}.  
\epf


\lem\label{lemma:orthogonality}
Let $\mc{T}$ be a triangulated category. Let $\mc{I}_1, \dots, \mc{I}_r$ be full triangulated subcategories which are mutually orthogonal, in the sense that $\hom(K, L)=0$ for all $K \in \mc{I}_k$ and $L \in \mc{I}_l$ if $k \neq l$. Then $\langle \mc{I}_1 \cup \dots \cup \mc{I}_r \rangle_n = \langle \langle \mc{I}_1 \rangle_n \cup \dots \cup \langle \mc{I}_r \rangle_n \rangle$ for all $n \in \N \cup \{\infty\}$. 
\qed
\elem

The following lemma is well-known; see e.g.\ \cite[Cor.\ 3.14]{rouquier}.
\lem\label{lemma:split-compact}
Let $\mc{T}$ be a triangulated category and let $\mc{T}^c$ denote the full triangulated subcategory of compact objects. Then given $\mc{I} \sub \mc{T}^c$, we have $\mc{T}^c \cap \langle \op{ads}(\mc{I}) \rangle_d = \langle \mc{I} \rangle_d$.
\qed
\elem

\subsection{$A_\infty$ categories}

\subsubsection{Basic definitions}
Unless otherwise indicated, we follow the conventions of \cite{seidel-book} and \cite{seidel-cat-dynamics} when discussing $A_\infty$ categories. We will work mostly with $\Z/2$-graded or $\Z$-graded $A_\infty$ categories, which are always assumed to be defined over a field $k$. However, the categorical reasoning in this paper applies with purely superficial changes to the setting of $\Z/n$-graded $A_\infty$ categories for any $n \in \N_+$ (the case $n=1$ corresponds to an ungraded $A_\infty$ category).

Given $A_\infty$ categories $\mc{C}_1$ and $\mc{C}_2$, an $(\mc{C}_1, \mc{C}_2)$-bimodule is an $A_\infty$ functor $\mc{C}_1^{\op{op}} \tms \mc{C}_2 \to \op{Ch} k$, where $\op{Ch} k$ is the dg category of chain complexes of $k$-vector spaces. Note that $(\mc{C}_1, \mc{C}_2)$-bimodules form a dg category, which is denoted by $(\mc{C}_1, \mc{C}_2)-\op{mod}$. As usual, two such bimodules are said to be quasi-isomorphic if they are isomorphic in the homotopy category $H^0((\mc{C}_1, \mc{C}_2)-\op{mod})$. Let $\Delta_{\mc{C}}$ be the diagonal bimodule over an $A_\infty$ category $\mc{C}$. A right (resp.\ left) module over $\mc{C}$ is a $(k \tms \mc{C})$-bimodule (resp. a $(\mc{C} \tms k)$-bimodule). The right (resp.\ left) modules over $\mc{C}$ form a dg category $\op{Mod} \mc{C}$ (resp.\ $\op{Mod} \mc{C}^{\op{op}}$).


Let $\{ \mc{C}_i \}_{i=1,2,3}$ be $A_\infty$ categories. Given a $(\mc{C}_1,\mc{C}_2)$-bimodule $\mc{P}$ and a $(\mc{C}_2, \mc{C}_3)$-bimodule $\mc{Q}$, we can form their convolution (bimodule tensor product) $\mc{P} \otimes_{\mc{C}_2} \mc{Q}$. Using the same notations, the convolution has the following properties:
\begin{itemize}
\item it is strictly associative; 
\item we have quasi-isomorphisms between bimodules $\Delta_{\mc{C}_1} \otimes_{\mc{C}_1} \mc{P} \to \mc{P}$ and $\mc{P} \otimes_{\mc{C}_2} \Delta_{\mc{C}_2} \to \mc{P}$.
\item the convolution $- \otimes_{\mc{C}_2} \mc{Q}$ induces a dg functor $(\mc{C}_1, \mc{C}_2)-\op{mod}$ to $(\mc{C}_1, \mc{C}_3)-\op{mod}$. Similarly for $\mc{Q} \otimes_{\mc{C}_3} -$.
\end{itemize}

There is a canonical Yoneda embedding $\mc{C} \to \op{Mod} \mc{C}$. Let $\op{Tw} \mc{C} \sub \op{Mod} \mc{C}$ be the closure of the image of the Yoneda embedding under taking mapping cones. Let $\op{Perf} \mc{C} \sub \op{Mod} \mc{C}$ be the closure of the image of the Yoneda embedding under taking mapping cones and summands. Note that we have
\eq\label{equation:perf-compact} \op{Perf} \mc{C} = (\op{Mod} \mc{C})^c, \eeq
where $(-)^c$ means taking compact objects; see \cite[(A.2) and Prop.\ A.3]{gps3}.

An $A_\infty$ category $\mc{C}$ is said to be \emph{pre-triangulated} if the canonical embedding $\mc{C}\hookrightarrow \op{Tw} \mc{C}$ is a quasi-equivalence. The homotopy category $H^0(\mc{C})$ of a pre-triangulated $A_\infty$ category is a triangulated category.

Given an $A_\infty$ category $\mc{C}$ and a set of objects $\mc{A}$, we can form the \emph{quotient} $A_\infty$ category $\mc{C}/\mc{A}$ (see \cite{lyubashenko-ovsienko, drinfeld}), which comes equipped with a canonical map $q: \mc{C} \to \mc{C}/\mc{A}$. Suppose $\mc{B}$ is another $A_\infty$ category. Given a functor $\mc{C} \to \mc{B}$ which sends $\mc{A}$ to acyclic objects, there is an induced functor $\mc{C}/\mc{A} \to \mc{B}$ (see \cite[Sec.\ A.7]{gps3}). We note that the canonical map $(\op{Tw} \mc{C})/ \mc{A} \to \op{Tw} (\mc{C}/\mc{A})$ is a quasi-equivalence, and the canonical map $\op{Perf} \mc{C}/ \mc{A} \to \op{Perf}(\mc{C}/\mc{A})$ is a Morita equivalence (it is not in general a quasi-equivalence; see \cite[Sec.\ 2]{p-s} for concrete counterexamples).  

\subsubsection{Equivalent constructions of the derived category of a (dg) $k$-algebra}

In the literature, one encounters various definitions of the derived category of a (dg) $k$-algebra. We therefore collect some standard facts which will be implicitly assumed in the sequel.

\lem\label{lemma:derived-cats-dg}
Let $\mc{A}$ be a dg $k$-algebra. The following constructions of the derived category $D(\mc{A})$ produce equivalent triangulated categories. 
\begin{enumerate} 
\item Consider the quotient (in the sense of Drinfeld \cite{drinfeld}) of the dg category of dg $\mc{A}$-modules by the dg subcategory of acyclic dg $\mc{A}$-modules. Then pass to $H^0(-)$. 
\item View $\mc{A}$ as an $A_\infty$ algebra and consider $\op{Mod} \mc{A}$, the $A_\infty$ category of $A_\infty$ $\mc{A}$-modules.  Now pass to $H^0(-)$. 
\end{enumerate}
\qed
\elem

\lem\label{lemma:derived-cats}
Let $A$ be a $k$-algebra. The following constructions of the derived category $D(A)$ produce equivalent triangulated categories. 
\begin{enumerate} 
\item\label{item:abelian-def} Consider the abelian category of modules over $A$ and let $K(A)$ be the homotopy category. Then take the quotient (in the sense of Verdier) of $K(A)$ by the subcategory of acyclic complexes. 
\item\label{item:dg-def} View $A$ as a dg algebra concentrated in degree zero with the trivial differential. Now apply either one of the equivalent constructions of \Cref{lemma:derived-cats-dg}. 
\end{enumerate}
\qed
\elem
In contrast to \Cref{lemma:derived-cats}, if $A$ is a graded $k$ algebra, it is in general \emph{not} the case that the derived category of the abelian category of graded $A$-modules coincides with the derived category of dg $A$-modules.

\rmk\label{remark:perf}
Let $A$ be a $k$-algebra. Viewing $A$ as a dg $k$-algebra, we can consider the $A_\infty$ category $\op{Perf} A$. Then $H^0(\op{Perf}A) \sub H^0(\op{Mod} A)= D(A)$ is equivalent to the smallest idempotent complete triangulated subcategory of $D(A)$ containing $A$ (where $A$ is viewed as a module over itself). It is shown in \cite[Prop.\ 07LT]{stacks-project} that the objects of this latter category are precisely the perfect complexes, i.e.\ the complexes of $A$-modules quasi-isomorphic to a bounded complex of finite projective $A$-modules.
\ermk

\subsection{The Rouquier dimension}

In his remarkable paper \cite{rouquier}, Rouquier introduced the following notion of dimension for triangulated categories. 

\defi\label{definition:rdim}
The \emph{Rouquier dimension} of a triangulated category $\mc{T}$ is the smallest $n \in \Z_{\geq -1} \cup \{ \infty\}$ such that there exists an object $G \in \mc{T}$ with $\mc{T}= \langle G \rangle_{n+1}$. We denote the Rouquier dimension of $\mc{T}$ by $\rdim{\mc{T}}$. 

The Rouquier dimension of an $A_\infty$ category $\mc{C}$ is defined to be the Rouquier dimension of $H^0(\op{Perf} \mc{C})$ and is denoted by $\rdim{\mc{C}}$. 
\edefi

We now collect some standard properties of the Rouquier dimension which will be needed in the sequel. The following lemma is a useful inequality for the Rouquier dimension of semi-orthogonal decompositions. 

\lem\label{lemma:sod-rdim}
Let $\mc{T}$ be a triangulated category. Let $\mc{I}_1,\dots,\mc{I}_m$ be thick full triangulated subcategories satisfying the following conditions:
\begin{itemize}
\item[(i)] for any $j>i$, we have $\hom(K, L)=0$ for all $K \in \mc{I}_j$ and $L \in \mc{I}_i$; 
\item[(ii)] $\langle \mc{I}_1 \cup \dots \cup \mc{I}_m \rangle_\infty = \mc{T}$. 
\end{itemize}

Then 
\eq\label{equation:sod-inequality} \rdim{\mc{T}} \leq \sum_{l=1}^m (\rdim{\mc{I}_l}+1) -1. \eeq
\elem
\eqref{equation:sod-inequality} is in general far from sharp. For example, if $\G$ is a quiver whose underlying graph is a Dynkin diagram of type ADE, then $\rdim (\op{Perf} k[\G])=0$ (see \Cref{proposition:elagin-gabriel}). However, $\op{Perf} k[\G]$ admits a full exceptional collection of length $|\G|$ (see \cite[Sec.\ 2.4]{elagin2020calculating}). 

\pf
Note that $\mc{T}= \langle \mc{I}_1 \cup \dots \cup \mc{I}_m \rangle_\infty = \langle \langle \mc{I}_1 \cup \dots \cup \mc{I}_{m-1} \rangle_\infty \cup \mc{I}_m \rangle_\infty$. Hence it is sufficient to treat the case where $m=2$. 

Suppose $\rdim{\mc{I}_1} = n_1$ and $\rdim{\mc{I}_2}=n_2$. This means that there exist objects $G \in \mc{I}_1$ and $H \in \mc{I}_2$ such that $\mc{I}_1 = \langle G \rangle_{n_1+1}$ and $\mc{I}_2= \langle H \rangle_{n_2+1}$. 

We now have, by \Cref{lemma:semi-orthogonality}, $\mc{T} = \langle \langle H \rangle_{n_2+1} * \langle G \rangle_{n_1+1} \rangle \sub  \langle \langle G \oplus H \rangle_{n_2+1} * \langle G \oplus H \rangle_{n_1+1} \rangle = \langle ([G \oplus H]_{n_2+1})_\sd *  ([G \oplus H]_{n_1+1})_\sd \rangle \sub \langle ([G \oplus H]_{n_2+1}*  [G \oplus H]_{n_1+1})_\sd \rangle = \langle ([G \oplus H]_{n_1+n_2+2})_\sd \rangle = \langle \langle G \oplus H \rangle_{n_1+n_2+2} \rangle$.
\epf

The following statements show that the Rouquier dimension is non-increasing under quotients of categories.

\lem[Lem.\ 3.4 in \cite{rouquier}]\label{lemma:dense-image}
Let $\mc{T}$ and $\mc{T}'$ be triangulated categories. If there exists an exact functor $\mc{T} \to \mc{T}'$ with dense image (i.e.\ every object in $\mc{T}'$ is isomorphic to a summand of an object in the image), then $\rdim{\mc{T}'} \leq \rdim{\mc{T}}$.
\qed
\elem

\cor \label{corollary:quot-rdim}
Let $\mc{C}$ be an $A_\infty$ category and let $\mc{A}$ denote a set of objects. Then $\rdim{\mc{C}/\mc{A}} \leq \rdim{\mc{C}}$.
\qed
\ecor

When a triangulated category comes from an $A_\infty$ category, the following lemma shows that the Rouquier dimension is bounded above by the length of the shortest resolution of the diagonal bimodule. 

\lem\label{lemma:resolution-diag}
Let $\mc{C}$ be an $A_\infty$ category and let $\Delta_\mc{C}$ denote the diagonal bimodule. Suppose there exists a resolution of length $\ell \in \N$
\eq 0 = \mc{M}_0 \to \mc{M}_1 \to \dots \to \mc{M}_\ell \eeq
of $(\mc{C}, \mc{C}$)-bimodules where 
\begin{itemize}
\item $\Delta_\mc{C}$ is a summand of $\mc{M}_\ell$; 
\item there exists perfect modules $\mc{P}^r_i \in \op{Perf} \mc{C}$ and $\mc{P}^l_i \in \op{Perf} \mc{C}^{\op{op}}$ such that $\op{cone}(\mc{M}_{i-1} \to \mc{M}_i) \in [\mc{P}^l_i \otimes \mc{P}^r_i]$ in $H^0( (\mc{C}, \mc{C})-\op{mod})$. 
\end{itemize}
Then $\rdim{\mc{C}} \leq \ell -1$.
\elem
\pf
Given $Q \in \op{Perf} \mc{C}$, we convolve $\mc{Q}$ with the above resolution (i.e.\ we apply the functor $\mc{Q} \otimes_{\mc{C}} -$). Note that $\mc{Q} \otimes_\mc{C} \mc{P}^l_i \otimes \mc{P}^r_i \simeq H^\bu(\mc{Q} \otimes_\mc{C} \mc{P}^l_i) \otimes  \mc{P}^r_i$, since we are working over a field. Writing $C_i:= H^\bu(\mc{Q} \otimes_\mc{C} \mc{P}^l_i)$, it follows that $\mc{Q} \in \langle \oplus_{i=1}^\ell  (C_i \otimes \mc{P}^r_i) \rangle_\ell \sub  \langle \op{ads}(\mc{G}) \rangle_\ell$, where $\mc{G}= \oplus_{i=1}^\ell \mc{P}^r_i$. Since $\mc{Q}$ and $\mc{G}$ are compact objects in $\op{Mod} \mc{C}$ (see \eqref{equation:perf-compact}), this implies that $\mc{Q} \in \langle \mc{G} \rangle_\ell$ by \Cref{lemma:split-compact}. 
\epf

\subsection{Homotopy colimits}

Given a diagram of $A_\infty$ categories $\{\mc{C}_\s \}_{\s \in \S}$ indexed by a poset $\S$, one can form the \emph{Grothendieck construction} (or semiorthogonal gluing) 
\eq \op{Groth}_{\s \in \S} \mc{C}_\s. \eeq 
This is an $A_\infty$ category with $\op{Ob} (\op{Groth}_{\s \in \S} \mc{C}_\s)= \bigsqcup_{\s \in \S} \op{Ob} \mc{C}_\s$. It satisfies the following important property: if $K \in \mc{C}_{\s}$ and $L \in \mc{C}_{\s'}$ are objects, then the space of morphism from $K$ to $L$ in $\op{Groth}_{\s \in \S} \mc{C}_\s$ is zero unless $\s \leq \s'$.

The category $\op{Groth}_{\s \in \S} \mc{C}_\s$ admits a distinguished collection of morphism $A_\S$ called \emph{adjacency morphisms}. The \emph{homotopy colimit} of the diagram $\{\mc{C}_\s \}_{\s \in \S}$ may be defined as the localization of the Grothendieck construction at the adjacency morphisms:
\eq \label{equation:homotopy-colim-def}
\op{hocolim}_{\s \in \S}:= \lt( \op{Groth}_{\s \in \S} \mc{C}_\s \rt)[A_\S]^{-1}.
\eeq

We refer the reader to \cite[A.4]{gps2} for details on the Grothendieck construction and homotopy colimit. We will not need any properties of these constructions beyond those mentioned above. (Note that \eqref{equation:homotopy-colim-def} should really be viewed as a construction which happens to compute the homotopy colimit, in the sense that it satisfies an appropriate universal property as discussed for instance in \cite[Sec.\ 1.2.13]{lurie-htt}. However, the reader is free to take it as a definition.)

We now collect some lemmas about the Rouquier dimension of homotopy colimits. We begin with a small piece of notation: given a finite poset $\S$ and an element $\s \in \S$, we let $h(\s) \in \N_+$ be the length of the longest chain $\s_1<\dots<\s_{h(\s)}=\s$ terminating at $\s$. Note that $h(\s)= h(\s')$ implies that $\s$ and $\s'$ are either incomparable or equal. Let $D$ be the maximum of $h(\s)$ for all $\s \in \S$.

Let us now consider the triangulated categories $\mc{T}= H^0(\op{Perf} \op{Groth}_{\s \in \S} \mc{C}_\s)$ and $\mc{T}'= H^0(\op{Perf} \op{hocolim}_{\s \in \S} \mc{C}_\s)$. We write  $\mc{I}_\ell = \langle \bigcup_{h(\s)=\ell} \mc{C}_\s \rangle \sub \mc{T}$. 

\lem\label{lem_rdim_max}
$\rdim{\mc{I}_\ell} = \op{max}_{h(\s)=\ell} \rdim \mc{C}_\s$. 
\elem
\pf 
Follows from \Cref{lemma:orthogonality} and the fact that $\mc{C}_\s, \mc{C}_{\s'}$ are orthogonal whenever $\s \neq \s'$.
\epf

\lem\label{lem_rdim_upper}
We have $\rdim{\mc{T}'} \leq \rdim{\mc{T}} \leq \sum_{\ell=1}^D (\rdim{\mc{I}_\ell} +1)-1$. 
\elem
\pf
The first inequality follows from \Cref{corollary:quot-rdim}. For the second one, observe that $\langle \bigcup_{\ell=1}^L \mc{I}_\ell \rangle_\infty = \mc{T}$. As noted above, we have $\hom_\mc{T}(K, L)=0$ if $K \in \mc{I}_j$ and $L \in \mc{I}_i$ for $j>i$. The claim thus follows from \Cref{lemma:sod-rdim}. 
\epf

\subsection{The Rouquier dimension of some triangulated categories}\label{subsection:first-examples}

We collect some computations from the literature to be used later. For many other interesting computations, we refer to \cite{elagin2020calculating}.

\subsubsection{Coherent sheaves}\label{subsubsection:coherent-shaves-computations}

We let $D^b \op{Coh}(Y)$ denote the bounded derived category of coherent sheaves over a scheme $Y$, which is always assumed to be separated and of finite type over a field.

The following facts are due to Rouquier \cite{rouquier}. 
\begin{itemize}
\item If $Y$ is a separated scheme of finite type over a perfect field, then $\rdim D^b \op{Coh}(Y)$ is finite; \cite[Thm.\ 7.38]{rouquier}. 
\item If $Y$ is a smooth quasi-projective scheme over a field, then $\rdim D^b \op{Coh}(Y) \leq 2 \op{dim}(Y)$; \cite[Prop.\ 7.9]{rouquier}.
\item If $Y$ is a reduced separated scheme of finite type over a field, then $\rdim D^b \op{Coh}(Y) \geq \op{dim}(Y)$; \cite[Prop\ 7.16]{rouquier}.
\item If $Y$ is a smooth affine scheme of finite type over a field, then $\rdim D^b \op{Coh}(Y) = \op{dim}(Y)$; \cite[Thm.\ 7.17]{rouquier}.
\end{itemize}

\begin{remark}
The idea of Rouquier's proof of the upper bound $\rdim D^b \op{Coh}(Y) \leq 2 \op{dim}(Y)$ above is to resolve the diagonal in $Y \times Y$ and conclude by an analog of \Cref{lemma:resolution-diag}. Such a resolution has $\op{dim}(X \times X)+1= 2 \op{dim}X +1$ terms, which explains the factor of $2$. 
\end{remark}

It follows from the above facts that $\rdim D^b \op{Coh}(Y)= \dim Y$ if $Y$ is affine. Orlov conjectured in \cite[Conj.\ 10]{orlov2009remarks} that this equality holds for any smooth quasi-projective variety. We will return to this conjecture in \Cref{subsection:orlov-conj}.

If $Y$ is a variety which is not assumed to be smooth, then the Rouquier dimension $D^b \op{Coh}(Y)$ can be bigger than its Krull dimension. The following example, which was communicated to us by W.\ Hara, shows that there exist projective curves whose which have arbitrarily large Rouquier dimension.

\ex\label{example:hara}
Fix a positive integer $m$. Wantanabe \cite[Thm.\ 1]{watanabe} constructed a $k$-algebra $(A, \fk{m})$ which is a complete intersection Noetherian local domain of Krull dimension $1$ and with $\op{dim}_{R/\fk{m}} \fk{m}/ \fk{m}^2= m$. Moreover, by construction, $A$ is the localization of a finite type $k$-algebra $B$ at a maximal ideal. 

Let $Y_0 = \op{Spec} B$ and let $Y \sub \bb{P}^n$ be the projective compactification of $Y_0$. We then have $\rdim D^b \op{Coh}(Y) \geq \rdim D^b \op{Coh}(Y_0) \geq \rdim D^bCoh(R)$, where the second inequality is by \cite[Lem.\ 4.2]{aihara-takahashi}. We now appeal to \cite[Cor.\ 5.10]{b-i-k-o}, which implies that $\rdim D^bCoh(R) \geq m- \op{dim} Y -1= m-2$. 
\eex

Ballard--Favero gave examples of singular varieties whose Rouquier dimension equals their Krull dimension (see \cite[Cor.\ 3.3]{ballard2012hochschild}). We will prove in \Cref{subsection:orlov-conj} that such an equality holds for all toric boundary divisors of smooth quasiprojective DM stacks over the complex numbers of dimension at most $3$ and some other singular varieties.

\subsubsection{Representation categories of quivers}

In this paper, a \emph{quiver} $\G$ is a directed graph which is connected, finite, and acyclic (i.e.\ having no loops and no cycles). 

Given a quiver $\G$ and a field $k$, we can consider the \emph{path algebra} $k[\G]$. As a vector space, $k[\G]$ is generated by words of composable arrows; in particular, it is finite dimensional as a $k$-vector space since $\G$ is acyclic. The product of two words is given by concatenation when this makes sense and is defined to be zero otherwise. The path algebra of a quiver is in general non-commutative. Unless otherwise indicated, we always view $k[\G]$ as an \emph{ungraded} algebra (although it admits a natural grading by path length). 

For a quiver $\G$, we will be interested in the $A_\infty$ category $\op{Perf} k[\G]$, as well as in the associated triangulated category $H^0(\op{Perf} k[\G])$. As noted in \Cref{remark:perf}, the later category is equivalent to the smallest idempotent complete triangulated subcategory of $D(k[\G])$ containing $k[\G]$.

A quiver is said to be of Dynkin ADE type if its underlying graph is a Dynkin diagram of type $A_n, D_n, E_6, E_7, E_8$. 

The following proposition is a special case of \cite[Prop.\ 4.3, 4.4]{elagin2020calculating} and is ultimately an application of Gabriel's theorem \cite[Chap.\ 4]{quivers-book}. (We warn the reader that the notation $\op{Perf}(A)$ in \cite{elagin2020calculating} denotes the triangulated category which we have been denoting by $H^0(\op{Perf} A)$.)
\prop\label{proposition:elagin-gabriel}
Let $\G$ be a quiver.
\begin{itemize}
\item If $\G$ is of Dynkin ADE type, then $\rdim (\op{Perf} k[\G])= 0$; 
\item if $\G$ is not of Dynkin ADE type, then $\rdim (\op{Perf} k[\G])=1$.
\end{itemize}
\eprop

\section{Wrapped Fukaya categories}

\subsection{Liouville manifolds and related objects}
\subsubsection{Basic definitions}
A Liouville vector field on a symplectic manifold $(X, \o)$ is a vector field $V$ which satisfies $\mc{L}_V \o=\o$. Dually, a Liouville $1$-form $\l$ on a symplectic manifold $(X, \o)$ is a primitive for $\o$.  

A \emph{Liouville cobordism} is an exact symplectic manifold $(X_0, \l_0)$ with boundary $\d X_0= \d_+X_0 \sqcup \d_- X_0$, such that the Liouville vector field is outward pointing along $\d_+X_0$ and inward pointing along $\d_-X_0$. If $\d_- X_0 = \emptyset$, such a cobordism is called a \emph{Liouville domain}.


A \emph{Liouville manifold} is an exact symplectic manifold $(X, \l)$ which is modeled near infinity on the symplectization of a closed contact manifold $(Y, \xi=\op{ker} \l_Y)$. More formally, this means that there exists a proper embedding 
\eq
e: (SY, \l_Y) \hookrightarrow (X, \l)
\eeq
which covers the complement of a compact set and satisfies $e^*\l=\l_Y$. Here $(SY, \l_Y)= \{\s \in T^*Y \mid \s|_\xi =0\}$ denotes the symplectization of $(Y, \xi)$, where $\l_Y$ is the (restriction of) the tautological $1$-form on $T^*Y$. A Liouville manifold can equivalently be viewed as the result of completing a Liouville domain by gluing on the positive symplectization of its boundary. The contact manifold $Y$ is written as $\partial_{\infty} X$ and is called the \emph{ideal contact boundary} of $X.$

A \emph{Liouville sector} is a Liouville manifold with boundary whose characteristic foliation on the boundary satisfies a niceness condition which ensures that Floer theory is well-behaved (see \cite[Def.\ 2.4]{gps1} for the precise definition).  A Liouville manifold is in particular a Liouville sector. One example of a Liouville sector which is not a Liouville manifold is the cotangent bundle of a compact manifold with boundary. Given a Liouville sector $(X,\l)$, a \emph{stop} is an arbitrary closed subset $\fk{f} \sub (\d_\infty X)^\circ$. Such a pair $(X, \fk{f})$ is often called a stopped Liouville sector.

A Liouville manifold $(X, \l)$ is said to be \emph{Weinstein} if the Liouville vector field $V$ is gradient-like with respect to a proper Morse function $\phi$. Note that according to this definition, being Weinstein is a \emph{property} of Liouville manifolds rather than an extra structure. A Liouville cobordism is said to be Weinstein if its completion is Weinstein; a Liouville sector is said to be Weinstein if its convex completion \cite[Sec.\ 2.7]{gps1} is Weinstein. A \emph{Liouville pair} consists of a pair $(X, A)$ such that $X$ is a Liouville manifold and $A$ (possibly empty) is a Liouville hypersurface embedded inside the ideal contact boundary $\partial_{\infty} X$, i.e. the induced contact form $\lambda_{\partial_{\infty} X}$ from the Liouville $1$--form $\lambda$ restricts to the Liouville $1$--form on $A$. The pair $(X, A)$ is called a \emph{Weinstein pair} of both $X$ and $A$ are Weinstein. A \emph{deformation} of Liouville manifolds/sectors/cobordisms is a one parameter family of such objects. 

\subsubsection{Skeleta} 
The \emph{skeleton} $\fk{c}_X \subset X$ of a Liouville manifold or sector $(X, \l)$ is the set of points which do not escape to infinity under the positive Liouville flow. (More formally: $p \in \fk{c}_X \subset X$ iff there exists a compact subset $K \subset X$ containing $p$ and such that $p$ does not leave $K$ under the positive Liouville flow).  The skeleton is sometimes called the \emph{core} in the literature.

The (relative) skeleton $\fk{c}_{X, \fk{f}}$ of a stopped Liouville manifold $(X, \fk{f})$ is defined to be the set of points which do not escape to $\d_\infty X- \fk{f}$ under the positive Liouville flow. If $(X, A)$ is a Liouville pair, we similarly define the (relative) skeleton $\fk{c}_{X, A}$ to be the set of points which do not escape to $\d_\infty X- \fk{c}_A$ under the positive Liouville flow.

The skeleton of an arbitrary Liouville manifold always has measure zero.  If $(X, \l)$ is Weinstein, then a choice of a gradient-like proper Morse function $\phi: X \to \R$ defines a handlebody decomposition of $X$. The skeleton is then precisely the union of the stable manifolds associated to this decomposition. In particular, the skeleton is stratified by isotropic submanifolds. In contrast, the skeleton of an arbitrary Liouville manifold can be much larger (see e.g.\ McDuff's examples in \cite{mcduff}).
%

\subsection{Invariants of Liouville manifolds}

The \emph{wrapped Fukaya category}  $\mc{W}(X)$ of a Liouville manifold $(X, \l)$ is an $A_\infty$ category. The objects of this category are exact, cylindrical Lagrangian submanifolds. The $A_\infty$ operations are constructed by counting pseudoholomorphic curves. The wrapped Fukaya category was first constructed by Abouzaid--Seidel \cite{abouzaid-seidel} and has become a fundamental invariant in the study of Liouville manifolds.

More generally, one can also consider the wrapped Fukaya category $\mc{W}(X,\fk{f})$ of a stopped Liouville sector. The relevant theory was developed by Ganatra--Pardon--Shende (see \cite[Sec.\ 2]{gps2} for an overview), following earlier work of Sylvan \cite{sylvan-stops}. 

The wrapped Fukaya category is invariant up to quasi-equivalence under deformations of the underlying Liouville manifold/sector. More generally, the wrapped Fukaya category of a stopped Liouville manifold is invariant under deformations which preserve the contactomorphism type of the complement of the stop \cite[Thm.\ 1.4]{gps2}.  

The \emph{symplectic cohomology} $SH^\bu(X)$ of a Liouville manifold $(X,\l)$ is a graded-commutative ring (meaning that $a b= (-1)^{|a| |b|} b a$). The generators can be interpreted as orbits for a Hamiltonian on $X$ which grows suitably fast at infinity, and the differential counts pseudoholomorphic curves. Our grading conventions for symplectic cohomology are uniquely determined by stipulating that the differential increases grading and that the unit lives in degree zero. They match those of e.g.\ \cite{seidel-biased-view} and \cite{ganatra}. 

Symplectic cohomology was originally introduced by Floer and Hofer \cite{floer-hofer} and subsequently recast by Viterbo \cite{viterbo}. We refer to \cite{seidel-biased-view} for a beautiful survey which also matches our sign and grading conventions.

\rmk[Orientations/gradings]
Note that one must fix additional orientations/grading data in order for the above invariants to be $\Z$-graded and defined over a field of characteristic zero. We generally avoid discussing signs/gradings in this paper (see instead \cite[Sec.\ 5.3]{gps3}).
\ermk

\rmk[Technical setups]
There are multiple technical setups for Floer theoretic invariants on Liouville manifolds in the literature. In particular, there are (at least) three distinct versions of the wrapped Fukaya category of a Liouville manifold. The first is the so-called ``linear" version of \cite{abouzaid-seidel}. The second is the so-called ``quadratic" version, which first appeared in \cite{abouzaid-generation}. The third is the so-called ``localization" version which was first implemented in the literature by Ganatra--Pardon--Shende in \cite{gps1}. Some standard references (such as \cite{abbondandolo-schwarz} and \cite{ganatra}) also construct Floer theoretic invariants using almost-complex structures of ``rescaled contact type", which is in contrast to the previously cited references which only use almost-complex structures of contact type. 

The fact that all of these setups yield isomorphic theories is a folklore result which is routinely assumed in the literature. However, to the best of our knowledge, detailed proofs for most of these equivalences have yet to appear (see however \cite[Prop.\ 2.6]{sylvan} concerning the comparison between the quadratic and localization versions). In this paper, the reader is free to assume throughout that the wrapped Fukaya category is defined according to the localization setup of \cite{gps1} (this is the only setup available for sectors). However, we will quote facts from \cite{abbondandolo-schwarz} and \cite{ganatra} in \Cref{subsection:string-top}, and we will sometimes also appeal to certain mirror symmetry results from the literature which are proved using other setups. 
\ermk

\subsection{Structural facts about wrapped Fukaya categories}

We collect some structural properties of wrapped Fukaya categories which will be used in the sequel. Unless otherwise indicated, the theory summarized here is due to Ganatra--Pardon--Shende \cite{gps1, gps2}. 

\defi[Def.\ 1.7 in \cite{gps2}]
A closed subset $\fk{c}$ of a symplectic manifold $W$ is said to be \emph{mostly Lagrangian} if it admits a decomposition $\fk{c}= \fk{c}^{\op{crit}} \cup \fk{c}^{subcrit}$, where $\fk{c}^{crit}$ is a Lagrangian submanifold (in general non-compact) and $\fk{c}^{\op{subcrit}}$ is closed and covered by the smooth image of a second countable manifold of dimension at most $\frac{1}{2} \op{dim} W-1$. We similarly define the notion of a \emph{mostly Legendrian} subset of a contact manifold.
\edefi

If $\fk{c}_X \subset X$ is mostly Lagrangian, a \emph{generalized cocore} associated to a connected component of $\fk{c}_X^{\op{crit}}$ is an exact (cylindrical at infinity) Lagrangian which intersects this component transversally in a single point.  There is a closely related notion for stops: suppose $\fk{c} \sub \d_\infty X$ is a stop and let $p \in \fk{c}^{\op{crit}}$. Then there is an associated Lagrangian submanifold $L_p \sub X$ called the \emph{linking disk}, which is well-defined up to isotopy and depends only on the connected component of $\fk{c}^{\op{crit}}$ containing $p$. 

\ex
If the cocores of $(X, \l)$ are properly embedded (this condition holds after a generic perturbation of $\l$), then the skeleton is mostly Lagrangian. Similarly, the skeleton of a Weinstein hypersurface $V \sub \d_\infty  X$ is mostly Legendrian after possibly perturbing $V$. See \cite[Sec.\ 1.3]{gps2}. 
\eex

If $X$ is a Liouville sector and $\fk{c} \sub \d_\infty X$ is a stop, then there is a natural functor (often called ``pushforward" or ``stop-removal")
\eq\label{equation:stop-removal}
\mc{W}(X, \fk{c}) \to \mc{W}(X).
\eeq
This functor is characterized as follows when $\fk{c}$ is mostly Legendrian.

\thm[see Thm.\ 1.20 in \cite{gps2}]\label{theorem:stop-removal}
If $\fk{c} \sub \d_\infty X$ is mostly Legendrian, then \eqref{equation:stop-removal} induces a quasi-equivalence
\eq 
\mc{W}(X, \fk{c})/ \mc{D} \xrightarrow{\sim} \mc{W}(X),
\eeq
where $\mc{D}$ is the collection of linking disks associated to $\fk{c}^{\op{crit}}$. 
\ethm

It is a fundamental fact due independently to Chaintraine--Dimitroglou Rizell--Gigghini--Golovko \cite{cdgg} and Ganatra--Pardon--Shende \cite[Thm.\ 1.13]{gps2} that the wrapped Fukaya category of a Weinstein manifold is generated by cocores. We will need a slightly more refined statement, which follows from essentially the same proof.  This is the content of \Cref{theorem:gps-intersections-enhanced} below. 

To set the stage, let $(X, \l)$ be a Liouville manifold. Suppose that the skeleton $\fk{c}_X=\fk{c}^{\op{crit}} \cup \fk{c}^{\op{subcrit}}$ is mostly Lagrangian, and that every connected component of $\fk{c}^{\op{crit}}$ admits a generalized cocore.  Let $L$ be an object of $\mc{W}(X)$ such that the underlying Lagrangian satisfies the following conditions:
\begin{itemize}
\item $L \cap \fk{c}_V \sub \fk{c}_V^{\op{crit}}$;
\item $L$ intersects $\fk{c}_V^{\op{crit}}$ transversally.
\end{itemize}

\thm \label{theorem:gps-intersections-enhanced}
With the above assumptions, suppose that $|L \cap \fk{c}_V| = r \in \N$. Then there exists a resolution
\eq 0=K_0 \to K_1 \to \dots \to K_{r-1} \to K_r=L \eeq
where each $\op{cone}(K_{i-1} \to K_i)$ is isomorphic in $H^0(\mc{W}(X))$ to a generalized cocore.
\ethm

\pf
The proof is essentially the same as that of \cite[Thm.\ 1.13]{gps2}. After possibly modifying $L$ outside a compact set, we may assume that $\l|_L$ has a compactly-supported primitive \cite[Lem.\ 7.2]{gps2}. Let $f: X \to \R$ be an extension of this primitive to $X$.

Consider the functors
\eq\label{equation:k-embedding}
\mc{W}(X) \hookrightarrow \mc{W}(X \tms \C, \fk{c}_X \tms \{ \pm \infty\}) \xrightarrow{\sim}  \mc{W}(X \tms \C, \fk{c}_X \tms \{+\infty\} \cup \fk{c}_X \tms \{e^{4\pi i /3} \cdot \infty\}),
\eeq
where the first arrow is the K\"{u}nneth embedding \cite[Thm.\ 1.5]{gps2} and the second arrow is induced by the isotopy of stops $\fk{c}_X \tms \{-\infty\} \rightsquigarrow \fk{c}_X \tms \{e^{\pi i \cdot (1+s)} \infty\})$ for $s \in [0,1/3]$, which is a quasi-equivalence by \cite[Thm.\ 1.4]{gps2}.

On objects, \eqref{equation:k-embedding} sends $L \mapsto L \tilde{\tms} i\R$, which is the \emph{cylindrization} of the product Lagrangian $L \tms i\R$ as defined in \cite[Sec.\ 7.2]{gps2}. The cylindrization is well-defined up to a contractible choice of Lagrangian isotopies. It will be convenient to construct a nice cylindrization. To this end, let $C \gg 1$ and $0<\e \ll 1$ be constants which will be fixed later. Let $\phi(r, \t)= \phi(r): \C \to \R$ be a function with the following properties: 
\begin{itemize}
\item the derivative of $\phi$ in the radial direction is non-negative;
\item $\phi(r)=0$ for $r \in [0, C-\e]$ and $\phi(r)=1$ if $r \geq C$.
\end{itemize}

To define the cylindrization, we consider the Liouville $1$-form $\l_X + \l_{\C} - d(\phi f)$, where $\l_{\C}= r^2 d\t$. The associated Liouville vector field is 
\eq\label{equation:liouville-deformed} 
(Z_X + \phi X_f) + (\frac{1}{2}r\d_r + f X_{\phi}).
\eeq 
Assuming $C$ is large enough, it was shown in \cite[Sec.\ 7.2]{gps2} that there is a unique Hamiltonian isotopy $\Phi$ which is the identity on a Liouville subdomain of $X \tms \C$ and sends $Z_X + \frac{1}{2}r\d_r$ (i.e.\ the Liouville vector field of $\l_X+ \l_{\C}$) to \eqref{equation:liouville-deformed}. 

For $t \in [0, 2/3]$, let $L_t:= L \tms e^{\pi i \cdot (t + 1/2)} \R$.  We define $\tilde{L}_t := \Phi(L_t)$ to be the cylindrization.

Consider now the contact manifold $\mc{O}= X \tms \{ r= C\}$, with contact form $\l_X+ C^2 d\t$. This defines coordinates on (a region of) $(\d_\infty(X \tms \C), \l_X + \l_\C)$.  We set $\L_t = \mc{O} \cap L_t$ and $\tilde{\L}_t = \mc{O} \cap \tilde{L}_t$. By construction, $\tilde{\L}_t \sub \mc{O}$ is Legendrian. Moreover, by taking $\e$ small enough, we can ensure that $\tilde{L}_t$ is arbitrarily close to the Legendrian lift of $L_t$ in any $C^k$ norm. In particular, the projection of $\tilde{L}_t$ onto $X \tms \{r=C, \t=\pi \cdot (t + 1/2)\}$ is arbitrarily close to $\L_t$. 

We now consider the isotopy $\tilde{L}_t$ for $t \in [0, 2/3]$. The intersection points of $\tilde{L}_t$ with the stop $\fk{c}_X \tms \{+\infty\}$ are in correspondence with the intersection points of $\L_t$ with $\fk{c}_X \tms \{\infty\}$, and hence with the intersection points $L \cap \fk{c}_X$. By assumption, these intersections are all transverse and contained in $\fk{c}_X^{\op{crit}}$. 

The rest of the argument is identical to the proof of \cite[Thm.\ 1.13]{gps2} and is therefore only sketched. We may assume (e.g.\ by a small perturbation of the stop) that the intersections happens at discrete times $0<t_1<\dots <t_r <2/3$. At each intersection point $t_i$, the wrapping exact triangle \cite[Thm.\ 1.10]{gps2} gives $L_{t_i+\delta} \to L_{t_i-\delta} \to D_i \xrightarrow{[1]}$, where $D_i$ is a linking disk. One can verify that $L_{t_r+\delta}$ is the zero object (for example, because it is contained in a Liouville sector deformation equivalent to $X \tms \C_{\op{Im} \geq 0}$ whose wrapped Fukaya category thus vanishes). Finally, it can be shown that the linking disks correspond to generalized cocores under \eqref{equation:k-embedding}. 
\epf

\rmk
We expect that a version of \Cref{theorem:gps-intersections-enhanced} can also be established using the methods of \cite{cdgg} (see in particular \cite[Prop.\ 9.3]{cdgg}). 
\ermk

\subsection{The Rouquier dimension of wrapped Fukaya categories}

We now arrive at the main topic of this paper.

\defi
Let $(X, \l)$ be a Liouville manifold. We let 
\eq \rdim{\mc{W}(X)} \in \Z_{\geq -1} \cup \{\infty \} \eeq
denote the Rouquier dimension of the wrapped Fukaya category $\mc{W}(X)$, as defined in \Cref{definition:rdim}. 

More generally, if $(X, \fk{f})$ is a stopped Liouville sector, we let $\rdim{\mc{W}(X,\fk{f})}$ denote the Rouquier dimension of its wrapped Fukaya category. 
\edefi

\rmk[Orientations/gradings]
Of course, the category $\mc{W}(X)$ depends on auxiliary orientations/grading data, so $\rdim{\mc{W}(X)}$ also depends on this data.  Thus $\rdim \mc{W}(-)$ should be viewed as an invariant of Liouville manifolds \emph{equipped with orientation/grading data}. 

To obtain an honest invariant of Liouville manifolds defined in maximum generality, one could simply choose to work exclusively with $\Z/2$-coefficients and $\Z/2$-gradings. On the other hand, for the purpose of this paper, we are mostly interested in the Rouquier dimension as a tool for symplectic and algebro-geometric applications. Thus we will typically restrict our discussion to whichever choice of coefficients/gradings is most convenient for computations or applications. 
\ermk

In the following examples, we collect some preliminary facts about the Rouquier dimension of wrapped Fukaya categories. 

\ex
If $X$ is Weinstein, then it can be shown using either \cite[Thm.\ 1.1]{cdgg} and \cite[Thm.\ 1.13]{gps2} that the diagonal $\Delta \sub \mc{W}(\ov{X} \tms X, -\l \oplus \l)$ is generated by product-cocores. It then follows by combining \Cref{theorem:bimod-functor} in the appendix with \Cref{lemma:resolution-diag} that 
\eq \rdim{\mc{W}(X)} < \infty. \eeq
We do not know if there exists a Liouville manifold whose wrapped Fukaya category has infinite Rouquier dimension.  
\eex

\ex
If $\mc{W}(X)=0$, then $\rdim{\mc{W}(X)}=-1$. This holds of course for subcritical Weinstein manifolds (and more generally for Weinstein manifolds which admit a handlebody presentation where all critical handles are attached along loose Legendians).
\eex

\ex[Point-like objects]\label{example:point-like-objects}
Suppose that $\mc{W}(X)$ is equipped with $\Z$-grading and admits a point-like object (e.g.\ an exact Lagrangian torus). Then it follows from \cite[Thm.\ 5.2]{elagin-lunts-koszul} that $\rdim \mc{W}(X) \geq \frac{1}{2} \op{dim}_{\R} X$. (For an $A_\infty$ category $\mc{C}$, we say that $L \in \mc{C}$ is \emph{point-like} if (i) $\mc{C}$-module $\mc{Y}^r_L$ is proper and; (ii) $\hom(L, L)= \Lambda^\bu V$ as $A_\infty$ algebras, where the right hand side is the exterior algebra (viewed as an $A_\infty$ algebra) of an $r$-dimensional graded vector space $V$ which is supported in degree $1$.  Strictly speaking, \cite[Thm.\ 5.2]{elagin-lunts-koszul} applies only to proper categories. However, one can verify that the proof only uses the weaker assumption (i).)
\eex

\ex[Homological mirror symmetry]
If $X$ happens to be mirror to an algebraic variety $X^\vee$, then we can immediately transfer the computations of $\rdim D^b \op{Coh}(X^\vee)$ discussed in \Cref{subsection:first-examples} into computations of $\rdim \mc{W}(X)$. A partial list of homological mirror symmetry results for Weinstein manifolds and sectors is \cites{hacking2020homological, lekili-polishchuk, lekili-ueda, gammage-bh, gammage2017mirror, pomerleano2021intrinsic}. We leave it to the reader to verify which results from \Cref{subsection:first-examples} apply in each of these instances.
\eex

\section{Symplectic flexibility and Orlov's conjecture}\label{sec_arbo}
Before delving into the technical details, let us illustrate the ideas of this section by considering the example $(X, \l) = (T^*M, pdq)$, where $M$ is an $n$-dimensional closed manifold. Choose a smooth triangulation (see the discussion in Section \ref{subsec_triang} for a precise definition) $\mc{S}$ of $M$ and denote by $\{ v_{\a} \}_{\a \in I}$ the set of $0$-simplices of $\mc{S}$. Denote by $\op{star}(v_{\a})$ the star of the vertex $v_{\alpha}$ and let $\ov{\op{star}(v_{\a})}$ be its closure. Then $\ov{\op{star}(v_{\a})}$ is an embedded polyhedron. We choose an inward smoothing of each $\ov{\op{star}(v_{\a})}$ to obtain a smooth manifold with boundary $M_{\a}$ such that the boundaries $\{ \partial M_{\a} \}_{\a \in I}$ intersect transversely. By \cite[Ex.\ 1.33]{gps2}, the collection of cotangent bundles $\{ T^* M_{\a} \}_{\a \in I}$ forms a \emph{sectorial cover} (see \cite[Def.\ 1.32]{gps2}) of the Weinstein manifold $(T^*M, pdq)$. Then \cite[Thm.\ 1.35]{gps2} asserts a pre-triangulated equivalence
\eq\label{eqn:cotan_descent}
\op{hocolim}_{\emptyset \neq J \subset I} \mc{W}(\bigcap_{\a \in J} T^* M_{\a}) \xrightarrow{\sim} \mc{W}(T^*M).
\eeq
The sectorial cover $\{ T^* M_{\a} \}_{\a \in I}$ admits the following special properties:
\begin{enumerate}
\item it has depth at most $n+1$ if the inward smoothing is small enough; i.e.\ if $\bigcap_{\a \in J} T^* M_{\a} \neq \emptyset$, then $|J| \leq n+1$;
\item each local piece $\bigcap_{\a \in J} T^* M_{\a}$ is Weinstein homotopic to the cotangent bundle of a closed $n$-dimensional ball $\ov{D}^n$.
\end{enumerate}
Note that $\mc{W}(T^* \ov{D}^n)$ is equivalent to $\op{Perf}(k[A_1])$ where $A_1$ is the quiver with only one vertex so $\rdim \mc{W}(T^* \ov{D}^n) = 0$ by Proposition \ref{proposition:elagin-gabriel}. Then Lemma \ref{lem_rdim_upper} and \ref{lem_rdim_max} and the fact that this cover has depth at most $n+1$ implies $\rdim \mc{W}(T^* M) \leq n$ because of \eqref{eqn:cotan_descent}.

We wish to apply the above ideas to general Weinstein manifolds or sectors with arboreal skeleta (see Section \ref{subsec_arbo}). For purely technical reasons, it is more convenient to implement this argument using microlocal sheaf theory. We therefore study the \emph{Nadler--Shende microlocal sheaf} (see Section \ref{subsec_nadler-shende}) as an alternative model of the wrapped Fukaya category, and appeal to \cite[Thm.\ 1.4]{gps3} to translate our results back into the framework of wrapped Fukaya categories. The Nadler--Shende sheaf satisfies a local-to-global descent formula similar to \eqref{eqn:cotan_descent}. Thus we can cover the skeleta by the stars of topological triangulations and conclude our result Theorem \ref{thm_arbo_sheaf_upper} as above.

We will use the notations $\op{skel}(X)$ and $\fk{c}_{X}$ interchangeably throughout this section to denote the skeleton of $(X, \l)$.

\subsection{Stratifications and Goresky's triangulation theorem}\label{subsec_triang}
A \emph{stratification} $\mc{S}$ of a topological space $X$ is a locally finite decomposition $\{X_\s\}_{\s \in \mc{S}}$ of $X$ into disjoint, locally closed subsets (called strata\footnote{The strata are required to be connected under some conventions. We find it more convenient to also allow disconnected strata -- however, this is just a matter of terminology which does not affect the substances of the arguments.}) which satisfies the frontier condition: if $\ov{X_\a} \cap X_\b \neq \emptyset$, then $X_\b \sub X_\a$. The pair $(X, \mc{S})$ is said to be a \emph{stratified space}. Note that $\mc{S}$ is naturally a partially ordered set, where we set $X_\a \leq X_\b$ if and only if $X_\a \sub \ov{X_\b}$. 

Given a stratum $X_\a$, the (open) star of $X_\a$ is defined to be the union of strata from $\mc{S}$ whose closures contain $X_\a$ and is denoted by $\op{star}(\a)$. The \emph{link} of $X_\a$ is denoted by $\op{link}(\a)$ and defined as the closure of $\op{star}(\a)$ minus its interior. 

Given a fixed topological space $X$, a stratification $\mc{S}'$ of $X$ is said to be a \emph{refinement} of some coarser stratification $\mc{S}$ of $X$ if every stratum in $\mc{S}'$ is contained in a stratum in $\mc{S}$ (equivalently, every stratum in $\mc{S}$ is a union of strata in $\mc{S}'$).

If $X$ is a topological subspace of a manifold (without loss of generality $X \sub \R^N$), then a stratification $\mc{S}$ of $X$ is said to be $C^\mu$ (for $\mu \in \mathbb{Z}_{\geq 0} \cup \{ \infty \}$) if each stratum is a locally closed $C^\mu$ submanifold.

\defi
A $C^1$ stratification $\mc{S}$ of $X \subset \R^N$ is called \emph{Whitney} if the following condition is satisfied.

\emph{Whitney (b)}: given any pair of strata $X_{\a}$ and $X_{\b}$ for $\a, \b \in \mc{S}$ with $X_\a < X_\b$, suppose that $x_i \in X_{\b}$ and $y_i \in X_{\a}$ are a sequences of points both converging to $y \in X_{\a}$. Suppose also that $T_{x_i} X_{\b}$ converges to a subspace $T \sub T_y\R^N$, and the secant lines $\ov{x_iy_i}$ also converge to a line $ \ell \sub T_y \R^N$. Then $\ell \sub T$.
\edefi

Next we review the notion of a triangulation. A $k$-dimensional \emph{simplex} is a subset of $\R^d$ (for some $d>0$) which is the convex hull of some finite set of independent points $p_0,\dots, p_k$.  The \emph{faces} of a simplex are defined in the obvious way. A \emph{simplicial complex} (see \cite[Def.\ 7.1]{munkres2016elementary}) is a set $\mc{K}$ of simplices in $\R^d$ (for some fixed $d>0$) which satisfies the following conditions: 
\begin{enumerate}
\item every face of a simplex $\s \in \mc{K}$ is also in $\mc{K}$;
\item the non-empty intersection of any two simplices $\s_1, \s_2 \in \mc{K}$ is a face of both $\s_1$ and $\s_2$;
\item local finiteness (i.e.\ each point of $|\mc{K}|$ is contained in finitely many simplices, where $|\mc{K}|$ is the topological space formed by the union of the simplices from $\mc{K}$). 
\end{enumerate}

Note that the space $|\mc{K}|$ admits a stratification whose strata are given by interiors of simplices from $\mc{K}$ and we denote this stratification by $\mc{S}_{\mc{K}}$. Each stratum of $\mc{S}_{\mc{K}}$ is naturally a smooth submanifold of $\R^d$ and it is easy to see that $\mc{S}_{\mc{K}}$ defines a Whitney stratification. 

The notion of a simplicial subcomplex is defined in the obvious way. 
\defi
Given a smooth manifold $M$, a \emph{smooth triangulation} of $M$ consists of the following data:
\begin{enumerate}
\item a simplicial complex $\mc{K}$ and a simplicial subcomplex $\mc{L} \subset \mc{K}$;
\item a homeomorphism $f: |\mc{K}| \setminus |\mc{L}| \to X$ such that for each simplex $\s \in \mc{K}$ and $x \in \s \setminus |\mc{L}|$, there is a neighborhood $U$ of $x$ in the plane containing $\s$ and a smooth embedding $\tilde{f}: U \to X$ extending $f|_{U \cap \sigma}$.
\end{enumerate}
\edefi
\defi\label{def_triangu}
For a Whitney stratified space $(X, \mc{S})$ with smooth strata, a triangulation of $(X, \mc{S})$ is given by a simplicial complex $\mc{K}$ and a homeomorphism $f: |\mc{K}| \to X$ such that for each stratum $X_{\alpha}$, the set $f^{-1}(X_{\alpha})$ is the underlying topological space of some simplicial subcomplex of $\mc{K}$ and $f: f^{-1}(X_{\alpha}) \rightarrow X_{\alpha}$ is a smooth triangulation.
\edefi
Given a triangulation of $(X, \mc{S})$ as in Definition \ref{def_triangu}, the stratification $\mc{S}_{\mc{K}}$ on $|\mc{K}|$ defines a (topological) stratification of $X$ and we call a stratification of this from \emph{triangulation}. It is easy to see that this stratification is a refinement of the stratification $\mc{S}$. The following theorem of Goresky asserts that it is always possible to refine a Whitney stratification into a triangulation.
\thm[Goresky \cite{goresky1978triangulation}]\label{theorem:goresky-triangulation}
Every Whitney stratified space $(X, \mc{S})$ with smooth strata admits a triangulation.
\ethm
Note that Theorem \ref{theorem:goresky-triangulation} does not prove that the stratification induced by the triangulation is Whitney. It is an open question whether the theorem can be strengthened in that way. 

For the purpose of our application, we introduce some more notions. Given two stratified topological spaces $(X, \mc{S})$ and $(X', \mc{S}')$, a \emph{stratified homeomorphism} $(X, \mc{S}) \to (X', \mc{S}')$ is a homeomorphism $X \to X'$ which sends strata to strata and defines a bijection between strata.
\defi
For a stratified topological space $(X, \mc{S})$, a \emph{topological triangulation} of $(X, \mc{S})$ is defined by a simplicial complex $\mc{K}$ and a homeomorphism $f: |\mc{K}| \to X$ such that $\mc{S}$ admits a refinement $\mc{S}'$ with the property that $f: (|\mc{K}|, \mc{S}_{\mc{K}}) \to (X, \mc{S}')$ defines a stratified homeomorphism.
\edefi
Then Theorem \ref{theorem:goresky-triangulation} implies the following result which suffices for our application.
\cor\label{cor_triang_exist}
Any Whitney stratified space with smooth strata admits a topological triangulation. \qed
\ecor

\subsection{Arboreal singularities}\label{subsec_arbo}
We review the necessary inputs from the arborealization program \cite{AGEN1, AGEN2, AGEN3} and develop some topological aspects of arboreal Lagrangians needed for our applications.

\subsubsection{Arborealization of polarized Weinstein manifolds}
A \emph{signed rooted tree} $\fT = (T, \rho, \epsilon)$ consists of an acyclic finite graph $T$, a distinguished vertex $\rho$ which is called the \emph{root}, and a decoration $\epsilon$ of edges which are not adjacent to $\rho$ by signs $\pm$. Let $n(\fT)$ be the number of vertices of $T$ which are not the root. To any such $\fT$, \cite[Def.\ 2.19]{AGEN1} associates a closed Lagrangian subset 
\eq L_{\fT} \subset T^*{\R^{n(\fT)}} \eeq 
stratified by isotropic submanifolds, where $T^*{\R^{n(\fT)}}$ is endowed with the standard symplectic structure.

\defi[Def.\ 3.1 in \cite{AGEN1}]\label{definition:arboreal-lagrangian}
A subset $L$ in a symplectic manifold $(M^{2n}, \o)$ is called a \emph{(boundaryless) arboreal Lagrangian} if, for any $p \in L$, there is a neighborhood $U$ containing $p$ and a symplectic embedding
\eq\label{eqn_arbo_model} (U, U \cap L) \hookrightarrow (T^*\R^{n(\fk{T})} \tms T^*\R^{n-n(\fk{T})}, L_{\fk{T}} \tms \R^{n-n(\fk{T})}),\eeq
for some rooted tree $\fk{T}$. Here $\R^{n-n(\fk{T})} \subset T^*\R^{n-n(\fk{T})}$ is the $0$-section.
We say that $L$ is an \emph{arboreal Lagrangian with boundary} (or just an \emph{arboreal Lagrangian}) if for any $q \in L$, there is a neighborhood $U$ containing $q$ and a symplectic embedding modeled either on \eqref{eqn_arbo_model} or
\eq\label{eqn_arbo_bound_model} (U, U \cap L) \hookrightarrow (T^*\R^{n(\fk{T})} \tms T^*\R^{n-n(\fk{T})-1} \tms T^*\R_{\geq 0}, L_{\mc{T}} \tms \R^{n-n(\fk{T})-1} \tms \R_{\geq 0}).\eeq
\edefi

The following theorem is due to \'{A}lvarez-Gavela--Eliashberg--Nadler (see Thm.\ 1.5 and Rmk.\ 1.6(iii) in \cite{AGEN3}).
\thm[Arborealization]\label{theorem:arborealization}
Let $(X, \l)$ be a Weinstein manifold which has a polarization. Let $A \sub \d_\infty X$ be a Weinstein hypersurface. Then there exists the following data:
\begin{enumerate}
\item\label{item:thm_arbo} A family of Liouville forms $(\lambda_t)_{t \in [0,1]}$ on $X$ with $\l_0=\l$ such that $(X, \lambda_t)$ and $(A, {\lambda_t}|_A)$ are both Weinstein homotopies, and such that $\op{skel}((X, \lambda_1), A)$ is arboreal with smooth boundary and possibly with ``corner-type" singularities (i.e.\ locally modeled on $(0_{\R^n} \cup conormal(\{x_n=0, x_{n-1}\geq 0\})) \sub T^*\R^n$); 
\item\label{item:cap-infinity} A Weinstein hypersurface $B \sub \d_\infty X$ such that $\op{skel}(B) \setminus \op{skel}(A, \l_1)$ is a smooth Legendrian with finitely many connected components, and such that $\op{skel}((X, \lambda_1), B)$ is arboreal. (To be clear, we do not allow any ``corner-type" singularities here).  
\end{enumerate}
\ethm

\rmk
When $A = \emptyset$, the skeleton $\op{skel}(X, \lambda_1) = \op{skel}((X, \lambda_1), A)$ from \eqref{item:thm_arbo} of \Cref{theorem:arborealization} can be made into an arboreal Lagrangian with \emph{smooth boundary}.
\ermk

\rmk
If $X$ has real dimension at most $4$, then both conclusions of \Cref{theorem:arborealization} hold without assuming the existence of a global Lagrangian plane field. (This is essentially a consequence of work of Starkston \cite{starkston2018arboreal}, although strictly speaking the result there is only stated under the assumption that $A= \emptyset$, and produces an arboreal skeleton with ``corners".)
\ermk

\subsubsection{The singularity-order stratification}
We introduce a stratification on arboreal Lagrangians according to the local behavior of the singularities.
\defi
Suppose $L \subset (M^{2n}, \omega)$ is an arboreal Lagrangian. Given a rooted tree $\fk{T}$ with $n(\fk{T}) \leq n$, define $L(\fk{T})$ to be the subset of $L$ such that for any $p \in L(\fk{T})$,  there is a neighborhood $U$ containing $p$ and a symplectic embedding modeled on \eqref{eqn_arbo_model}.
For any rooted tree $\fk{T}$ with $n(\fk{T}) \leq n-1$, define $L^{\prime}(\fk{T})$ to be the subset of $L$ such that for any $q \in L^{\prime}(\fk{T})$, there is a neighborhood $U$ containing $q$ and a symplectic embedding modeled on \eqref{eqn_arbo_bound_model}.

Then the following decomposition, which is denoted by $\mc{S}_{\op{arb}}$,
\eq 
L = \bigcup_{\fk{T}, n(\fk{T}) \leq n} L(\fk{T}) \cup \bigcup_{\fk{T}, n(\fk{T}) \leq n-1} L^{\prime}(\fk{T})
\eeq
defines a $C^{\infty}$ stratification of $L$ called \emph{the singularity-order stratification}.
\edefi

\prop\label{proposition:whitney-sing-order}
For any arboreal Lagrangian $L \subset (M, \omega)$, the singularity-order stratification $\mc{S}_{\op{arb}}$ on $L$ is Whitney.
\eprop

The rest of this subsection is devoted to the proof of Proposition \ref{proposition:whitney-sing-order}. 

Suppose we have a possibly singular Legendrian $\Lambda \subset S^* \R^n$ endowed with a $C^{\infty}$ Whitney stratification $\mc{S}_{\Lambda}$ by isotropic submanifolds. When $\Lambda$ is connected, assume that (i) the front projection restricted to each stratum from $\mc{S}_{\Lambda}$ defines an embedding into the $0$-section $\R^n \subset T^* \R^n$; (ii) different strata have disjoint images under the projection $\pi$. If $\Lambda$ is disconnected, we require that (i) each component satisfies the previous conditions and; (ii) given connected components $\Lambda_1, \dots, \Lambda_k$ together with strata $Y_1 \subset \Lambda_1, \dots, Y_k \subset \Lambda_k$, the restriction of the front projection $\pi: S^* \R^n \rightarrow \R^n$ along $Y_1 \cup \cdots \cup Y_k$ is self-transverse.

Then we can consider the singular Lagrangian $\op{Cone}(\Lambda) \cup \R^n \subset T^* \R^n$ (where $\op{Cone}(\Lambda)$ is the Liouville cone and we write $\R^n \equiv 0_{\R^n}$). The Lagrangian $\op{Cone}(\Lambda) \cup \R^n$ admits a natural stratification $\mc{S}_{\op{Cone}(\Lambda) \cup \R^n}$ as follows. Given any singular isotropic $\Lambda'$, note that $\pi(\Lambda') \subset \op{Cone}(\Lambda')$ and define $C^{\circ}(\Lambda') := \op{Cone}(\Lambda') \setminus \pi(\Lambda')$. Then the union of all $C^{\circ}(Y)$ for $Y \in \mc{S}_{\Lambda}$ defines a decomposition of $C^{\circ}(\Lambda)$ and we declare each $C^{\circ}(Y)$ to be a stratum from $\mc{S}_{\op{Cone}(\Lambda) \cup \R^n}$. For the $0$-section $\R^n$, define
\eq (\op{Cone}(\Lambda) \cup \R^n)_{Y_1, \dots, Y_k} = \bigcap_{i=1}^k \pi(Y_i) \setminus \lt( \bigcup_{Y'} \pi(Y')  \cap \bigcap_{i=1}^k  \pi(Y_i) \rt) \eeq
where $Y_1 \subset \Lambda_1, \dots, Y_k \subset \Lambda_k$ are strata of distinct components $\Lambda_1, \dots, \Lambda_k$ and $Y'$ is any stratum of $\Lambda$ different from $Y_1, \dots, Y_k$. Finally, define $(\R^n)_{\emptyset}$ to be the complement of the union of all $(\op{Cone}(\Lambda) \cup \R^n)_{Y_1, \dots, Y_k}$ in the $0$-section. Then we obtain a decomposition
\eq \R^n = (\R^n)_{\emptyset} \cup \bigcup_{Y_1, \dots, Y_k} (\op{Cone}(\Lambda) \cup \R^n)_{Y_1, \dots, Y_k} \eeq
and the decomposition $\mc{S}_{\op{Cone}(\Lambda) \cup \R^n}$ is simply
\eq \op{Cone}(\Lambda) \cup \R^n = \bigcup_{Y} C^{\circ}(Y) \cup  (\R^n)_{\emptyset} \cup \bigcup_{Y_1, \dots, Y_k} (\op{Cone}(\Lambda) \cup \R^n)_{Y_1, \dots, Y_k}. \eeq
$\mc{S}_{\op{Cone}(\Lambda) \cup \R^n}$ indeed defines a stratification as the frontier condition holds.

\lem\label{lem_cone_whit}
The stratification $\mc{S}_{\op{Cone}(\Lambda) \cup \R^n}$ is Whitney.
\elem
\pf
It suffices to check condition \emph{Whitney (b)} for any pair of strata $X_{\a}$ and $X_{\b}$ with $X_\a < X_\b$ from $\mc{S}_{\op{Cone}(\Lambda) \cup \R^n}$. 

If both $X_{\a}$ and $X_{\b}$ are strata of the form $C^{\circ}(Y)$, \emph{Whitney (b)} holds because $\mc{S}_{\Lambda}$ is a Whitney stratification.

If $X_{\a} \subset \R^n$ and $X_{\b}$ is of the form $C^{\circ}(Y)$, then $X_{\a}$ is necessarily of the form $(\op{Cone}(\Lambda) \cup \R^n)_{Y_1, \dots, Y_k}$ such that there exists $Y_i$ with $Y_i \subset \ov{Y}$. \emph{Whitney (b)} holds by observing that the canonical stratification on $[0, \infty)$ by $\{ \{ 0 \}, (0, \infty) \}$ is Whitney.

If both $X_{\a}$ and $X_{\b}$ are contained in the $0$-section $\R^n$, there are two cases. When $\Lambda$ only has one connected component, the projection $\pi$ is assumed to be a topological embedding along $\Lambda$ and is a smooth immersion when restricted to each stratum. Then \emph{Whitney (b)} holds because it holds for the stratification $\mc{S}_{\Lambda}$. If $\Lambda$ has multiple connected components, we make the following observation. Given two smooth submanifolds $N_1$ and $N_2$ of a smooth manifold $M$, if $N_1$ and $N_2$ intersect transversely, then the decompositions $N_1 \cup N_2 = (N_1 \setminus N_2) \cup (N_1 \cap N_2) \cup (N_2 \setminus N_1)$ and $M = (M \setminus N_1 \cup N_2) \cup (N_1 \setminus N_2) \cup (N_1 \cap N_2) \cup (N_2 \setminus N_1)$ are all Whitney stratifications. The obvious generalizations for multiple submanifolds are easy to verify. Then \emph{Whitney (b)} for the case when $\Lambda$ has multiple connected components just follows from this observation, the discussion for one single component and the transversality assumption of the front projection $\pi$ for strata from different connected components. (Note that it is possible that the front projections of components of $\Lambda$ are disjoint from each other, but \emph{Whitney (b)} holds nevertheless similar to the discussion of $\Lambda$ with only $1$ component.)
\epf

\pf[Proof of Proposition \ref{proposition:whitney-sing-order}]
This is because of Lemma \ref{lem_cone_whit} and the axiomatic characterization of arboreal Lagrangians \cite[Def.\ 1.1]{AGEN1}. Note that the stratification constructed in Lemma \ref{lem_cone_whit} coincides with the singularity-order stratification using the inductive description of arboreal Lagrangians from \cite[Def.\ 1.1]{AGEN1}.
\epf

\subsection{Upper bounds on the Rouquier dimensions}\label{subsec_nadler-shende}
Now we combine microlocal sheaf theory computations with the main result of \cite{gps3} to obtain upper bounds for the Rouquier dimension of wrapped Fukaya categories.

\subsubsection{Microlocal sheaves}
Given a topological space $X$, one can consider sheaves/cosheaves valued in a symmetric monoidal stable $\infty$-category on $X$. We refer to \cite[Sec.\ 3.1]{shende-takeda}, \cite[Rmk.\ 5.1]{nadler-shende} for references to foundational material on sheaves of categories.

If $(X, \mc{S})$ is a stratified space, then a sheaf of (stable $\infty$-)categories on $X$ is said to be \emph{constructible} with respect to $\mc{S}$ if its restriction to each stratum is locally constant. Thus constructibility with respect to $\mc{S}$ is a \emph{property}. Under (rather weak) hypotheses on $(X, \mc{S})$, it is a theorem that the full subcategory of $\mc{S}$-constructible sheaves is equivalent to a certain representation category of the ``infinity exit-path category"; see \cite[Appendix A]{lurie-ha}.

We record the following standard lemma: 
\lem\label{lemma: constructible-stalk}
Let $\mc{F}$ be a sheaf of (stable $\infty$-)categories on a simplicial complex $|\mc{K}|$ which is constructible with respect to the canonical stratification $\mc{S}_{\mc{K}}$. Given $\s \in \mc{K}$, we have $\mc{F}(\op{star}(\s)) = \mc{F}_c$, for any $c \in |\s|$. 
\elem
\pf
There is a homeomorphism $\op{star}(s)= [0,1) \tms \op{link}(\s) /(\{0\} \tms \op{link}(\s))$ which is compatible with the stratification. So $p$ admits a neighborhood base by open sets $U_t:= [0,t) \tms \op{link}(\s) /(\{0\} \tms \op{link}(\s))$. The restriction maps are isomorphisms since $\mc{F}$ is constructible: indeed, the stratification is independent of $t$ and the restriction maps are isomorphisms on each stratum. See also \cite[Sec.\ 4.2.1]{shende-takeda}. 
\epf

Suppose now that $(X, \l)$ is a Weinstein manifold of dimension $2n$, equipped with a polarization $\xi \subset TX$ (or more generally, with ``Maslov data" \cite[Sec.\ 10]{nadler-shende}). Given any conic subset $\L \sub X$, there is a sheaf of stable $\infty$-categories $\mu sh_{\L}(-)$ on $\L$, called the \emph{Nadler--Shende sheaf}, defined in \cites{shende-weinstein, nadler-shende}.\footnote{This sheaf of categories is often called the \emph{Kashiwara--Schapira stack} in the literature. We find this terminology misleading since $\mu sh_\L(-)$ is not a stack and was not, to our knowledge, defined by Kashiwara--Schapira (at least not in its current form).}  For the purposes of this paper, we will always assume that $\mu sh$ takes values in the dg derived category of $\Z$ or $\Z / 2$-graded $k$-modules.

We emphasize that although the \emph{construction} of $\mu sh_{\L}$ uses symplectic topology (i.e.\ it uses the fact that $\L$ is a conical subset of a Liouville manifold), $\mu sh_{\L}$ itself is not a symplectic object: it is just a sheaf of categories on the topological space $\L$.  

Here are some facts about the Nadler--Shende sheaf for Weinstein pairs with arboreal skeleton. Denote by $\mu sh_{\L}^c(-)$ the assignment $U \mapsto \mush_\L(U)^c$, where $\mush_\L(U)^c$ is the full subcategory of compact objects of $\mush_{\L}(U)$.
\prop\label{prop:microlocal-arboreal-properties}
Let $(X, A)$ be a Weinstein pair ($A$ could be empty) and suppose that $\op{skel}(X,A) = \fk{c}_{X, A} \sub X$ is an arboreal Lagrangian. 
\begin{enumerate}
\item\label{item:arboreal-constructible} The sheaf $\mu sh_{\fk{c}_{X, A}}$ is \emph{constructible} with respect to the singularity-order stratification.
\item\label{item:stalks} Suppose $p \in \fk{c}_{X, A} \subset X$ has a neighborhood $U$ containing $p$ and a symplectic embedding modeled on \eqref{eqn_arbo_model}. Let $T$ be the underlying tree of $\fk{T}$ and we endow $T$ with the natural quiver structure $\overrightarrow{T}$ with arrows pointing towards the root. Then the stalk of $\mu sh_{\fk{c}_{X, A}}^c (-)$ at $p$ satisfies
\eq 
\mu sh_{\fk{c}_X}^c(-)_p \simeq \op{Perf}(k[\overrightarrow{T}]).
\eeq
In particular, the stalks are independent of the signs on the edges of $\fk{T}$.
\item\label{item:compactly-generated} The category $\mu sh_{\fk{c}_{X, A}}(U)$ is compactly-generated \emph{for any} open set $U \sub \fk{c}_{X, A}$. Hence $\mu sh^c_{\fk{c}_{X, A}}(-)$ forms a cosheaf of categories on $\fk{c}_{X, A}$. In particular, if $\{U_\s\}_{\s \in \S}$ is an open cover of $\L$ closed under intersections and indexed by a finite poset $\S$, then 
\eq\label{eqn_local_global} \mu sh_{\fk{c}_{X, A}}^c(\fk{c}_{X, A}) \simeq \op{hocolim}_{\s \in \S} \mu sh_{\fk{c}_{X, A}}^c(U_\s).\eeq
\end{enumerate}
\eprop
\pf
\eqref{item:arboreal-constructible} follows from combining \Cref{definition:arboreal-lagrangian} with \cite[Lem.\ 5.3]{nadler-shende}. \eqref{item:stalks} was proved by \cite[Thm.\ 4.4]{nadler2017arboreal} in a somewhat different language since the definition of $\mu sh^{c}_{\L}$ was not in the literature (cf.\ also \cite[Thm.\ 4]{shende-weinstein} for our more modern formulation.)
For \eqref{item:compactly-generated}, combine \Cref{proposition:whitney-sing-order} and \cite[Rmk.\ 1.6, 1.7]{nadler-shende}. (More precisely, for any open set $U$, it follows from \Cref{proposition:whitney-sing-order} that $U \cap \fk{c}_{X, A}$ is Whitney stratified by isotropics. Hence it follows from \cite[Rmk.\ 1.6, 1.7]{nadler-shende} that $\mush_{\fk{c}_{X, A}}(U)$ is compactly-generated by microstalks, which are corepresentable.) Note that this property holds with much weaker assumptions on the skeleton. 
\epf

It turns out these properties of microlocal sheaves give rather strong control on the Rouquier dimensions of the sheaf of global sections, after using the topological results from Section \ref{subsec_arbo}.
\thm\label{thm_arbo_sheaf_upper}
Suppose that $(X^{2n}, A)$ is a Weinstein pair ($A$ could be empty), 
and $\op{skel}(X,A) = \fk{c}_{X, A} \sub X$ is an arboreal Lagrangian. Then
\eq\label{eqn_rouq_micro_small} \rdim \mu sh_{\fk{c}_{X, A}}^c(\fk{c}_{X, A}) \leq n \text{ if } n \leq 3 \text{ and} \eeq
\eq\label{eqn_rouq_micro} \rdim \mu sh_{\fk{c}_{X, A}}^c(\fk{c}_{X, A}) \leq 2n-3 \text{ for } n \geq 4. \eeq
\ethm
\pf
We prove the theorem assuming $A \neq \emptyset$ and the case $A = \emptyset$ follows with a similar (but simpler) proof. Note that $\fk{c}_{X,A}$ is the mapping cone $\fk{c}_{A} \times [0, \infty)$ glued to $\fk{c}_X$ along $\fk{c}_{A} \times \{0\}$ using the Liouville flow. Choose $N > 0$ and consider the ``truncation" $\Lambda_N := \fk{c}_X \cup \fk{c}_{A} \times [0, N]$. Note that $\Lambda_N$ is an arboreal Lagrangian with boundary. By Proposition \ref{proposition:whitney-sing-order}, the singularity-order stratification $\mc{S}_{\op{arb}}$ on $\Lambda_N$ defines a Whitney stratification. Note that $\mc{S}_{\op{arb}}$ restricts to a Whitney stratification $\mc{S}_{\op{arb}}^{\partial}$ along the ``boundary" $\fk{c}_A \times \{ N \} \subset \Lambda_N$.

Using Corollary \ref{cor_triang_exist}, we can find a topological triangulation $\mc{S}_{\mc{K}}$ of $(\Lambda_N, \mc{S}_{\op{arb}})$ which restricts to a topological triangulation $\mc{S}_{\mc{K}}^{\partial}$ of $(\fk{c}_A = \fk{c}_A \times \{ N \}, \mc{S}_{\op{arb}}^{\partial})$. In fact, we can extend the stratification $\mc{S}_{\mc{K}}$ induced by $\mc{K}$ to a stratification $\tilde{\mc{S}}_{\mc{K}}$ of the non-compact space $\fk{c}_{X, A}$ simply by introducing strata of the form $X_{\alpha} \times (N, \infty)$ where $X_\a$ is a stratum from the triangulation $\mc{S}_{\mc{K}}^{\partial}$.

Then the stars of the stratification $\tilde{\mc{S}}_{\mc{K}}$ on $\fk{c}_{X, A}$ define an open cover $\{ U_{\sigma} \}_{\sigma}$. This cover has the following nice properties.
\begin{enumerate}
\item This cover has depth at most $n+1$: given any chain of inclusions $U_{\sigma_1} \subsetneq \cdots \subsetneq U_{\sigma_k}$, we have $k \leq n+1$.
\item Given any $\op{star}(X_{\alpha})$ for a stratum $X_{\alpha}$ from $\tilde{\mc{S}}_{\mc{K}}$, the category $\mu sh_{\fk{c}_{X, A}}^c(X_{\alpha})$ is equivalent to $\op{Perf}(k[\overrightarrow{T}])$ for some quiver $\overrightarrow{T}$ with at most $n+1$ vertices. Indeed, this follows from \Cref{prop:microlocal-arboreal-properties}\eqref{item:stalks} and Lemma \ref{lemma: constructible-stalk}.
\end{enumerate}

The theorem is then a consequence of the local-to-global formula \eqref{eqn_local_global}, \Cref{lem_rdim_upper} and \Cref{lem_rdim_max}. Indeed, note that any quiver $\overrightarrow{T}$ with at most $4$ vertices must be of Dynkin ADE type so $\rdim \op{Perf}(k[\overrightarrow{T}]) = 0$ by Proposition \ref{proposition:elagin-gabriel}. This proves \eqref{eqn_rouq_micro_small}. For $n \geq 4$, note that for any chain of inclusions $U_{\sigma_1} \subsetneq \cdots \subsetneq U_{\sigma_k}$, there are least $4$ of the members satisfying $\mu sh_{\fk{c}_{X, A}}^c(U_{\s_i}) = 0$ and \eqref{eqn_rouq_micro} is proved again by using Proposition \ref{proposition:elagin-gabriel}.
\epf

\cor\label{cor_rouq_upper}
Suppose $(X^{2n}, A)$ is a Weinstein pair and $(X, \l)$ is equipped with a polarization. Then
\eq\label{eqn_rouq_fuk_small} \rdim \mc{W}(X,A) \leq n \text{ if } n \leq 3 \text{ and} \eeq
\eq\label{eqn_rouq_fuk} \rdim \mc{W}(X,A) \leq 2n-3 \text{ for } n \geq 4. \eeq
\ecor
\pf
By \cite[Cor.\ 7.28]{gps3}, we can assume that the Weinstein manifold $(X, \l)$ and Weinstein hypersurface $A$ are real analytic and the relative core $\fk{c}_{X, A}$ is subanalytic singular isotropic. Then \cite[Thm.\ 1.4]{gps3} asserts the equivalence
\eq \op{Perf}\mc{W}(X,A) \simeq \mu sh_{\fk{c}_{X, A}}^c(\fk{c}_{X, A}). \eeq
Because $X$ is endowed with a polarization, we can find a Weinstein homotopy $\{ \lambda_t \}_{0 \leq t \leq 1}$ and a Weinstein hypersurface $B$ with $\l_0 = \l$ such that for the Weinstein structure $\lambda_1$ we have $\op{skel}(B) \setminus \op{skel}(A)$ being a smooth Legendrian and $\fk{c}(X,B) = \op{skel}((X, \l_1), B)$ is arboreal by Theorem \ref{theorem:arborealization}. Then using Lemma \ref{corollary:quot-rdim} and the stop removal result \cite[Thm.\ 1.5 (2)]{nadler-shende}, we see that for the Weinstein structure $\l_1$,
\eq \rdim \mu sh^c_{\fk{c}_{X,A}}(\fk{c}_{X,A}) \leq \rdim \mu sh^c_{\fk{c}_{X,B}} (\fk{c}_{X,B}) \eeq
because $\mu sh^c_{\fk{c}_{X,A}}(\fk{c}_{X,A})$ is a quotient of $\mu sh^c_{\fk{c}_{X,B}} (\fk{c}_{X,B})$. The corollary is thus proved by Theorem \ref{thm_arbo_sheaf_upper} the invariance of $\mu sh^c_{\fk{c}_{X,A}}(\fk{c}_{X,A})$ under Weinstein homotopy \cite[Thm.\ 1.9(2), Rmk.\ 8.16]{nadler-shende}.
\epf

\rmk
Our proof of \Cref{cor_rouq_upper} crucially relies on the arborealization result of \'Alvarez-Gavela--Eliashberg--Nadler stated in  \Cref{theorem:arborealization}. However, as remarked by the anonymous referee, we could instead have appealed to an earlier (and easier) result of Nadler \cite{nadler-non-char}. Nadler's result implies that there exists some arboreal space $\Lambda'$ such that $\mu sh_{\fk{c}_{X,A}}(\fk{c}_{X,A})= \mu sh_{\Lambda'}(\Lambda')$. However, in contrast to \cite{AGEN3}, Nadler's result \emph{does not} assert that $\Lambda'$ is the skeleton $(X,A)$ for some deformed Liouville structure. Note that \cite{nadler-non-char} does not require any polarizability hypotheses (which, however, are still needed when appealing to \cite{gps3}, which is necessary for turning the sheaf-theoretic results into statements for Weinstein manifolds). 
\ermk

\subsection{Applications to Orlov's conjecture}\label{subsection:orlov-conj}
Now we can apply the results obtained above to resolve more cases of Orlov's conjecture \ref{conj_orlov_intro}.
\pf[Proof of Theorem \ref{thm:orlov_gen_state_intro}]
Given a quasi-projective variety $Y$ over $\mathbb{C}$, the assumption on homological mirror symmetry identifies $D^{b} \op{Coh}(Y)$ with the derived (partially) wrapped Fukaya category of some Weinstein pair satisfying the assumptions of Corollary \ref{cor_rouq_upper}. Therefore we have $\rdim D^{b} \op{Coh}(Y) \leq \dim_{\mathbb{C}} Y$. On the other hand, \cite[Prop 7.17]{rouquier} tells us that $\rdim D^{b} \op{Coh}(Y) \geq \dim_{\mathbb{C}} Y$. Therefore the Orlov conjecture holds.
\epf

\rmk
When $n \geq 4$ under the same assumptions, we obtain the bound
\eq \rdim D^{b} \op{Coh}(Y) \leq 2n-3. \eeq
This improves the general upper bound $\rdim D^{b} \op{Coh}(Y) \leq 2n$ as in \cite[Prop.\ 7.9]{rouquier}.
\ermk

Next we list the new instances of Orlov's conjecture covered by Theorem \ref{thm:orlov_gen_state_intro}.

\ex[Log Calabi--Yau surfaces]\label{ex_logCY}
In dimension $2$, we can consider non-toric examples arising in the work of Keating \cite{keating}, as we now explain. 

Let $Y_{p,q,r}$ be the complex algebraic variety constructed as follows. We consider the toric anti-canonical divisor of $D \sub \bb{P}^2$, which is a wheel of three $\bb{P}^1$. Pick $p$ copies of a point lying on the first component, $q$ copies of a point lying on the second component and $r$ copies lying on the third. Assume these points are colinear, lying in the interior of each irreducible component, and that $1/p+1/q+1/r<1$. Now blow up $\bb{CP}^2$ at these points. (By definition, the blowup along a point $x \in D$ with multiplicity $m>1$ means that we first blow up at $x$, then we blow up at the the intersection of the exceptional divisor with the strict transform of $D$, and so forth $m$ times).  

According to the main result of \cite{keating}, there are equivalences of categories
\eq\label{equation:keating-examples}
D^b \op{Coh}(Y_{p,q,r}) \cong  \op{Perf} \mc{W}(\mc{T}_{p,q,r}, V), \; \; D^b \op{Coh}(Y_{p,q,r}-D) \cong  \op{Perf} \mc{W}(\mc{T}_{p,q,r})
\eeq
where 
\begin{itemize}
\item $\mc{T}_{p,q,r}$ is a Weinstein $4$-manifold (a certain Milnor fiber in $\C^3$);
\item $V \sub \d_\infty \mc{T}_{p,q,r}$ is a Weinstein hypersurface;
\item (by abuse of notation) $D$ denotes the strict transform of the toric-anticanonical divisor on $\bb{P}^2$.
\end{itemize}

In particular, $\mc{T}_{p,q,r}$ is a complete intersection, so it admits a polarization. Therefore, we see Conjecture \ref{conj_orlov_intro} holds for $Y_{p,q,r}$ and $Y_{p,q,r}-D$. 
\eex

\ex[Three-dimensional toric varieties]\label{ex_toric}
Homological mirror symmetry for toric varieties (on the B-side) is proved by Fang--Liu--Treumann--Zaslow in \cite{FLTZ}. We follow the more modern fomulation from \cite[Cor.\ 6.16]{gps3}. For any toric variety $\mathbf{T}$ of dimension $n$ defined over $\mathbb{C}$, there exists a singular Legendrian $\partial_{\infty} \Lambda_{\mathbf{T}} \subset \partial_{\infty} T^* T^n$ such that
\eq  \op{Coh}(\mathbf{T}) \cong \mc{W}(T^* T^n, \partial_{\infty} \Lambda_{\mathbf{T}}). \eeq
According to the work of Zhou \cite{zhou2020}, which generalizes the arguments in \cite{gammage2017mirror}, the singular Legendrian $\partial_{\infty} \Lambda_{\mathbf{T}}$ arises as the skeleton of a Weinstein hypersurface $F_{\mathbf{T}}$. Note that the ambient Weinstein manifold $(\mathbb{C}^*)^n = T^* T^n$ comes with a natural polarization therefore Theorem \ref{thm:orlov_gen_state_intro} applies. We conclude that Conjecture \ref{conj_orlov_intro} is true for all $3$-dimensional toric varieties, using the fact that $\mc{W}(T^* T^n, \partial_{\infty} \Lambda_{\mathbf{T}})$ is equivalent to $\mc{W}(T^* T^n, F_{\mathbf{T}})$.
\eex




Using similar techniques, one can also obtain upper bounds for the Rouquier dimension of derived categories of coherent sheaves on singular algebraic varieties. Although Orlov's conjecture \ref{conj_orlov_intro} is formulated for smooth quasi-projective varieties, and is known to be false in general in the singular case (see \Cref{example:hara}), its extension to possibly singular varieties/stacks holds for the following examples. 

\ex[Toric boundary divisors]
We recall the following result of homological mirror symmetry by Gammage--Shende \cite[Thm.\ 1.0.1]{gammage2017mirror}. Denote by $\partial \mathbf{T}_{\Sigma}$ the toric boundary divisor of a smooth toric stack $\mathbf{T}$ over $\mathbb{C}$. Then there exists a Laurent polynomial $W_{\mathbf{T}}: (\mathbb{C}^*)^n \rightarrow \mathbb{C}$ and a natural Weinstein structure on the generic fiber $F_{W_{\mathbf{T}}}$ of $W_{\mathbf{T}}$ such that
$$\op{Coh}(\partial \mathbf{T}_{\Sigma}) \cong \mc{W}(F_{W_{\mathbf{T}}}).$$
(The special case where $\mathbf{T}= \mathbb{A}^n$ recovers mirror symmetry for the $(n-1)$-dimensional pair of pants, which was proved independently by Lekili--Polishchuk \cite[Cor.\ 1.1.2]{lekili-polishchuk}; see also \cite{nadler-wrapped-microlocal}.)
Note that $F_{W_{\mathbf{T}}}$ is a hypersurface in the affine variety $(\mathbb{C}^*)^n$ so it admits a polarization. When $\dim_{\mathbb{C}}(\partial \mathbf{T}_{\Sigma}) \leq 3$, Theorem \ref{thm:orlov_gen_state_intro} tells us 
\eq
\rdim D^{b}\op{Coh}(\partial \mathbf{T}_{\Sigma}) = \dim_{\mathbb{C}}(\partial \mathbf{T}_{\Sigma}).
\eeq
\eex


\ex[Mirrors of Milnor fibers of weighted homogeneous polynomials]
A weighted homogeneous polynomial $\mathbf{w} \in \mathbb{C}[x_1, \dots, x_n]$ with an isolated singularity at the origin is called \emph{invertible} if there exists a matrix $(a_{ij})_{i,j=1}^n$ with integer-valued entries and non-zero determinant such that 
$$\mathbf{w} = \sum_{i=1}^n \prod_{j=1}^n x_j^{a_{ij}}.$$
The corresponding weight system $(d_1, \dots, d_n, h)$ is determined uniquely by requiring 
$$w(t^{d_1}x_1, \dots, t^{d_n}x_n) = t^h w(x_1, \dots, x_n)$$ 
for all $t \in \mathbb{G}_m$ and $\op{gcd}(d_1, \dots, d_n, h) = 1$. The \emph{transpose} of $\mathbf{w}$ is defined to be
$$\check{\mathbf{w}} = \sum_{i=1}^n \prod_{j=1}^n x_j^{a_{ji}}$$
and let $(\check{d}_1, \dots, \check{d}_n, \check{h})$ be the weight system of $\check{\mathbf{w}}$.
Denote by $\check{V}_{\check{\mathbf{w}}}$ (the Liouville completion of) the Milnor fiber of $\check{\mathbf{w}}$, which is a Weinstein manifold coming with a global Lagrangian plane field. The symmetry group of $\mathbf{w}$ is defined by
$$\Gamma_{\mathbf{w}} := \{ (t_1, \dots, t_n) \in (\mathbb{G}_m)^n | t_1^{a_{11}} \cdots t_n^{a_{1n}} = \cdots = t_1^{a_{n1}} \cdots t_n^{a_{nn}} \}.$$
It is conjectured in \cite[Conj.\ 1.4]{lekili-ueda} and proved in \cite[Thm.\ 1.2]{lekili2020homological} and \cite[Thm.\ 1.6]{gammage-bh} that
$$\op{mf}(\mathbb{A}^{n+1}, \mathbf{w} + x_0 \cdots x_n, \Gamma_{\mathbf{w}}) \cong \mc{W}(\check{V}_{\check{\mathbf{w}}}),$$
where $\op{mf}(\mathbb{A}^{n+1}, \mathbf{w} + x_0 \cdots x_n, \Gamma_{\mathbf{w}})$ is the (dg) category of equivariant matrix factorizations. If $\check{h} - \sum_{i=1}^n \check{d}_i$ is positive, a theorem of Orlov \cite[Thm.\ 3.11]{orlov} identifies $\op{mf}(\mathbb{A}^{n+1}, \mathbf{w} + x_0 \cdots x_n, \Gamma_{\mathbf{w}})$ with $\op{Coh}(\mathcal{Z}_{\mathbf{w}})$. Here $\mathcal{Z}_{\mathbf{w}}$ is the stack
$$[(\op{Spec} \mathbb{C}[x_0,\dots, x_n] / \mathbf{w} + x_0 x_1 \cdots x_n) / \Gamma_{\mathbf{w}}].$$
For $n \leq 4$, Theorem \ref{thm:orlov_gen_state_intro} applies so we conclude that
\eq
\rdim D^{b} \op{Coh}(\mathcal{Z}_{\mathbf{w}}) = \dim_{\mathbb{C}} \mathcal{Z}_{\mathbf{w}}.
\eeq
\eex

\subsection{Applications to non-commutative motives}\label{subsection:motives}
A functor from the category of small $A_\infty$ categories over a field $k$ to some additive category $\mc{D}$ is called an \emph{additive invariant} if: 
\begin{itemize} 
\item[(i)] it sends Morita equivalences to isomorphisms;
\item[(ii)] it sends semiorthogonal decompositions to direct sums.
\end{itemize}
Well-known examples of additive invariants include Hochschild homology, algebraic $K$-theory, cyclic homology, topological Hochschild homology, etc. We recommend \cite[Sec.\ 2.2]{tabuada-motives} for details and many more examples. As explained in \cite[Sec.\ 2.3]{tabuada-motives}, all additive invariants factor through a \emph{universal} additive invariant, which takes values in the additive category $\op{Hmo}_0(k)$ of \emph{non-commutative Chow motives}; the objects in this category are small dg/$A_\infty$ categories, morphisms are obtained by taking the Grothendieck group of a certain category of bimodules. 

As a byproduct of the argument of the previous section, we obtain the following (to us) rather surprising result: 
\cor\label{corollary:additive-invariant}
Let $(X, A)$ be a polarizable Weinstein pair such that $\mc{W}(X, A)$ is proper (it is always smooth).  Let $U(-)$ be the universal additive invariant. Then $U(\mc{W}(X,A))$ is a summand of $\oplus_{\s \in \S} U(\op{Perf} k[T_\s])$, where $\S$ is a finite set and $T_\s$ is a tree quiver. 
\ecor
\pf 
The proof of \Cref{thm_arbo_sheaf_upper} and \Cref{cor_rouq_upper} furnishes a localization map $p: \mc{G}:= \op{Groth}_{\s \in \S} \op{Perf}(k [T_{\s}]) \to \mc{W}(X, A)$, where $\S$ is a finite poset. 
Since $p$ is a localization, the pullback $p^*$ on modules is fully faithful. Moreover, since $\op{Perf}\mc{W}(X,A)$ and $\mc{G}$ are smooth and proper, $p^*$ takes $\op{Perf}\mc{W}(X,A)$ into $\op{Perf}\mc{G} \subset \op{Mod}\mc{G}$. This implies that there is a semi-orthogonal decomposition $\op{Perf} \mc{G} = \op{Perf} \langle \op{Perf}\mc{W}(X,A), \mc{B} \rangle$, where $\mc{B}$ is right orthogonal to $p^*$ and $\mc{B}$ is smooth and proper; cf.\ \cite[Sec.\ 4]{tabuada-k}. Hence $U(\op{Perf} \mc{W}(X,A))) \oplus U(\mc{B})= U(\mc{G})= \oplus_{\s \in \S} U(\op{Perf} k[T_\s])$. 
\epf

Given a tree quiver $T$, it is well-known that $\op{Perf} k[T]$ admits a full exceptional collection indexed by the vertices of $T$ (see e.g.\ \cite[Sec.\ 2.4]{elagin2020calculating}). Hence $U(\op{Perf} k[T])= U(k)\oplus \dots \oplus U(k)$ ($|T|$ times).

\cor \label{corollary:additive-invariant2}
Under the assumptions of \Cref{corollary:additive-invariant}, $U(\mc{W}(X,A))$ is a direct summand of $U(k)^m$, for some $m \geq 1$. 
\qed
\ecor

Of course, by virtue of $U(-)$ being the \emph{universal} additive invariant, the statement of \Cref{corollary:additive-invariant2} remains true if one substitutes for $U(-)$ one's favorite additive invariant, such as $K_\bu(-), HH_\bu(-), THH(-)$, etc. 

\ex
It is a standard fact that $HH_\bu(\op{Perf} k[T])$ is concentrated in degree zero (see e.g.\ \cite[Sec.\ 5.3]{shende-takeda}). Hence, under the assumptions of \Cref{corollary:additive-invariant}, $HH_\bu(\mc{W}(X, A))$ is concentrated in degree $0$. As a sanity check, suppose that $k=\C$ and let $(T^* T^n, \partial_{\infty} \Lambda_{\mathbf{T}})$ be mirror to a smooth toric variety $X$ as in \Cref{ex_toric}. By the HKR theorem, $HH_\bu(D^b\op{Coh}(X))= \oplus_p H^{p, p-*}(X)$.  Hence we learn from the fact that $HH_\bu(\mc{W}(T^* T^n, \partial_{\infty} \Lambda_{\mathbf{T}}))$ is concentrated in degree zero that $h^{p,q}=0$ whenever $p \neq q$. This is indeed true: it is a consequence of the general fact that the cohomology of a smooth toric variety coincides with its Chow ring.
\eex

Let us mention another consequence of \Cref{corollary:additive-invariant}. There are certain non-commutative analogs of classical conjectures in arithmetic geometry which were introduced by Tabuada (building on proposals of Kontsevich); see \cite{tabuada-weil}.  These are conjectures about smooth and proper dg-categories. However, as shown by Tabuada, they recover classical conjectures in arithmetic geometry when specialized to $\op{Perf}(X)$ for an appropriate variety $X$.  By combing \Cref{corollary:additive-invariant} and \cite[Cor.\ 12.2]{tabuada-weil}, we learn:
\cor
Assume $k$ is a finite field of positive characteristic. If $(X,A)$ is a polarizable Weinstein pair such that $\mc{W}(X, A)$ is proper, then $\mc{W}(X, A)$ satisfies the non-commutative Weil conjecture and the strong non-commutative Tate conjecture. 
\qed
\ecor

\section{Central actions on the wrapped Fukaya category}\label{section:lower-bounds-action}

The main goal of this section is to prove \Cref{corollary:main-central-bound}, which gives a lower bound on the Rouquier dimension of certain wrapped Fukaya categories. We do this by considering a central action of $HH^\bu(\mc{W}(X))$ on $\mc{W}(X)$ and adapting methods of Bergh--Iyengar--Krause--Opperman \cite{b-i-k-o} to deduce lower bounds on the Rouquier dimension. 

\subsection{Central actions on triangulated categories}

Let $\mc{T}= (\mc{T}, \S, \Delta)$ be a triangulated category. Let $Z^d(\mc{T})$ be the set of natural transformations $\eta: \op{id} \to \S^d$ such that $\eta \circ \S = (-1)^d \S \circ \eta$. Observe that $Z^\bu(\mc{T})= \oplus_{d \in \Z} Z^d(\mc{T})$ inherits the structure of a graded-commutative ring by composition. We call $Z^\bu(\mc{T})$ the \emph{graded center} of $\mc{T}$.

We will be mainly interested in the case where $\mc{T}= H^0(\mc{C})$ for $\mc{C}$ being a pre-triangulated $A_\infty$ category. In that case, the graded center of $H^0(\mc{C})$ is precisely the set of graded natural transformations of the identity in $H^\bu(\mc{C})$; cf.\ \cite[Sec.\ 4.2]{lowen}. Concretely, the elements of the graded center of $H^0(\mc{C})$ are products 
\eq \zeta= (\zeta_K)_{K \in \mc{C}} \in \prod_{K \in \mc{C}} H^\bu(K,K) \eeq 
subject to the naturality condition that the following diagram commutes
\eq\label{equation:central-action-diagram}
\begin{tikzcd}
H^\bu(K, L) \ar{d}{1 \otimes \zeta_K} \ar{r}{\zeta_L \otimes 1} & H^\bu(L, L) \otimes H^\bu(K, L) \ar{d}{} \\
H^\bu(K, L) \otimes H^\bu(K, K) \ar{r} & H^\bu(K, L).
\end{tikzcd}
\eeq

Given a graded-commutative ring $R$, a \emph{central action} of $R$ on the triangulated category $\mc{T}$ is a morphism of graded-commutative rings $R \to Z^\bu(\mc{T})$. 

Any triangulated category $\mc{T}$ which comes from a pre-triangulated $A_\infty$ category $\mc{C}$ (i.e.\ which is of the form $\mc{T}= H^0(\mc{C})$) admits a canonical central action
\eq\label{equation:char-mor} \chi_{\mc{C}}: HH^\bu(\mc{C}, \mc{C}) \to Z^\bu(\mc{T}), \eeq
where $HH^\bu(\mc{C}, \mc{C})$ is the Hoschschild cohomology of $\mc{C}$. In general, this map is neither injective nor surjective. (Recall that the Hochschild cohomology of an $A_\infty$ category is naturally a graded-commutative ring. We refer to \cite[Sec.\ 2.13]{ritter-smith} for a very detailed overview of Hochschild cohomology of $A_\infty$ categories.)

The morphism \eqref{equation:char-mor} which realizes this canonical action is called the \emph{characteristic morphism}. To describe the characteristic morphism concretely, let us recall the bar model for the Hochschild cochain complex.

\eq
CC^k(\mc{C}, \mc{C}) = \prod_{\begin{smallmatrix}l\geq 0\\K_0,\ldots,K_l \in \mc{C} \end{smallmatrix}} \op{Hom}_{k-l}(\hom_\mc{C}(X_{l-1}, X_l) \otimes \dots \otimes \hom_{\mc{C}}(X_0, X_1), \hom_\mc{C}(X_0, X_l))
\eeq
where $\op{Hom}_s(-,-)$ is the space of degree $s$ morphisms of graded vector spaces.

The differential applied to an element $\phi$ of degree $d$ is
\begin{align}
\delta(\varphi)(a_l,  \dots,  a_1) &= \sum_{r,s} (-1)^{\clubsuit_1}   \mu^{l-r+s}_\mc{C}(a_l,\dots, a_{r+1},\varphi^{r-s+1}(a_r,\dots, a_s), a_{s-1},\dots,a_1) \\ 
&- (-1)^d \sum_{r,s} (-1)^{\clubsuit_2} \varphi^{l-r+s}(a_l,\dots,a_{r+1},\mu^{r-s+1}_\mc{C}(a_r,\dots,a_s),a_{s-1},\dots, a_1),
\end{align}
where $\clubsuit_1:= d (|a_{s-1}| + \dots + |a_1| - (s-1))$ and $\clubsuit_2=(|a_{s-1}| + \dots + |a_1| - (s-1))$.

Given any object $K \in \mc{C}$, there is a natural projection map
\begin{align}
\pi_K: CC^k(\mc{C}, \mc{C}) &\to \op{Hom}_k(k, \hom_\mc{C}(K, K)) \equiv \hom_\mc{C}(K, K) \\
\varphi \mapsto \varphi^0_K \equiv \varphi^0_K(1).
\end{align}

It can be verified (cf.\ \cite[Sec.\ 2.14]{ritter-smith}) that the projection map is a morphism of rings. The characteristic morphism is then precisely 
\eq  \chi_\mc{C}(\varphi)= (\varphi^0_K(1))_{K \in \mc{C}}. \eeq

The fact that this is a central action follows from the definition of the Hochschild differential. Indeed, given a closed element $x \in \hom(K, L)$, we have $0=\delta(\varphi)(x)=  \mu^2_\mc{C}(\varphi_L^0(1), x) -  \mu^2_\mc{C}(x, \varphi_K^0(1))$. This matches \eqref{equation:central-action-diagram}. 

\rmk
The Hochschild cohomology of an $A_\infty$ category is canonically equivalent to the space of natural transformations of the identity functor, i.e.\ the ``derived center" of the category. The characteristic morphism thus realizes a map from the derived center of an $A_\infty$ category to the center of its derived category. For more on this perspective, see \cite{neumann-szymik} and the references therein.
\ermk

\subsection{Central actions and the Rouquier dimension}\label{subsection:commutative-alg-part}

The goal of this section is to prove the following theorem. We refer to \Cref{subsubsection:biko-adapt} for motivation and a discussion of the proof. 

\thm\label{theorem:main-biko-adapt}
Let $\mc{T}$ be a $\Z/m$-graded triangulated category $(1 \leq m \leq \infty)$ and let $R$ be a finite type $k$-algebra acting centrally on $\mc{T}$. Suppose that there exists a (split-)generator $G \in \mc{T}$ such that $\op{Hom}_{\mc{T}}^\bu(G,G)$ is a Noetherian $R$-module. Then $\rdim \mc{T} \geq \op{dim}_R \op{Hom}_{\mc{T}}^\bu(G,G).$
\ethm

\subsubsection{Recollections from commutative algebra}

Given a module $M$ over a commutative ring $R$, the \emph{annihilator} is the ideal $\op{Ann}_R(M):= \{ r \in R \mid r \cdot m =0 \text{ for all } m \in M\} \sub R$.  

The \emph{Krull dimension} of a commutative ring $R$ is defined as the length of the longest strictly increasing chain of prime ideals $\fk{p}_0 \subsetneq \dots \subsetneq \fk{p}_l$ and is denoted by $\op{dim}(R)$. Given a finitely generated module $M$ over $R$, one defines the Krull dimension of $M$ to be $\op{dim}_R(M):= \op{dim}(R / \op{Ann}_R(M))$. Finally, the Krull dimension of a (Zariski) closed subset of $\op{Spec} R$ is the length of the longest strictly increasing chain of closed subsets. When the ring $R$ is clear from the context, we usually write the Krull dimension as $\op{dim}(M) : = \op{dim}_R(M)$.

The \emph{support} of a module $M$ over $R$ is the set $\op{Supp}_R(M):= \{ \fk{p} \in \op{Spec} R \mid M_{\fk{p}} \neq 0\}$. The support enjoys the following standard properties:

\begin{fact}[Properties of the support]\label{fact:support} Let $R$ be a commutative ring.
\begin{enumerate}
\item\label{item:supp-closed} If $M$ is a finitely generated $R$ module, then $\op{Supp}_R(M)= V(\op{Ann}_R (M)):= \{\fk{p} \in \op{Spec} R \mid \op{Ann}_R(M) \sub \fk{p}\}$. In particular, $\op{dim}_R(M)= \op{dim}(\op{Supp}_R(M))$.
\item If $0 \to M_1 \to M \to M_2 \to 0$ is an exact sequence of $R$-modules, then $\op{Supp}_R(M) = \op{Supp}_R(M_1) \cup \op{Supp}_R(M_2)$. 
\item\label{item:supp-tensor}  If $M, N$ are finitely generated $R$-modules, then $\op{Supp}_R( M \otimes_R N)= \op{Supp}_R(M) \cap \op{Supp}_R(N)$. 
\item\label{item:supp-quotient} If $M$ is a finitely generated $R$-module and $I \sub R$ is an ideal, then $\op{Supp}_{R/I}(M/I)= V(I) \cap \op{Supp}_R(M)$. 
\end{enumerate}
\end{fact}

We also need to review the notion of regular sequences.
\defi 
Let $R$ be a commutative ring and let $M$ be an $R$-module. A sequence $x_1,\dots,x_n$ is \emph{$M$-regular} (or just \emph{regular}) if for each $k \leq n$, we have that $x_k$ is not a zero-divisor in $M/(x_1,\dots,x_{k-1})M$, and that $(x_1,\dots,x_k)M \neq M$.  Given an ideal $I \sub R$, we let $\op{depth}_I M$ denote the length of the longest $M$-regular sequence contained in $I$. (If $R$ is local and $I \sub R$ is the maximal ideal, we usually just write $\op{depth}_I M$ as $\op{depth} M$). 
\edefi

\begin{fact}[Regular sequences]\label{fact:regular-sequence} \label{lemma:depth-equal} Let $R$ be a commutative ring and let $M$ be a module.
\begin{enumerate}
\item\label{item:depth-local} Suppose $R$ is Noetherian and $M$ is finitely generated. If $\fk{m} \in \op{Spec} R$ is maximal, then 
\eq 
\op{depth}_{\fk{m}} M = \op{depth}_{\fk{m} R_\fk{m}} M_{\fk{m}}.
\eeq
\item\label{item:reg-sequence-power} Let $x_1,\dots,x_n$ be elements of $R$ and let and $a_1,\dots, a_n>0$ be integers. Then $x_1,…,x_n$ is an $M$-regular sequence if and only if $x_1^{a_1},\dots, x_n^{a_n}$ is an $M$-regular sequence. 
\item\label{item:reg-sequence-localization} Suppose that $x_1,\dots,x_n$ is an $M$-regular sequence. Suppose that there exists $\fk{p} \in \op{Supp}_R(M)$ such that $x_i \in \fk{p}$. Then $x_1,\dots,x_n$, viewed as a sequence of elements of $R_{\fk{p}}$, is an $M_{\fk{p}}$ regular sequence. 
\end{enumerate}
\end{fact}
\pf
\eqref{item:depth-local} follows from \cite[Cor.\ 1.27]{mathew-notes}. \eqref{item:reg-sequence-power} is \cite[\href{https://stacks.math.columbia.edu/tag/0AUH}{Tag 0AUH}]{stacks-project}. For \eqref{item:reg-sequence-localization}, note that since multiplication by $x_i$ is injective on $M/(x_1,\dots, x_{i-1})M$, the same is true after localizing (recall that localization is exact). Finally, it follows from Nakayama's lemma that $M_\fk{p}/(x_1,\dots, x_n) M_{\fk{p}} = 0$ iff $M_{\fk{p}}=0$.
\epf



\lem\label{lemma:dimension-localization}
Let $R$ be a finite type algebra over a field $k$. Suppose that $R$ is equidimensional (meaning that all its minimal primes have the same dimension). Given a minimal prime $\fk{p} \subset R$, if $f \in R \setminus \fk{p}$, then $\op{dim}(R)= \op{dim}(R_f)$.  
\elem
\pf  
Observe that $\op{dim}(R_f)= \op{dim}(R_f/ \fk{p})$. We have a canonical isomorphism $R_f/ \fk{p} = (R/\fk{p})_f$. Finally, since $R/\fk{p}$ and $(R/\fk{p})_f$ are integral domains and finite type algebras over $k$, it follows from \cite[Thm.\ 11.2.1]{vakil} that $\op{dim}(R)= \op{dim}(R/\fk{p})= \op{tr. deg} K(R/ \fk{p}) = \op{tr. deg} K((R/ \fk{p})_f) = \op{dim}((R/\fk{p})_f)= \op{dim}(R_f)$. 
\epf

Note that the hypothesis that $R$ is of finite type is essential. (For example, $R= k[x]_{(x)}$ has Krull dimension $1$. However, $x$ is not contained in the unique minimal prime $(0)$ and $R_x= k(x)$ has dimension $0$.)

\cor\label{corollary:dim-localization}
Let $R$ be a finite type algebra over a field $k$ and let $M$ be a finitely generated $R$-module. Let $\fk{p} \in \op{Ass}_R(M)$ be an associated prime with the property that $\op{dim}(M/ \fk{p}) = \op{dim}(M)$. If $f \in R \setminus \fk{p}$ and $f$ is contained in every other associated prime of $M$, then:
\begin{itemize}
\item $\op{Ass}_{R_f}(M_f)$ has exactly one associated prime (hence $\op{Supp}_{R_f}(M_f)$ has exactly one minimal prime). 
\item $\op{dim}(M)= \op{dim}(M_f)$. 
\end{itemize}
\ecor

\pf  
It follows from \cite[\href{https://stacks.math.columbia.edu/tag/05BZ}{Tag 05BZ}]{stacks-project} that $\fk{p}$ is the unique associated prime of  $M_f$. Hence \cite[\href{https://stacks.math.columbia.edu/tag/02CE}{Tag 02CE}]{stacks-project} implies that it is also the unique minimal prime. This is the first claim.

We have that $\op{dim}(M_f)= \op{dim}(M_f/ \fk{p})= \op{dim}((M/ \fk{p})_f)$. It is enough to show that $\op{dim}(M/\fk{p})= \op{dim}((M/\fk{p})_f)$. It follows from \Cref{fact:support} that $\op{Supp}(M/\fk{p})= V(\fk{p}) = \op{Spec}(R/ \fk{p})$, while the support of $(M/\fk{p})_f$ is $V(\fk{p}) \cap D(f)= \op{Spec}((R/\fk{p})_f)$ (where $D(f):= \{ \fk{q} \in \op{Spec} R \mid f \notin \fk{q} \}$). Hence the second claim follows from \Cref{lemma:dimension-localization} applied to $R/\fk{p}$. 
\epf

\prop\label{proposition:exists-reg-sequence}
Let $R$ be a finite type algebra over a field $k$ and let $M$ be a finite $R$-module of Krull dimension $c \geq 0$. Then there exists an $M$-regular sequence of length $c$. 
\eprop
\pf
Choose a minimal associated prime $\fk{p} \in \op{Ass}(M)$ such that $M/\fk{p}$ has the same dimension as $M$ (see \cite[\href{https://stacks.math.columbia.edu/tag/02CE}{Tag 02CE}]{stacks-project}). Since $M$ is finite, $\op{Ass}_R(M)$ is finite \cite[\href{https://stacks.math.columbia.edu/tag/00LC}{Tag 00LC}]{stacks-project}. Choose an element $f_1 \in R$ such that $f_1 \notin \fk{p}$ but $f_1 \in \fk{q}$ for all $\fk{q} \in \op{Ass}(M) \setminus \fk{p}$. According to \Cref{corollary:dim-localization}, $\op{dim}(M_{f_1})= \op{dim}(M)=c$. 

Now choose a tower of primes $\fk{p}R_{f_1}=\fk{p}_0 \sub... \sub \fk{p}_c$ in $R_{f_1}$ of maximal length. Pick $r_1 \in \fk{p}_1 \setminus \fk{p}_0$. Then $r_1$ is not a zero-divisor of $M_{f_1}$ (indeed, by \cite[\href{https://stacks.math.columbia.edu/tag/00LD}{Tag 00LD}]{stacks-project} $\fk{p}_0$ contains all zero divisors). According to \Cref{fact:support}\eqref{item:supp-quotient}, the dimension of the quotient module is at least $c-1$. 

We now replace $M$ with $M_{f_1}/r_1$ and $R$ with $R_{f_1}/r_1$, and repeat the same argument. After $c$ steps, we end up with $(\dots(M_{f_1}/ r_1)_{f_2} \dots/r_{c-1})_{f_c}/r_c$. This is the same as $M_{f_1\dots f_c}/(r_1,\dots, r_c)$. (Note that this module has dimension zero, so it is not the zero module!) So $r_1,\dots, r_c$ is a regular sequence for $M_{f_1\dots f_c}$. 

Since $M_{f_1\dots f_c}$ is not zero, it has a maximal ideal $\fk{m}$ in its support. We now localize at $\fk{m}$. Observe first that $(M_{f_1\dots f_c})_\fk{m}= M_{\fk{m}}$ (indeed, it follows from \Cref{fact:support}\eqref{item:supp-tensor} that $\fk{m} \in D(f_1\dots f_m)$, so $f_1\dots f_m$ is a unit in $R_\fk{m}$). Moreover, it follows from \Cref{fact:regular-sequence}\eqref{item:reg-sequence-localization} that $r_1,\dots,r_c$ is still a regular sequence for the $R_{\fk{m}}$ module $M_{\fk{m}}$. The conclusion now follows from \Cref{fact:regular-sequence}\eqref{item:depth-local}. 
\epf

\rmk
In fact, more is true: with the notation of the proof of \Cref{proposition:exists-reg-sequence}, one can apply \cite[Sec.\ 3.1, Ex.\ 23]{kaplansky} to show that there exists an $M$-regular sequence generating $\fk{m}$ which is \emph{stable under permutations}. We thank Srikanth Iyengar for this reference.
\ermk

\subsubsection{Support varieties}

We begin with some standard facts about Noetherian modules. First, recall that a module over a Noetherian ring is Noetherian iff it is finitely generated. Also, Noetherian modules are well behaved under exact sequences, meaning that if $0 \to M_1 \to M \to M_2 \to 0$ is an exact sequence of $R$-modules, then $M$ is Noetherian iff $M_1, M_2$ are Noetherian. 

\lem\label{lemma:noetherian-stable}
Assume that $R$ is Noetherian and acts centrally on $\mc{T}$. Given any object $X \in \mc{T}$, the subcategory $\{Y \in \mc{T} \mid \op{Hom}_\mc{T}^\bu(X, Y)\; \text{is Noetherian over } R \} \sub \mc{T}$ is a thick triangulated subcategory. Similarly for $\{Y \in \mc{T} \mid \op{Hom}_\mc{T}^\bu(Y, X)\; \text{is Noetherian over } R \} \sub \mc{T}$.
\elem
\pf
We only treat the first case as the second one is similar. Consider an exact triangle $Y_1 \to Q \to Y_2 \to $ where $\op{Hom}_\mc{T}^\bu(X, Y_i)$ is Noetherian over $R$. Then we have an exact sequence $0 \to M \to \op{Hom}_\mc{T}^\bu(X, Q) \to N \to 0$ where $M =\op{Hom}_\mc{T}^\bu(X, Y_1) / \op{ker} (\op{Hom}_\mc{T}^\bu(X, Y_1) \to \op{Hom}_\mc{T}^\bu(X, Q))$ and $N = \op{im} ( \op{Hom}_\mc{T}^\bu(X, Q) \to \op{Hom}_\mc{T}^\bu(X, Y_2) )$. Observe that $M, N$ are Noetherian over $R$. Hence so is the middle term. 

Suppose that $\op{Hom}_\mc{T}^\bu(X, Y_1) \hookrightarrow \op{Hom}_\mc{T}^\bu(X,Y_2)$, e.g.\ if $Y_1$ is a direct summand of $Y_2$. If $\op{Hom}_\mc{T}^\bu(X,Y_2)$ is Noetherian over $R$, it follows that $\op{Hom}_\mc{T}^\bu(X, Y_1)$ is too. 
\epf

We also need the following lemma, which can be deduced from \Cref{fact:support}. 

\lem\label{corollary:support-generator}
Let $R$ be a commutative Noetherian ring which acts centrally on $\mc{T}$. Fix an object $X \in \mc{T}$. Then for any subset $\mc{U} \sub \op{Spec} R$, the subcategories $\{Y \in \mc{T} \mid \op{Supp}_R\op{Hom}_\mc{T}^\bu(X, Y) \in \mc{U} \}$ and $\{Y \in \mc{T} \mid \op{Supp}_R \op{Hom}_{\mc{T}}^\bu(Y, X) \in \mc{U} \}$ are thick triangulated subcategories. 
\qed
\elem

As a corollary of \Cref{lemma:noetherian-stable} and \Cref{corollary:support-generator}, we get:

\cor\label{corollary:noetherian-support-thick}
Let $R$ be a commutative Noetherian ring which acts centrally on $\mc{T}$. Suppose that $G \in \mc{T}$ is a generator and let $X, Y \in \mc{T}$ be arbitrary objects.
\begin{itemize}
\item If $\op{Hom}_\mc{T}^\bu(G,G)$ is a Noetherian $R$-module, then $\op{Hom}_\mc{T}^\bu(X,Y)$ is also Noetherian.
\item We have $\op{Supp}_R(\op{Hom}_\mc{T}^\bu(X,Y)) \sub \op{Supp}_R(\op{Hom}_\mc{T}^\bu(G,G))$. In particular, all generators have the same support. 
\end{itemize}
\qed
\ecor

\subsubsection{Proof of \Cref{theorem:main-biko-adapt}}
A standard (and possibly the only) tool for obtaining lower bounds on the Rouquier dimension of a triangulated category is the so-called ``ghost lemma". 

\lem[Ghost lemma; see Lem.\ 4.2 in \cite{elagin-lunts-koszul}]\label{new-attempt3.2}
Let $\mc{T}$ be a triangulated category and let $F, G$ be objects. Suppose there exist morphisms
\eq K_c \xrightarrow{\theta_c} K_{c-1} \xrightarrow{\theta_{c-1}} \dots \xrightarrow{\theta_1} K_0
\eeq
such that the following conditions hold:
\begin{enumerate}
\item $\op{Hom}^\bu_\mc{T}(G, \theta_i)= 0$ for each $i=1, \dots,c$; 
\item $\op{Hom}^\bu_\mc{T}(F, \theta_1 \dots \theta_c) \neq 0$.
\end{enumerate}
Then one has $F \notin \langle G \rangle_c$. 
\elem

Following \cite{b-i-k-o}, we recall the notion of a Koszul object. Suppose that $R$ is a commutative ring acting centrally on a triangulated category $\mc{T}$. Given $X \in \mc{T}$ and $r \in R$, we say that $\kos{X}{r}$ is a \emph{Koszul object} for $r$ if it fits into the exact triangle
\eq\label{equation:main-triangle} X \xrightarrow{r} X \to \kos{X}{r} \to \S X. \eeq
Note that $\kos{X}{r}$ is well-defined up to isomorphism. 

By splicing together exact triangles, we get 
\eq\label{equation:spliced-triangle} X \xrightarrow{r} X \to \kos{X}{r} \to \S X \xrightarrow{-r} \S X.\eeq 

More generally, given a sequence of elements $r_1,\dots, r_n \in R$, we inductively define $\kos{X}{\bsr} :=\kos{X}{(r_1,\dots, r_n)}$ as the unique object fitting into the triangle
\eq\label{equation:main-triangle-full} \kos{X}{(r_1,\dots,r_{n-1})} \xrightarrow{r_n} \kos{X}{(r_1,\dots,r_{n-1})} \to \kos{X}{(r_1,\dots,r_n)} \to \S (\kos{X}{(r_1,\dots,r_{n-1})}). \eeq

\lem\label{lemma:hom-to-m}
Suppose that $a_1,\dots,a_n \in R$ is an $M$-regular sequence, where $M= \op{Hom}_\mc{T}^\bu(X,X)$. Then there are natural isomorphisms of $R$-modules:
\begin{enumerate}
\item\label{item:contravariant-isom} $\op{Hom}_\mc{T}^\bu(\kos{X}{(a_1,\dots, a_n)}, X) \simeq M[-n]/(r_1,\dots,r_n);$ 
\item\label{item:covariant-isom} $\op{Hom}_\mc{T}^\bu(X, \kos{X}{(a_1,\dots, a_n)}) \simeq M/(r_1,\dots, r_n)$
\end{enumerate}
\elem
\pf
We only treat (i) because (ii) is similar. Suppose first $n=1$. Apply the functor $\op{Hom}_\mc{T}^\bu(-, X)$ to \eqref{equation:spliced-triangle} (with $r=a_1$) gives an exact sequence of $R$-modules
\eq
M[-1] \xrightarrow{-a_1} M[-1] \to \op{Hom}_{\mc{T}}^\bu(\kos{X}{a_1}, X) \to M \xrightarrow{a_1} M,
\eeq
from which the desired claim follows. In general, suppose that the claim holds for all regular sequences of length at most $n-1$. Consider the spliced triangle 
\eq
\begin{aligned} \kos{X}{(a_1,\dots, a_{n-1})} \xrightarrow{a_n} \kos{X}{(a_1,\dots, a_{n-1})} \to \kos{X}{(a_1,\dots, a_n)} \\ 
 \to \S (\kos{X}{(a_1,\dots, a_{n-1})}) \to  \S (\kos{X}{(a_1,\dots, a_{n-1})}). 
\end{aligned}
\eeq
and apply the same argument. 
\epf

We will also need the following lemma.
\lem[cf.\ Lem.\ 5.11(1) in \cite{benson2008local} and Lem.\ 3.1 in \cite{b-i-k-o}]\label{lemma:vanishing-large-power}
Let $n \geq 1$ be an integer and set $s=2^n$. For any sequence of elements $r_1,\dots, r_n \in R$, we have $r_i^s \cdot \op{Hom}_{\mc{T}}^\bu(\kos{X}{\bsr}, -) = 0 =   \op{Hom}_{\mc{T}}^\bu(-, \kos{X} {\bsr}) \cdot r_i^s$. 
\qed
\elem

\pf[Proof of \Cref{theorem:main-biko-adapt}]
By assumption, $N:= \op{Hom}_\mc{T}^\bu(G,G)$ is a Noetherian $R$-module. Let $n = \op{dim}_R(N)$. If $G'$ is another generator, then it follows from \Cref{corollary:noetherian-support-thick} that $\op{Hom}_\mc{T}^\bu(G', G')$ also a Noetherian $R$-module of Krull dimension $n$. Hence, we can restrict our attention to $G$. Our goal is to construct an object which is not contained in $\langle G \rangle_n$. 

By \Cref{proposition:exists-reg-sequence}, there exists an $N$-regular sequence $r_1,\dots,r_n$. Set $s=2^n$ and note that $r_1^s,\dots,r_n^s$ is still a regular sequence (\Cref{fact:regular-sequence}). For $l \leq n$, consider the exact triangles (obtained by rotating \eqref{equation:main-triangle-full} and setting $X=G$)
\eq  \S^{-1} (\kos{G}{(r_1^s,\dots,r_l^s)}) \xrightarrow{-w[-1]} \kos{G}{(r_1^s,\dots,r_{l-1}^s)} \xrightarrow{r_l^s} \kos{G}{(r_1^s,\dots,r_{l-1}^s)} \to \kos{G}{(r_1^s,\dots,r_l^s)}, \eeq
where $w$ is the connecting map. Suspending $l-1$ more times, we get 
\eq\label{equation:theta-def} 
\begin{aligned}
\S^{-l} (\kos{G}{(r_1^s,\dots,r_l^s)}) \xrightarrow{\theta_l} \S^{-(l-1)}(\kos{G}{(r_1^s,\dots,r_{l-1}^s)})  \\
\xrightarrow{\pm r_l^s} \S^{-(l-1)} (\kos{G}{(r_1^s,\dots,r_{l-1}^s)}) \to \S^{-(l-1)} (\kos{G}{(r_1^s,\dots,r_l^s)}). 
\end{aligned}
\eeq
(note that each application of $\S^{-1}$ has the effect of changing the sign of $\cdot r_l^s$; here $\theta_l$ is defined as the morphism obtained by applying $\S^{-(l-1)}(-)$ to $-w[-1]$).

Applying the functor $\op{Hom}_\mc{T}^\bu(G, -)$, it follows from \Cref{lemma:hom-to-m}\eqref{item:covariant-isom} that we have a natural isomorphism $\op{Hom}_\mc{T}^\bu(G, \S^{-(l-1)}(\kos{G}{r_1^s,\dots, r_{l-1}^s})) \simeq N/(r_1^s,\dots, r_{l-1}^s)[-(l-1)]$. By regularity, $\pm r_l^s$ acts injectively, and it follows that $\op{Hom}_\mc{T}^\bu(G, \theta_l)=0$ by \eqref{equation:theta-def}. 

Next, it follows from \Cref{lemma:vanishing-large-power} that $\pm r_l^s$ acts by zero on $\op{Hom}(\kos{G}{\bsr}, \S^{-(l-1)}(\kos{G}{(r_1^s, \dots, r_{l-1}^s)})$, which implies that $\op{Hom}_\mc{T}^\bu(\kos{G}{\bsr}, \theta_l)$ is surjective by \eqref{equation:theta-def}. But \Cref{lemma:hom-to-m}\eqref{item:contravariant-isom} implies that $\op{Hom}_\mc{T}^\bu(\kos{G}{\bsr}, G) \simeq N[-n]/(r_1,\dots,r_n)$. By regularity, $N[-n]/(r_1,\dots,r_n) \neq 0$. Since this module is the target of the surjective map $\op{Hom}_\mc{T}^\bu(G, \theta_1 \dots \theta_n)$, this map is nonzero. We conclude by \Cref{new-attempt3.2} (with $\kos{G}{\bsr}$ in place of $F$) that $\kos{G}{\bsr} \notin \langle G \rangle_n$. 
\epf

\subsection{Central actions on Fukaya categories}

We now specialize the discussion of the previous sections to wrapped Fukaya categories. For simplicity, we assume that the wrapped Fukaya category is defined with field coefficients and $\Z/2$-gradings.

Given a Liouville manifold $(X, \l)$, we wish to consider the canonical central action of $HH^\bu(\mc{W}(X), \mc{W}(X))$ on the triangulated category $H^0(\op{Perf} \mc{W}(X))$ furnished by the characteristic morphism. The following lemma states that this characteristic morphism admits a purely symplectic topological interpretation, in terms of the familiar closed-open string map which is constructed for instance in \cite[Sec.\ 5.4]{ganatra}. 

\lem\label{lemma:oc-diagram}
For any object $K \in \mc{W}(X)$, the following diagram of unital rings commutes
\eq
\begin{tikzcd}
SH^\bu(X) \ar{r}{\mc{CO}} \ar[rd, "\mc{CO}(-)^0_K(1)" '] & HH^\bu(\mc{W}(X), \mc{W}(X)) \ar[d] \ar[d] & HH^\bu(\op{Perf} \mc{W}(X), \op{Perf} \mc{W}(X)) \ar[l, "\simeq" '] \ar[dl] \\
                                  &  HW^\bu(K, K)                                   &
\end{tikzcd}
\eeq
\elem
\pf
The fact that $\mc{CO}$ is a ring morphism is proved in \cite[Prop.\ 5.3]{ganatra}, and the left triangle is then tautologically commutative. The top right arrow is induced by the canonical (Yoneda) embedding $\mc{W}(X) \to \op{Perf} \mc{W}(X)$, and the left triangle is also tautologically commutative; cf.\ \cite[Lem.\ 2.11]{ritter-smith}. The fact that the top right arrow is an isomorphism of graded commutative rings is a manifestation of the familiar Morita invariance of Hochschild cohomology (see e.g.\ \cite[Rmk.\ 9.5]{seidel-cat-dynamics}). 
\epf


\cor\label{corollary:main-central-bound}
Let $R \sub SH^0(X)$ be a subring which is of finite type as an algebra over $k$. Suppose that there exists a (split-)generator $K \in \mc{W}(X)$ such that $\Hom_{\mc{T}}^\bu(K,K)$ is a Noetherian $R$-module, where $R$ acts via the closed-open map. Then 
\eq \rdim \mc{W}(X) \geq \dim_R \Hom_{\mc{T}}^\bu(K,K). \eeq
\ecor

\pf
Combine \Cref{theorem:main-biko-adapt} and \Cref{lemma:oc-diagram}.
\epf

We will see in the next section how to perform concrete computations via string topology. For now, let us record the following rather trivial example.  
\ex\label{example:computation-n-torus}
If $X= T^*\bb{T}^n$ and $F$ is a cotangent fiber, then we have $\mc{CO}: SH^\bu(X)= \Lambda^\bu[v_1,\dots, v_n] \otimes k[x_1^\pm, \dots, x_n^\pm] \to k[x_1^\pm, \dots, x_n^\pm] = HW^\bu(F, F)$, where $|v_i|=1, |x_i|=0$ and $v$ are in degree $1$ and $\mc{CO}(1 \otimes x_i)= x_i$.  Hence (taking $R=k[x_1^\pm, \dots, x_n^\pm]$, we find that $\mc{W}(T^*\bb{T}^n) \geq n$, where $\mc{W}(-)$ is defined with $k$-coefficients). 
\eex

For a related (but significantly more complicated) example, see \Cref{example:intro-hom-section}.

\subsection{Computations for cotangent bundles}\label{subsection:string-top}
We decribe lower bounds for the dimension of the wrapped Fukaya category of cotangent bundles. There are two approaches which one can take: 
\begin{enumerate}
\item By Abouzaid's generation result \cite{abouzaid-fiber} (and under mild assumptions on $M$), we have $\mc{W}(T^*M) \simeq \op{Perf}(C_{-\bu}(\Omega M))$. Now we apply \Cref{theorem:main-biko-adapt} directly to $\op{Perf}(C_{-\bu}(\Omega M))$ and we are reduced to a purely topological problem, namely, to describe the action of $HH^\bu(\op{Perf}(C_{-\bu}(\Omega M)))$ on $C_{-\bu}(\Omega M)$ via the characteristic morphism. In particular, there is no need to go through symplectic cohomology and closed-open maps.
\item We apply \Cref{corollary:main-central-bound}. By work of Abbondandolo--Schwarz \cite{abbondandolo-schwarz}, the action of symplectic cohomology on the wrapped Floer cohomology of a fiber can be described in terms of string topology, which can then be computed via algebraic topology.
\end{enumerate}
Both approaches ultimately reduce to the same computations in algebraic topology, so it is entirely a matter of taste which one to take. We opt for the second one here.

We will always assume in this section that the wrapped Fukaya category is defined with $\Q$ coefficients and $\Z/2$-gradings (so in particular, all Liouville manifolds carry the appropriate orientations/grading data). Similarly, all singular homology groups considered in this subsection are assumed to have coefficients in $\Q$ and all graded objects (such as homology/cohomology grousp) are assumed to be $\Z/2$-graded (although the gradings all come from a $\Z$-grading, so the reader is welcome to work with $\Z$-gradings if they prefer).

Given a closed, oriented manifold $M$ of dimension $n$, the (shifted) homology of its free loop space $H_{\bu+n}(\mc{L}M)$ carries the structure of a graded-commutative algebra (in fact a BV algebra). The product is called the \emph{Chas--Sullivan product} and was introduced by Chas and Sullivan in \cite{chas-sullivan}. At an intuitive level, this product is constructed as follows: choose cycles $\s_p \in C_p(\mc{L}M)$ and $\s_q \in C_q(\mc{L}M)$, and consider the intersection $\op{ev}(\s_p) \cap \op{ev}(\s_q)$, where $\op{ev}: \mc{L}M \to M$ is the evaluation map sending a loop $\g: S^1 = \R / \Z \to M$ to $\g(1) \in M$. Along this intersection, we can concatenate loops coming from $\s_p$ with loops from $\s_q$. This defines a $p+q-n$ chain in $\mc{L}M$, which is precisely the desired product of $\s_p$ and $\s_q$. 

Next, it will be useful to recall the construction of the \emph{umkehr} (or ``wrong-way") map in differential topology. Let $\mc{M}$ be a smooth manifold and let $\mc{N} \sub \mc{M}$ be a closed, co-oriented submanifold of codimension $d$. Let $\tau_{\mc{M}/\mc{N}} \to \mc{N}$ be the normal bundle of $\mc{N} \sub \mc{M}$ and let $\mc{U} \sub \mc{M}$ be a tubular neighborhood.  We have a sequence of maps
\eq\label{equation:umkehr}
H_\bu(\mc{M}) \to H_\bu(\mc{M}, \mc{M}-\mc{N}) \to H_\bu(\mc{U}, \mc{U}- \mc{N}) \to H_\bu(\tau_{\mc{M}/\mc{N}}, \tau_{\mc{M}/\mc{N}} -0_{\mc{N}}) \to H_{\bu-d}(\mc{N}),
\eeq
where we are considering singular homology. (The first map is the natural projection, the second one is obtained by excision and the fourth one is the Thom isomorphism.) The composition \eqref{equation:umkehr} is called the \emph{umkehr} map associated to the embedding $\mc{N} \sub \mc{M}$. Note that it is not necessary for $\mc{M}, \mc{N}$ to be finite dimensional in this construction. It works as long as $\mc{N} \sub \mc{M}$ has finite codimension and admits a tubular neighborhood.

We will in fact mostly be interested in the umkehr map associated to the embedding $\O_pM \hookrightarrow \mc{L}M$, where $M$ is a closed manifold of dimension $n<\infty$ with a basepoint $p \in M$.%
%
%

The associated umkehr map 
\eq \mc{I}: H_{\bu+n}(\mc{L}M) \to H_\bu(\O_p M) \eeq 
is called the \emph{intersection morphism}.  It was also first introduced (under a slightly different guise) by Chas and Sullivan \cite{chas-sullivan} and plays an important role in string topology. The construction of the intersection morphism via the umkehr map is found in \cite[Sec.\ 9]{q-string-jems}. 

The connection to symplectic topology arises via the work of Abbondandolo--Schwarz \cite{abbondandolo-schwarz}.

\thm[Abbondandolo--Schwarz, \cite{abbondandolo-schwarz}]\label{theorem:a-s}
Let $M$ be a closed manifold of dimension $n$. Let $F_p \sub T^*M$ be a cotangent fiber. Then the following diagram commutes
\eq\label{equation:a-s-diagram} 
\begin{tikzcd}
H_{-\bu+n}(\mc{L}M) \ar{r}{\simeq} \ar{d}{\mathcal{I}}& SH^\bu(T^*M) \ar{d}{\mc{CO}(-)^0_{F_p}(1)} \\
H_{-\bu}(\O_pM) \ar{r}{\simeq} & HW^\bu(F_p, F_p)
\end{tikzcd}
\eeq
\ethm
\pf
All of this is contained in \cite{abbondandolo-schwarz}. However, since the notation and terminology in \textit{loc}.\ \textit{cit}.\ differs from ours, we give a brief overview for the reader. 

Given a manifold $M$, the wrapped Floer cohomology of a cotangent fiber is denoted by $HF^\Omega(T^*M)$ in \cite{abbondandolo-schwarz} and called ``Floer homology for Hamiltonian orbits with Dirichlet boundary conditions". The intersection morphism that we denote by $\mc{I}$ is denoted by $i_!$. The authors describe how to realize this map in Morse theory in \cite[Sec.\ 2.2]{abbondandolo-schwarz}, which also depends on material from the appendix. The closed-open map is denoted by $I_!$ and constructed in \cite[Sec.\ 3.5]{abbondandolo-schwarz}. The top horizontal arrow in \eqref{equation:a-s-diagram} corresponds to \cite[(3)]{abbondandolo-schwarz}, while the bottom horizontal arrow is \cite[Thm.\ B]{abbondandolo-schwarz}. The fact that \eqref{equation:a-s-diagram} commutes amounts to the assertion that the diagram \cite[(7)]{abbondandolo-schwarz} has a Floer theoretic counterpart; as is often the case in Floer theory, the proof consists in showing that the two ways of traveling around the diagram correspond to the two ends of a cobordism of moduli spaces (see \cite[Thm.\ 4.11]{abbondandolo-schwarz}). 
\epf

We now discuss some concrete computations. We begin by recording the following lemma.

\lem\label{lemma:products-umkehr}
Let $M$ and $N$ be closed oriented manifolds of dimension $m$ and $n$ respectively. Then the following diagram of vector spaces commutes, where the vertical arrow is the usual Eilenberg--Zilber map:
\eq
\begin{tikzcd}
H_{\bu+m}(\mc{L}M) \otimes H_{\bu+n}(\mc{L}N) \ar{r}{\mc{I} \otimes \mc{I}} \ar{d}{\simeq} & H_\bu(\O_*M) \otimes H_\bu(\O_*N) \ar{d}{\simeq}  \\
H_{\bu+m+n}(\mc{L}(M \tms N)) \ar{r}{\mc{I}}  & H_\bu(\O_* (M \tms N)) 
\end{tikzcd}
\eeq
The vertical arrows are in fact isomorphisms of algebras with respect to the Chas--Sullivan and Pontryagin products. 
\end{lemma}
\pf
There are canonical splittings of the based loop and free loop spaces which are compatible with the inclusions (of the former into the latter). It follows from the K\"unneth formula and the definition of the umkehr map that the umkehr map takes products of spaces to tensor products. Finally, it is apparent from the above description of the Chas--Sullivan product that it is compatible with the K\"{u}nneth isomorphism, which means that the vertical left arrow is indeed an isomorphism of algebras. Similar considerations apply for the Pontryagin product. 
\epf

The following theorem is extremely useful, as it essentially reduces the computation of the intersection morphism $\mc{I}$ to rational homotopy theory.

\thm[F\'{e}lix--Thomas--Vigu\'{e}-Poirrier; see Thm.\ F. in \cite{q-string-jems}]\label{theorem:f-t-v-diagram}
Let $M$ be a simply-connected closed oriented manifold of dimension $n$ with basepoint $p$. Then there exists a commutative diagram of morphisms of graded-commutative algebras
\eq\label{equation:hochschild-loop}
\begin{tikzcd}
HH^\bu(C^\bu(M), C^\bu(M)) \ar{d}{HH^\bu(C^\bu(M), \e)} & H_{\bu+n}(\mc{L}M) \ar{l}{\simeq} \ar{d}{\mc{I}} \\
HH^\bu(C^\bu(M), \Q) & H_\bu(\Omega_pM) \ar{l}{\simeq},
\end{tikzcd}
\eeq
where $\e: C^\bu(M) \to \Q$ is the augmentation induced by the inclusion $p \hookrightarrow M$. 
\ethm

Given a graded vector space $V= \{V_i\}_{i\in \Z}$, let $T(V):= \oplus_{i \geq 0} V^{\otimes i}$ be the \emph{tensor algebra} on $V$.  Let $I \sub T(V)$ be the ideal generated by elements of the form $x \otimes y - (-1)^{|x| |y|} y \otimes x$ with $x, y \in V$. Then $\Lambda V:= T(V)/ I$ be the \emph{free commutative graded algebra} on $V$. See \cite[Sec.\ 3(a)]{q-homotopy-book}. If $V$ is generated by a single element $v$, we write $T(V), \Lambda(V)$ as $T(v), \Lambda v$ respectively. 

Let us now record some useful facts drawn from \cite[Sec.\ 6.1]{f-t-v}. First of all, $HH^\bu(C^\bu(S^n), \Q)= T(v)$ with $|v|= n-1$. If $n>1$ is odd, then $HH^\bu(C^\bu(S^n), C^\bu(S^n)) \simeq \Lambda u \otimes T(v)$, where $|u|=-n$ and $|v|= n-1$. We have $\mc{I}= \e \otimes 1: \Lambda u \otimes T(v) \to T(v)$, where we implicitly identify both columns of \eqref{equation:hochschild-loop} under the horizontal isomorphisms.

If $n>1$ is even, then 
\eq HH^\bu(C^\bu(S^n), C^\bu(S^n)) \simeq \Lambda b \otimes k[a,c]/(2ac, a^2, ab)\eeq
with $|a|=-n, |b|=-1$ and $|c|= 2n-2$. We have $\mc{I}(a)=\mc{I}(b)=0$ and $\mc{I}(c)= v^2$. 

\ex[Products of odd dimensional spheres]\label{example:products-odd-spheres}
Let $M= S^{n_1}\tms \dots \tms S^{n_k}$ be a product of odd dimensional spheres with $n_i \geq 3$. Fix a basepoint $p \in M$. If $k=1$, then it follows from the computations above that there is an inclusion $\Q[x] \hookrightarrow H_{\bu + n_1}(\mc{L} S^{n_1})$ sending $x \mapsto 1 \otimes v$, such that the composition $\mc{Q}[x] \to H_{\bu + n_1}(\mc{L} S^{n_1}) \to H_{\bu }(\O S^{n_1})$ is an isomorphism.

For $k>1$, similar considerations along with \Cref{lemma:products-umkehr} furnish maps
\eq \Q[x_1,\dots,x_k] \hookrightarrow H_{\bu + n_1}(\mc{L} M) \to H_{\bu }(\O M) \eeq
where $|x_i|= n_i-1$ and the composition is an isomorphism.

By combining \Cref{theorem:a-s} and \Cref{corollary:main-central-bound}, it follows that 
\eq \rdim \mc{W}(T^*M) \geq k. \eeq

\eex

\ex[Compact Lie groups]
Let $G$ be a simply-connected compact Lie group. Then $G$ has a Sullivan minimal model of the form $(\Lambda V, d=0)$, where $V=\{V^i\}_{i \in \N}$ is a graded $\Q$-vector space supported in odd degrees and satisfying $V^1=0$. Thus $G$ has the rational homotopy type of a product of odd dimensional spheres.   

The dimension of $V$ is precisely the rank of the group $G$ \cite[Thm.\ 3.33]{f-o-t}. It follows from \Cref{example:products-odd-spheres} that we have
\eq \rdim \mc{W}(T^*G) \geq \op{rank} G. \eeq
\eex

For the reader's convenience, we list the simply-connected simple compact Lie groups in \Cref{figure:lie-groups}.  It follows from the classification theorem that all other simply-connected compact Lie groups are obtained from these by taking products. For $i=6,7,8$, we write $\tilde{E}_i$ to denote the universal cover of the exceptional Lie group $E_i$.

\begin{figure}[h]
\begin{tabular}{ |c||c|c|c| } 
 \hline
 Simply-connected compact Lie group & Real dimension & Rank & $\rdim \mc{W}(T^*G)$  \\ 
 \hline 
$\op{Sp}(n), n \geq 1$ & $n(2n+1)$  & $n$ & $\geq n$ \\ 
$\op{SU}(n), n \geq 3$ & $n^2-1$ & $n-1$ & $\geq n-1$ \\ 
$\op{Spin}(n), n \geq 7$ & $n(n-1)/2$ & $\lfloor n/2 \rfloor$ & $\geq \lfloor n/2 \rfloor$ \\ 
$G_2$ & $14$ & 2 & $\geq 2$ \\ 
$F_4$ & $52$ & $4$ & $\geq 4$ \\ 
$\tilde{E}_6$ & $156$ & $6$ & $\geq 6$ \\ 
$\tilde{E}_7$ & $266$ & $7$ & $\geq 7$ \\ 
$\tilde{E}_8$ & $496$ & $8$ & $\geq 8$ \\ 
 \hline
\end{tabular}
\caption{The simply-connected compact Lie groups}
\label{figure:lie-groups}
\end{figure}


\ex[Products of arbitrary spheres]\label{example:products-all-spheres}
We again assume $M= S^{n_1} \tms \dots \tms S^{n_k}$. Suppose first that $k=1$ and $n_1>1$ is even.  Then there is a natural inclusion $\Q[x] \hookrightarrow \L b \otimes \Q[a,c]/ (2ac, a^2, ab)$ with $x \mapsto c$. Composing with the intersection morphism gives the map $\Q[c] \to \Q[v]$ sending $c \mapsto v^2$. More generally, we can apply \Cref{lemma:products-umkehr}, to get a map 
\eq \Q[x_1,\dots, x_k] \to  H_{\bu + n_1}(\mc{L} M)  \to \Q[v_1,\dots,v_k], \eeq 
where 
\begin{itemize}
\item $|x_i|= n_i-1$ and $x_i \mapsto v_i$ if $n_i$ is odd;
\item $|x_i|= 2n_i-2$ and $x_i \mapsto v_i^2$ if $n_i$ is even.
\end{itemize}
Observe that $\Q[v_1,\dots,v_k]$ is a Noetherian module over $\Q[x_1,\dots, x_k]$ (indeed, it is finitely generated as a module over a Noetherian ring).  By combining \Cref{theorem:a-s} and \Cref{corollary:main-central-bound}, it follows that 
\eq \rdim \mc{W}(T^*M) \geq k. \eeq
\eex


\ex[Complex projective spaces]
Let us also briefly mention the case of complex projective spaces. The cohomology of $\O \bb{CP}^n$ is computed in \cite[Thm.\ 1.3]{seeliger}. It contains a free subalgebra $\Q[y]$ where $|y|=2n$. On the other hand, as explained in e.g.\ \cite[Sec.\ 3]{seeliger}, we have $\O \bb{CP}^n= S^1 \tms \O S^{2n+1}$ and hence $H_\bu(\O \bb{CP}^n) = H_\bu(S^1) \otimes H_\bu(\O S^{2n+1})$. It is shown in \cite[Ex.\ 5.7]{briggs-gelinas}  that the characteristic morphism hits the class $1 \otimes v$ of degree $2n$ in $H_\bu(\O \bb{CP}^n)$. For degree reasons, the composition of $\Q[y] \to H_{\bu+2n}(\mc{L}\bb{CP}^n) \to H_\bu(\O \bb{CP}^n)$ must send $y$ to $1 \otimes v$. 

By a similar argument to \Cref{example:products-odd-spheres}, we find that if $M=M_1 \tms \dots \tms M_k$ where each $M_i$ is a sphere of dimension at least $2$ or a complex projective space of arbitrary dimension, then 
\eq \rdim \mc{W}(M) \geq k. \eeq
\eex

\section{Applications to symplectic geometry}

We now describe how the Rouquier dimension of wrapped Fukaya categories is related to certain quantities and operations of geometric interest. 

\subsection{Embedding monotonicity}

An important fact about the Rouquier dimension of Weinstein manifolds is that it is monotonic under Liouville embeddings which are Weinstein up to deformation. To explain this precisely, we begin with the following definition.

\defi\label{definition:liouville-embedding}
A \emph{Liouville embedding} of Liouville domains is a codimension $0$ smooth embedding 
\eq i: (X_0^{\op{in}}, \l_0^{\op{in}}) \hookrightarrow (X_0, \l_0) \eeq
where $i^*\l_0= e^{\rho} \l_0^{\op{in}} + df$ for some $\rho \in \R$ and $f: X_0^{\op{in}} \to \R$.
\edefi

We refer to \cite[Sec.\ 2.1.3]{ganatra-siegel} for a detailed discussion of Liouville embeddings. As explained in \cite[Lem.\ 2.2]{ganatra-siegel}, any Liouville embedding as in \Cref{definition:liouville-embedding} extends to a family of Liouville embeddings $i^t: (X_0^{\op{in}}, \l_0^{\op{in}}) \hookrightarrow (X_0, \l_0^t)$ such that $i=i^0$ and $i^1$ is a strict Liouville embedding, i.e.\ $(i^1)^* \l_0^1 = \l_0^{\op{in}}$. 

\prop\label{proposition:embeddings}
Let $i: (X_0^{\op{in}}, \l^{\op{in}} ) \sub (X_0, \l)$ be a Liouville embedding and suppose that $(X_0^{\op{in}}, \l^{\op{in}})$ and $(X_0- i(X_0^{\op{in}}), \l)$ are Weinstein up to deformation. Then 
\eq \rdim \mc{W}(X^{\op{in}})  \leq \rdim \mc{W}(X) \eeq
where $X^{\op{in}}, X$ are the completions of $X_0^{\op{in}}, X_0$ respectively.
\eprop
\pf 
Since the wrapped Fukaya category is invariant under deformation, we may assume that the embedding is strict, i.e. $i^*\l= \l^{\op{in}}$. We may therefore appeal to \cite[Prop.\ 11.2]{gps2}, which constructs a Viterbo restriction functor 
\eq \op{Tw} \mc{W}(X) \to \op{Tw} \mc{W}(X^{\op{in}}). \eeq
As explained in \cite[Prop.\ 11.2]{gps2}, the hypothesis that the $X_0^{\op{in}}$ and $X_0- i(X_0^{\op{in}})$ are Weinstein implies that there is in fact a quasi-equivalence $\op{Tw} \mc{W}(X) /\mc{D} \to \op{Tw} \mc{W}(X^{\op{in}})$, where $\mc{D}$ is a set of objects. The claim now follows from \Cref{corollary:quot-rdim} upon passing to $\op{Perf}(-)$. 
\epf

In the statement of \Cref{proposition:embeddings}, note that the orientation/grading data on $\mc{W}(X^{\op{in}})$ is implicitly induced from the orientation/grading data on $X$. 

It is natural to ask whether one can drop the hypothesis that $X_0- i(X_0^{\op{in}})$ is Weinstein \Cref{proposition:embeddings}. We note that Sylvan has proved \cite{sylvan} that the Viterbo restriction functor in this setting functor is a \emph{homological epimorphism}, although it is not clear if this helps.

By combining various upper and lower bounds for $\rdim\mc{W}(-)$ discussed in the introduction and below, it is easy to write down examples of putative embeddings $X_0^{\op{in}} \hookrightarrow X_0$ which are obstructed by $\rdim \mc{W}(-)$. For instance:
\ex\label{example:contrived-embedding}
Let $X$ be a Weinstein manifold of real dimension $2(n^2-1)$ admitting a Lefschetz fibration with Weinstein fibers having $<n-2$ critical points. Then $X$ does not contain $T^*SU(n)$ as a Weinstein subdomain.
\eex
There is of course nothing special about \Cref{example:contrived-embedding}, and there is a plethora of similar examples. It seems unlikely to the authors that examples such as these could be handled using invariants already in the literature. From this perspective, $\rdim \mc{W}(-)$ is a useful addition to the symplectic topologist's toolkit of embedding obstructions. On the other hand, $\rdim \mc{W}(-)$ has many blind spots. The most obvious defect of $\rdim \mc{W}(-)$ as an embedding obstruction is that it outputs a very limited range of values: indeed, it seems plausible that $\rdim \mc{W}(X^{2n}) \in [0, n]$ (see \Cref{cor_rouq_upper} for a somewhat weaker statement).



\subsection{Lefschetz fibrations}

We now discuss upper bounds for the Rouquier dimension of wrapped Fukaya categories in terms of the number of critical points of a Lefschetz fibration. 

\defi[see Sec.\ 2 in \cite{maydanskiy-seidel}]
Given a Liouville manifold $X$, a \emph{Lefschetz fibration} is a smooth map $f: X\to \C$ having Morse-type critical points and such that $\op{ker} df$ is symplectic away from the critical points. This data is required to satisfy some additional technical conditions at infinity detailed in \cite[Sec.\ 2]{maydanskiy-seidel}. A Lefschetz fibration is said to have \emph{Weinstein fibers} if the general fiber is a Weinstein manifold up to deformation. 
\edefi

There are various essentially equivalent notions of a Lefschetz fibration in the literature. We have adopted the setup of \cite[Sec.\ 2]{maydanskiy-seidel} for consistency with \cite[Sec.\ 1.1]{gps1}. 

\defi
Given a Liouville manifold $X$ (considered up to Liouville homotopy), we let 
\eq \op{Lef}_w(X) \in \N \cup \{\infty\} \eeq 
denote the minimal number of critical points of a Lefschetz fibration on $X$ with Weinstein fibers. We set $\op{Lef}_w(X)= \infty$ if no such fibration exists.
\edefi

It was proved by Giroux--Pardon \cite[Thm.\ 1.10]{giroux-pardon} using quantitative transversality techniques that any Weinstein manifold admits a Lefschetz fibration with Weinstein fibers. Conversely, if a Liouville manifold admits a Lefschetz fibration with Weinstein fibers, then it is Weinstein (this is tautological if one uses \cite[Def.\ 1.9]{giroux-pardon}). Hence $\op{Lef}_w(X) < \infty$ iff $X$ admits a Weinstein structure. 

Concerning lower bounds, there is a purely topological bound $\op{rk} H_n(X; \Z) \leq \op{Lef}_w(X)$, where $X$ has real dimension $2n$. It is also immediate that $\op{Lef}_w(X)=0$ iff $X$ is subcritical, meaning that it splits as a product of Liouville manifolds $(X, \l)= (F_0 \tms \C, \l_0 \tms \l_{\C})$. 

The following proposition, which was explained to us by Y.\ Bar\i \c{s} Kartal, states that the Rouquier dimension gives another lower bound.

\prop\label{prop:lef-rdim}
Let $X$ be a Liouville manifold. Then 
\eq\label{equation:lefschetz-bound}  \op{Lef}_w(X) \geq \rdim{\mc{W}(X)} +1. \eeq
\eprop

\pf
We may assume without loss of generality that $0<\op{Lef}_w(X) = k <\infty$. Fix a Lefschetz fibration $f: X \to \C$ with $k$ critical points and choose a set of thimbles $\{ L_1,\dots,L_k \}$. These form an exceptional collection for the partially wrapped Fukaya category $\mc{W}(X, f^{-1}({\infty}))$; see \cite[Ex.\ 1.4]{gps1}. According to \cite[Cor.\ 1.17]{gps2}, this exceptional collection is \emph{full}, meaning that it generates $\mc{W}(X, f^{-1}(\{\infty\}))$. It follows by \Cref{lemma:sod-rdim} that $\rdim \mc{W}(X, f^{-1}(\{\infty\})) \leq k-1$. 

It follows from \cite[Cor.\ 3.9]{gps2} that $\mc{W}(X, f^{-1}(\{\infty\})) = \mc{W}(X, \fk{c})$, where $\fk{c}$ is the skeleton of $f^{-1}(\infty)$. The assumption that $f^{-1}(\{\infty\})$ is Weinstein implies that $\fk{c}$ is mostly Legendrian \cite[see Def.\ 1.7]{gps2}. It then follows from \Cref{theorem:stop-removal}, that the natural map $\mc{W}(X, f^{-1}(\{\infty\})) \to \mc{W}(X)$ is a quotient of $A_\infty$ categories.  \Cref{corollary:quot-rdim} thus implies that $\rdim \mc{W}(X)  \leq k-1$.
\epf

We do not expect \eqref{equation:lefschetz-bound} to be sharp except in very special cases. 

\ex\label{example:cotangent-sphere-lefschetz}
The affine quadric $\{z_1^2 + \dots + z_n^2 +z_{n+1}^2 = 1\}$ with the Liouville structure $\sum y_i dx_i$ is isomorphic to $T^*S^n$ as a Liouville manifold. The projection $(z_1,\dots, z_{n+1}) \mapsto z_{n+1}$ is a Lefschetz fibration with two critical points. Hence $\rdim \mc{W}(T^*S^n) \leq 1$. By \Cref{example:products-all-spheres}, this is an equality if $\mc{W}(T^*S^n)$ is defined with $\Z/2$-gradings and $\Q$-coefficients.
\eex

Although the definition of $\op{Lef}_w(X)$ only makes sense on Liouville manifolds, it is not a priori clear to what extent this invariant depends on the Liouville structure (recall Lazarev's result ``$\op{WMor}(X)=\op{Mor}(X)$" discussed in the introduction, which implies that $\op{WMor}(X)$ does not depend on the Liouville structure). The following corollary shows that the natural complex analog of ``$\op{WMor}(X)=\op{Mor}(X)$" is false.
\cor\label{corollary:lef-distinguish}
There exists a Weinstein manifold $T^*S^3_{\op{exotic}}$ which is formally isotopic to $T^*S^3$, but such that $\op{Lef}_w(T^* S^3) \neq \op{Lef}_w(T^* S^3_{\op{exotic}})$.
\ecor
\pf
\cite[Thm.\ 4.7]{e-g-l} produces a Weinstein manifold $T^*S^3_{\op{exotic}}$ which contains a \emph{flexible} Lagrangins $L \simeq \bb{T}^3$. This means by definition that L is \emph{regular} (see \cite[Def.\ 3.1]{e-g-l}), and hence the cobordism $T^*S^3_{\op{exotic}} - T^*_{\leq \e}L$ is Weinstein. Now combine \Cref{proposition:embeddings} and \Cref{prop:lef-rdim}. It follows that $\op{Lef}_w(T^*S^3_{\op{exotic}}) \geq \rdim \mc{W}(T^*S^3_{\op{exotic}}) +1\geq \rdim \mc{W}(T^*\bb{T}^3)+1 \geq 4$ by \Cref{example:computation-n-torus}. On the other hand, we saw in \Cref{example:cotangent-sphere-lefschetz} that $\op{Lef}_w(T^*S^3) \leq 2$. 
\epf

\subsection{Rigidity of skeleta}

Finally, we prove the quantitative symplectic rigidity result for the skeleton of a Weinstein manifold which was advertised in the introduction. 

\prop\label{prop_int_rdim}
Let $(X, \l)$ be a Liouville manifold.  Suppose that the skeleton $\fk{c}_X$ is mostly Lagrangian and that every component of $\fk{c}_X^{\op{crit}}$ admits a generalized cocore.  

Let $\phi: X \to X$ be a Hamiltonian symplectomorphism. If $\fk{c}_X \cap \phi(\fk{c}_X)= \fk{c}_X^{\op{crit}} \cap \phi(\fk{c}_X^{\op{crit}})$ and these intersection points are transverse (a generic condition on $\phi$), then 
\eq\label{equation:rdim-lower-bound-intersection}
 |\fk{c}_X \cap \phi(\fk{c}_X)| \geq \rdim{\mc{W}(X)} +1.
\eeq
\eprop

\pf
We consider the Liouville manifold $(\ov{X} \tms X, -\l \oplus \phi_* \l)$ which has skeleton $\fk{c}_X \tms \phi(\fk{c}_X)$. Let $\Delta \sub \ov{X} \tms X$ be the diagonal. By assumption, $\Delta$ intersects $\fk{c}_X \tms \phi(\fk{c}_X)$ at finitely many points. Suppose there are $r \geq 0$ such points.

Consider the diffeomorphism 
\begin{align} (\ov{X} \tms X, -\l \oplus \phi_* \l) &\to (\ov{X} \tms X, -\l \oplus \l) \\
(x_1, x_2) &\mapsto (x_1, \phi^{-1}(x_2)),
\end{align}
which sends $\Delta$ to the Lagrangian $L= \{(x, \phi^{-1}(x)\}$. 

We now apply \Cref{theorem:gps-intersections-enhanced}. This produces a resolution of $L$ of length $r \geq 0$ by generalized cocores of $\fk{c}_{\ov{X} \tms X}= \fk{c}_X \tms \fk{c}_X$. Observe that one can choose these generalized cocores to be products of generalized cocores of $\fk{c}_X$. Since $L$ is isotopic to $\Delta$ through a compactly-supported Lagrangian isotopy, they define isomorphic objects in $\mc{W}(\ov{X} \tms X)$.  Applying \Cref{theorem:bimod-functor}, we obtain a resolution of the diagonal bimodule of length $r$ by Yoneda bimodules. The conclusion now follows from \Cref{lemma:resolution-diag}. 
\epf

\rmk
The assumption that $\fk{c}_X^{\op{crit}}$ and $\phi(\fk{c}_X^{\op{crit}})$ intersect transversally is used in the proof of \Cref{prop_int_rdim} when appealing to \Cref{theorem:gps-intersections-enhanced}. However, we do not know whether this assumption is necessary for the conclusion of the theorem to hold.
\ermk

The proof of  \Cref{prop_int_rdim} can actually be adapted to derive lower bounds for $ |\fk{c}_X \cap \phi(\fk{c}_X)|$ which do not involve the Rouquier dimension of $\mc{W}(X)$.

For example, suppose that $\mc{W}(X)$ admits a nonzero proper module. Then we claim $ |\fk{c}_X \cap \phi(\fk{c}_X)|>1$. Indeed, if $r=0$, then $\mc{W}(X)=0$. If $r=1$, then the proof of \Cref{prop_int_rdim} implies that the diagonal bimodule is isomorphic to $\mc{Y}^l_K \otimes_k \mc{Y}^r_L$ for some objects $K, L$ of $\mc{W}(X)$. Convolving with $\mc{P}$, we find that $H^\bu(P, P)= H^\bu(K, K) \otimes H^\bu(\mc{Y}^r_L, P)$. Hence $H^\bu(K, K)$ is finite-dimensional. On the other hand, the convolution argument of \Cref{lemma:resolution-diag} proves that $K$ split-generates $\mc{W}(X)$. Thus $\mc{W}(X)$ is proper, which contradicts the main result of \cite{ganatra-non-proper}.

Similarly, one can show that $ |\fk{c}_X \cap \phi(\fk{c}_X)|$ is strictly greater than the dimension of the image of $K_0(\op{Prop} \mc{C})$ in $K_0(\op{Perf} \mc{C})$.

\ex
If $X$ contains an orientable closed exact Lagrangian, then this defines a proper module over $\mc{W}(X)$ (where we work with $\Z/2$ gradings and $\Z/2$ coefficients). Hence $ |\fk{c}_X \cap \phi(\fk{c}_X)| \geq 2$. 
\eex

\section{Appendix}\label{section:appendix}
The purpose of this appendix is to prove \Cref{theorem:bimod-functor}. Our proof is essentially a modification of \cite{gps2} and \cite[Sec.\ 9]{ganatra}, but we include it for completeness. 

\rmk
After the original version of this paper first appeared, Ganatra--Pardon--Shende posted a substantially revised version of \cite{gps2} from which a proof of \Cref{theorem:bimod-functor} can straightforwardly be extracted. 
\ermk

Let $(X, \l)$ be a Liouville manifold and let $\mc{W}(X)$ be its wrapped Fukaya category. Let $\mc{W}(\ov{X} \tms X)$ denote the wrapped Fukaya category of the product $(\ov{X} \tms X, -\l \oplus \l)$. 

If $K$ and $L$ are cylindrical Lagrangians in $X$, their product $K \times L$ is not a cylindrical Lagrangian in $\ov{X} \times X$ in general. However, following \cite[Section 7.2]{gps2}, the product can be ``cylindrized". This results in a cylindrical Lagrangian $K \tilde{\times} L$  which behaves like $K \times L$ in the Floer-theoretic sense; see Lemma \ref{lem_cyl}. If $K, L$ are invariant under the Liouville flow (e.g.\ they could be cocores), then $K \tms L = K \tilde{\tms} L$. 

Let $\mc{W}_{\op{cyl}}(\ov{X} \tms X) \subset \mc{W}(\ov{X} \tms X)$ denote the smallest full, replete subcategory of $\mc{W}(\ov{X} \tms X)$ containing the diagonal $\Delta \subset \ov{X} \tms X$ and all cylindrized products (i.e.\ all objects of the form $K \tilde{\tms} L$).


\thm\label{theorem:bimod-functor}
There exists a fully faithful $A_\infty$ functor 
\eq \mc{F}: \mc{W}_{\op{cyl}}(\ov{X} \tms X) \to (\mc{W}(X), \mc{W}(X))-\op{mod} \eeq which satisfies the following properties: 
\begin{itemize}
\item the image of the diagonal $\Delta \sub \ov{X} \tms X$ is the diagonal bimodule;
\item the image of any product Lagrangian $K \tms L$ is the Yoneda bimodule $\mc{Y}^l_K \otimes_k \mc{Y}^r_L$.
\end{itemize}
\ethm

\cor\label{corollary:appendix-weinstein}
If $X$ is Weinstein, then $\mc{F}$ extends to a fully-faithful functor $\op{Tw}\mc{W}(\ov{X} \tms X) \to (\mc{W}(X), \mc{W}(X))-\op{mod}.$
\qed
\ecor

The proof of \Cref{theorem:bimod-functor} will occupy the remainder of this appendix.

\defi[Wrapping category of $(X, \l)$]
Define the \emph{wrapping category} $((X, \l) \leadsto -)^+$ of $(X, \l)$ as follows. The objects of $((X, \l) \leadsto -)^+$ are Hamiltonian isotopies of $X$ generated by Hamiltonians $\{ H_t \}_{0 \leq t \leq 1}$ such that $H_t$ is linear at infinity.  Given two objects $\{ H_t \}_{0 \leq t \leq 1}$ and $\{ H^{\prime}_t \}_{0 \leq t \leq 1}$ in $((X, \l) \leadsto -)^+$, denote by $\phi_{H_t}$ and $\phi_{H^{\prime}_t}$ respectively the Hamiltonian diffeomorphisms generated by them. Morphisms from $\{ H_t \}_{0 \leq t \leq 1}$ to $\{ H^{\prime}_t \}_{0 \leq t \leq 1}$ are given by Hamiltonian diffeomorphisms satisfying $\phi_{H^{\prime}_t} = \phi_{H_t} \sharp \phi_{\tilde{H}_t}$ and the generating Hamiltonians $\{ \tilde{H}_t \}_{0 \leq t \leq 1}$ are \emph{positive} and linear at infinity. Here $\phi_{H_t} \sharp \phi_{\tilde{H}_t}$ means the Hamiltonian diffeomorphism generated by $H_t \sharp \tilde{H}_t (x) = H_t(x) + \tilde{H}_t((\Phi^t_{H_t})^{-1}(x))$ and $\Phi^t_{H_t}$ is the time-$t$ flow of the Hamiltonian vector field of $H_t$.
\edefi

Recall that a category $\mc{C}$ is \emph{filtered} if
\begin{enumerate}
\item $\mc{C}$ is non-empty;
\item for every pair of objects $x, y \in \mc{C}$, there exists an object $z \in \mc{C}$ and morphisms $x \rightarrow z, y \rightarrow z$;
\item for every pair of morphisms $f,g: x \rightarrow y$, there exists an object $z \in \mc{C}$ and a morphism $h: y \rightarrow z$ such that $h \circ f = h \circ g$.
\end{enumerate}
A filtered category $\mc{C}$ is has \emph{countable cofinality} if there exists a cofinal functor $\mathbb{Z}_{\geq 0} \rightarrow \mc{C}$.

\lem
The category $((X, \l) \leadsto -)^+$ is filtered and has countable cofinality.
\elem
\pf
The proof that the category is filtered proceeds exactly as \cite[Lem.\ 3.27]{gps1}, except replacing the Lagrangian isotopies by Hamiltonian isotopies. That $((X, \l) \leadsto -)^+$ has countable cofinality can be proved by choosing a countable family of linear Hamiltonians with slope going to infinity.
\epf

\lem\label{lem_app_diag}
Consider the diagonal $\Delta \subset \ov{X} \times X$ as a cylindrical Lagrangian. Let $(\Delta \leadsto -)^+$ be the positive wrapping category of $\Delta$ in the sense of \cite[Sec.\ 3.4]{gps1}. Denote by $\pi_1$ and $\pi_2$  the projection from $\ov{X} \times X$ to the first and second factors respectively. Then the functor from $((X, \l) \leadsto -)^+$ to $(\Delta \leadsto -)^+$ by taking a generating Hamiltonian $H_{t}$ to the Lagrangian isotopy of $\Delta$ which is generated by $\pi_1^{*}H_{1-t/2} + \pi_2^{*}H_{t/2}$ is cofinal.
\elem
\pf
This holds by inspecting the definitions.
\epf

We introduce a category $\widetilde{\mc{O}}^{\pd}_{\ov{X} \times X}$ following \cite[Sec.\ 8.2]{gps2}. Fix a countable collection of Lagrangians $\tilde{I}^{\pd}_{\ov{X} \times X}$ in $\ov{X} \times X$ consisting of $\Delta$ and product Lagrangians $K \times L \subset \ov{X} \times X$ which represent every pair of isotopy classes. For the diagonal $\Delta$, choose a countable cofinal family $\Delta = \Delta^{(0)} \leadsto \Delta^{(1)} \leadsto \cdots$ in $(\Delta \leadsto -)^+$ induced from $((X, \l) \leadsto -)^+$ as in Lemma \ref{lem_app_diag}. Also choose cofinal wrappings $K = K^{(0)} \leadsto K^{(1)} \leadsto \cdots$ and $L = L^{(0)} \leadsto L^{(1)} \leadsto \cdots$. Then the category $\widetilde{\mc{O}}^{\pd}_{\ov{X} \times X}$ has objects given by $\mathbb{Z}_{\geq 0} \times \tilde{I}^{\pd}_{\ov{X} \times X}$, which are either of the form $\Delta^{(i)}$ or $K^{(i)} \times L^{(i)}$. The objects of $\widetilde{\mc{O}}^{\pd}_{\ov{X} \times X}$ admits a total order induced from $\mathbb{Z}_{\geq 0}$.

Note that by choosing the wrappings generically, we can assume that for any pair of Lagrangians in $\widetilde{\mc{O}}^{\pd}_{\ov{X} \times X}$, they intersect transversally and are disjoint at infinity with distance uniformly bounded from below. Using generic product almost complex structures, as they induce uniformly bounded geometry, one can define $A_{\infty}$ operations on $\widetilde{\mc{O}}^{\pd}_{\ov{X} \times X}$ by counting suitable holomorphic discs based on the compactness arguments in \cite[Prop.\ 3.19]{gps1}. Of course, the morphism spaces are given by the Floer cochain complexes. The wrappings define a collection of continuation elements \cite[Def.\ 3.25]{gps1} $C_{\pd}$ inside the Floer cohomology groups $HF^{0}(K^{(i+1)} \times L^{(i+1)}, K^{(i)} \times L^{(i)})$ and $HF^{0}(\Delta^{(i+1)}, \Delta^{(i)})$. Then we define an $A_{\infty}$ category by localizing $\widetilde{\mc{O}}^{\pd}_{\ov{X} \times X}$ along $C_{\pd}$: 
\eq
\widetilde{\mc{W}}^{\pd}(\ov{X} \times X) := \widetilde{\mc{O}}^{\pd}_{\ov{X} \times X}[C^{-1}_{\pd}].
\eeq

By construction and \cite[Lem.\ 3.37]{gps1}, the following isomorphisms hold in $\widetilde{\mc{W}}^{\pd}(\ov{X} \times X)$ (this is parallel to the computations in \cite[Sec. 8.1]{ganatra} under the ``quadratic" setup):
\eq
\op{hom}^{\bu}(\Delta, \Delta) \cong SH^{\bu}(X),
\eeq
\eq\label{eqn_diag_floer}
\op{hom}^{\bu}(K \times L, \Delta) \cong HW^{\bu}(L, K), \op{hom}^{\bu}(\Delta, K \times L) \cong HW^{\bu}(K, L), 
\eeq
\eq\label{eqn_prod_floer}
\op{hom}^{\bu}(K \times L, K' \times L') \cong HW^{\bu}(K, K') \otimes HW^{\bu}(L, L').
\eeq

\prop\label{prop_prod_mod}
There exists a fully faithful $A_{\infty}$ functor $\mc{F}^{\pd}: \widetilde{\mc{W}}^{\pd}(\ov{X} \times X) \to (\mc{W}(X), \mc{W}(X))-\op{mod}$ which satisfies the following properties: 
\begin{itemize}
\item the image of the diagonal $\Delta \sub \ov{X} \tms X$ is the diagonal bimodule;
\item the image of any product Lagrangian $K \tms L$ is the Yoneda bimodule $\mc{Y}^l_K \otimes_k \mc{Y}^r_L$.
\end{itemize}
\eprop
\pf
Recall that the wrapped Fukaya category $\mc{W}(X)$ is defined as a localization $\mc{O}_X [C^{-1}]$. Here $\mc{O}_{X}$ is an $A_{\infty}$ category whose objects are given by $\mathbb{Z}_{\geq 0} \times I_{X}$, where $I_X$ is a countable collection of cylindrical Lagrangians in $X$ representing all isotopy classes and for each Lagrangian $L \in I_X$ we have a cofinal wrapping sequence $L = L^{(0)} \leadsto L^{(1)} \leadsto \cdots$. The set of morphisms $C$ consists of continuation elements in $HF^{0}(L^{(i+1)}, L^{(i)})$. For details, see \cite[Sec.\ 5.5]{gps2}.

Without loss of generality, we can assume that the product Lagrangians in the set $\tilde{I}^{\pd}_{\ov{X} \times X}$ in $\ov{X} \times X$ are all given by product of Lagrangians from $I_X$. Then we can define an $A_{\infty}$ functor
\eq \widetilde{\mc{O}}^{\pd}_{\ov{X} \times X} \rightarrow (\mc{O}_X, \mc{O}_X)-\op{mod} \eeq
by considering quilted holomorphic maps as in \cite[Sec.\ 9]{ganatra}. More precisely, for any Lagrangian $\mc{L} \in \widetilde{\mc{O}}^{\pd}_{\ov{X} \times X}$, it defines an $A_{\infty}$ bimodule $\mc{M}_{\mc{L}}$ over $\mc{O}_X$ whose morphism space associates to a pair of Lagrangians $K, L \in \mc{O}_X$ the Floer cochain complex
\eq \op{hom}_{\widetilde{\mc{O}}^{\pd}_{\ov{X} \times X}}^{\bu}(K \times L, \mc{L}). \eeq
The $A_{\infty}$ operations are defined by counting quilted holomorphic maps into $\ov{X} \times X$ with seam along the Lagrangian $\mc{L}$, see \cite[Figure 19]{gps2} without inserting marked point along the seam on the left picture and replacing $X^{-}, Y$ both by $X$. The almost complex structures are of product type and are chosen generically. Once again, the compactness arguments from \cite[Prop.\ 3.19]{gps1} and \cite[Sec.\ 8.2]{gps2} guarantee that the relevant moduli spaces which define $\mc{M}_{\mc{L}}$ and verify the $A_{\infty}$ relations are compact. To complete the definition of the $A_{\infty}$ functor $\widetilde{\mc{O}}^{\pd}_{\ov{X} \times X} \rightarrow (\mc{O}_X, \mc{O}_X)-\op{mod}$, it suffices to consider quilted holomorphic maps with marked points inserted along the seam, which are exactly the left picture of \cite[Figure 19]{gps2} after replacing both $X^{-}$ and $Y$ by $X$. Note that we allow Lagrangians $\Delta^{(i)}$ as markings on the seam, but they do not impose further difficulties in the construction. By definition, this functor takes the (strict) diagonal Lagrangian $\Delta$ to the diagonal bimodule $\Delta_{\mc{O}_X}$ and takes the product Lagrangian $K \times L$ to the Yoneda bimodule $\mc{Y}^l_K \otimes_k \mc{Y}^r_L$ over the category $\mc{O}_X$.

The localization functor $\mc{O}_X \rightarrow \mc{W}(X)$ induces an $A_{\infty}$ functor
\eq (\mc{O}_X, \mc{O}_X)-\op{mod} \to (\mc{W}(X), \mc{W}(X))-\op{mod} \eeq
and its composition with the functor defined in the previous paragraph defines an $A_{\infty}$ functor
\eq \widetilde{\mc{O}}^{\pd}_{\ov{X} \times X} \rightarrow (\mc{W}(X), \mc{W}(X))-\op{mod}. \eeq
Using the definition of continuation elements in $C_{\pd}$, we see that the above functor descends to an $A_{\infty}$ functor
\eq \mc{F}^{\pd}: \widetilde{\mc{W}}^{\pd}(\ov{X} \times X) \to (\mc{W}(X), \mc{W}(X))-\op{mod} \eeq
by the universal property of localization. The statements are proved because of \eqref{eqn_diag_floer} and \eqref{eqn_prod_floer}. 
\epf

\lem\label{lem_cyl}
There exists a fully faithful ``cylindrization" $A_{\infty}$ functor $\widetilde{\mc{W}}^{\pd}(\ov{X} \times X) \to \mc{W}(\ov{X} \tms X)$ which takes the diagonal to the diagonal and takes products $K \tms L$ to their cylindrization $K \tilde{\tms} L$. In particular, this functor is a quasi-equivalence onto $\mc{W}_{\op{cyl}}(\ov{X} \tms X)$. 
\elem
\pf
Let $\mc{O}_{\ov{X} \times X}$ be the pre-localized category whose quotient defines $\mc{W}(\ov{X} \tms X)$. The $A_{\infty}$ functor
\eq \widetilde{\mc{O}}^{\pd}_{\ov{X} \times X} \to \mc{O}_{\ov{X} \times X} \eeq
is defined using the extension of the $(\mc{O}^{\pd}_{\ov{X} \times X},  \mc{O}_{\ov{X} \times X})$--bimodule $\mc{B}$ from \cite[Sec.\ 8.2]{gps2} by adding the diagonal Lagrangian. Note that by our choice of wrappings of the diagonal, all $\Delta^{(i)}$ for $i \geq 0$ are cylindrical and obviously the diagonal is taken to the diagonal. The comparisons between Floer cochain complexes as in \cite[Sec.\ 8.2]{gps2} carry over without change so the induced functor from $\mc{B}$ is an equivalence onto its image. The desired functor is constructed by localization and it is a quasi-equivalence onto $\mc{W}_{\op{cyl}}(\ov{X} \tms X)$ by definition.
\epf

\pf[Proof of Theorem \ref{theorem:bimod-functor}]
The functor $\mc{F}$ is constructed as the composition of the inverse of the embedding from Lemma \ref{lem_cyl} and the functor $\mc{F}^{\pd}$. It satisfies the desired properties by Proposition \ref{prop_prod_mod} and Lemma \ref{lem_cyl}. 
\epf

\end{document}